\documentclass[twoside]{book}

\def\thetitle{An invitation to Alexandrov geometry: CAT(0) spaces}
\def\theauthors{Stephanie Alexander, Vitali Kapovitch, Anton Petrunin}

\usepackage{lectures}
\usepackage[colorlinks=true,
citecolor=black,
linkcolor=black,
anchorcolor=black,
filecolor=black,
menucolor=black,
urlcolor=black,
pdftitle={\thetitle},
pdfsubject={Geometry},
pdfauthor={\theauthors},
unicode,
linktocpage
]{hyperref}
\makeindex

\begin{document}

\frontmatter
\title{An invitation to Alexandrov geometry: CAT(0) spaces}
\author{
Stephanie Alexander,\\
Vitali Kapovitch,\\ 
Anton Petrunin
}
\date{}
\maketitle
\thispagestyle{empty}
\mainmatter
\newpage
\cleardoublepage
\phantomsection
\pdfbookmark[0]{\contentsname}{bm:toc}
\tableofcontents

\vfill

%%%%%%%%%%%%%%%%%%%%%%%%%%%%
%%!TEX root = invitation-CAT.tex
\chapter*{Preface}
\addcontentsline{toc}{chapter}{Preface}
\thispagestyle{myheadings}
\markboth{PREFACE}{PREFACE}

This is a substantially updated edition of our earlier book~\cite{alexander-kapovitch-petrunin-2019}, which itself grew out of our more technical text \cite{alexander-kapovitch-petrunin-2025}.
The changes were made based on our teaching experience in recent years and the many comments we received.

We discuss $\CAT(0)$ spaces --- metric spaces with non-positive curvature in the sense of Alexandrov; these spaces can be loosely described as a non-linear generalization of a Hilbert space.
Our goal is to show the beauty and power of Alexandrov geometry by reaching interesting applications and theorems with a minimum of preparation.

\section*{Notes on the contents}

Exercises that are used in the text are marked with an exclamation point: \textbf{Exercise!}

\begin{wrapfigure}{r}{45mm}
\vskip-8mm
\centering
\begin{tikzpicture}[->,>=stealth',shorten >=1pt,auto,scale=1.4,
  thick,main node/.style={circle,draw,font=\sffamily\bfseries,minimum size=8mm}]

  \node[main node] (1) at (0,1) {\ref{prelim}};
  \node[main node] (2) at (0,0) {\ref{chapter:gluing}};
  \node[main node] (3) at (-{sin(360/7)},{cos(360/7)}) {\ref{chap:billiards}};
  \node[main node] (4) at (-{sin(2*360/7)},{cos(2*360/7)}) {\ref{sec:resh-kirz}};
  \node[main node] (5) at (-{sin(3*360/7)},{cos(3*360/7)}) {\ref{chapter:globalization}};
  \node[main node] (6) at (-{sin(4*360/7)},{cos(4*360/7)}) {\ref{chap:poly}};
  \node[main node] (7) at (-{sin(5*360/7)},{cos(5*360/7)}) {\ref{chapter:shefel}};
  \node[main node] (8) at (-{sin(6*360/7)},{cos(6*360/7)}) {\ref{chap:barycenter}};

  \path[every node/.style={font=\sffamily\small}]
   (1) edge[dashed] node{}(2)
   (2) edge node{}(3)
   (2) edge node{}(4)
   (2) edge node{}(5)
   (5) edge node{}(6)
   (6) edge node{}(7)
   (2) edge node{}(8);
\end{tikzpicture}
\end{wrapfigure}

\medskip

In Lecture~\ref{prelim}, we discuss necessary preliminaries.
It should be used as a quick reference.

In Lecture~\ref{chapter:gluing}, we discuss Reshetnyak's gluing theorem.
Further, we apply it in Lecture~\ref{chap:billiards} to a problem in billiards which was solved by Dmitri Burago, Serge Ferleger, and Alexey Kononenko.

Lecture~\ref{sec:resh-kirz} is devoted to Reshetnyak's majorization theorem.
It is illustrated by several applications concerning convexity and geodesics.

In Lecture~\ref{chapter:globalization}, we discuss the Hadamard--Cartan globalization theorem.
Then we apply it in Lecture~\ref{chap:poly} both to describe polyhedral $\CAT(0)$ spaces and to construct exotic aspherical manifolds introduced by Michael Davis.

In Lecture~\ref{chapter:shefel}, we discuss examples of Alexandrov spaces with curvature bounded above.
It is largely based on the work of Samuel Shefel on nonsmooth saddle surfaces.

Finally, in Lecture~\ref{chap:barycenter} we discuss barycenters and their applications to dimension theory; this part is based on the work of Bruce Kleiner.

\medskip

Here is a list of some sources providing a good introduction to Alexandrov spaces with curvature bounded above, which we recommend for further information;
we will not assume familiarity with any of these sources.

\begin{itemize}
\item The book by Martin Bridson and Andr\'e Haefliger \cite{bridson-haefliger} gives the most comprehensive introduction available today;
\item The lecture notes of Werner Ballmann \cite{ballmann-1995,ballmann-2004} include a brief
and clear
introduction;
\item Chapter 9 in the book \cite{burago-burago-ivanov}
 gives another reader-friendly introduction, by Yuriy Burago, Dmitry Burago, and Sergei Ivanov;
\item The book by Jürgen Jost \cite{jost:book} gives a more analytic viewpoint on the subject;
\item Our book \cite{alexander-kapovitch-petrunin-2025} gives a comprehensive treatment of both curvature bounds in the sense of Alexandrov.
\end{itemize}

\section*{Early history of Alexandrov geometry}

The idea that the essence of curvature lies in a condition on quadruples of points apparently originated with Abraham Wald.
It is found in his publication on ``coordinate-free differential geometry'' \cite{wald} written under the supervision of Karl Menger.
(This paper remained almost unnoticed until the 1980s.)
In 1941, similar definitions were rediscovered independently by
Alexandr Danilovich Alexandrov %Alexandr is the right spelling
\cite{alexandrov-1941-def}.
In Alexandrov's work, the first fruitful applications of this approach were given.
Mainly:
\begin{itemize}
\item Alexandrov's embedding theorem --- 
\textit{metrics of non-negative curvature on the sphere, and only they, are isometric to closed convex surfaces in Euclidean 3-space}. 
\item Gluing theorem, which tells when the sphere obtained by gluing two discs along their boundaries has non-negative curvature in the sense of Alexandrov.
\end{itemize}
These two results together gave  a very intuitive geometric tool for studying  embeddings and bendings of surfaces in  Euclidean space, and changed this subject dramatically.
They formed the foundation of the branch of geometry now called \emph{Alexandrov geometry}.

The study of  spaces with curvature bounded above started later.
The first paper on the subject was written by Alexandrov \cite{alexandrov-1951}; it appeared in 1951.
It was based on the work of Herbert Busemann, who studied spaces satisfying a weaker condition \cite{busemann-CBA}.

Yuri Grigorievich Reshetnyak proved fundamental results about general spaces with curvature bounded above, the most important of which are  his gluing and majorization theorems.
An equally important theorem is the Hadamard--Cartan theorem (globalization theorem).
All these theorems and their history are discussed in our lectures.

Surfaces with upper curvature bounds were studied extensively in the 1950s and 1960s, and are by now well understood; see the survey \cite{reshetnyak:survey} and the references therein.

\section*{Manifesto of Alexandrov geometry}

Alexandrov geometry can use ``back to Euclid'' as a slogan.
Alexandrov spaces are defined via axioms similar to those given by Euclid,
but certain  equalities are changed to inequalities. 
Depending on the sign of the inequalities, we get Alexandrov spaces with \emph{curvature bounded above} or \emph{curvature bounded below}.
The definitions of the two classes of spaces are similar, but their properties and known applications are quite different.

Consider the space $\mathcal{M}_4$ of all isometry classes of 4-point metric spaces.
Each element in $\mathcal{M}_4$ can be described by 6 numbers 
 --- the distances between all 6 pairs of its points, say $\ell_{i,j}$ for $1\le i< j\le 4$, modulo permutations of the set of indices $(1,2,3,4)$.
These 6 numbers are subject to 12 triangle inequalities; that is,
\[\ell_{i,j}+\ell_{j,k}\ge \ell_{i,k}\]
holds for all $i$, $j$ and $k$, where we assume that $\ell_{j,i}=\ell_{i,j}$, and $\ell_{i,i}=0$.

\begin{wrapfigure}{o}{33mm}
\vskip-3mm
\centering
\includegraphics{mppics/pic-700}
\end{wrapfigure}

Consider the subset $\mathcal{E}_4\subset \mathcal{M}_4$ of all isometry classes of 4-point metric spaces that admit isometric embeddings into Euclidean space.
The complement $\mathcal{M}_4\setminus \mathcal{E}_4$ has two connected components.

\begin{thm}{Exercise}\label{ex:two-components-of-M4}
Prove the latter statement.
\end{thm}

One of the components will be denoted by $\mathcal{P}_4$ and the other by~$\mathcal{N}_4$.
Here $\mathcal{P}$ and $\mathcal{N}$ stand for {}\emph{positive} 
and {}\emph{negative curvature}, because spheres have no quadruples of type $\mathcal{N}_4$ and
hyperbolic spaces
have no quadruples of type~$\mathcal{P}_4$.

A metric space, with length metric, 
that has no quadruples of points of type $\mathcal{P}_4$ (respectively $\mathcal{N}_4$)
is called an Alexandrov space with non-positive (respectively non-negative) curvature.

Here is an exercise, solving which would force the reader to rebuild a considerable part of Alexandrov geometry.
It might be helpful to spend some time thinking about this exercise before proceeding.
(The length metric is defined in Section~\ref{sec:intrinsic}.)

\begin{thm}{Advanced exercise}\label{ex:convex-set}
Assume $\spc{X}$ is a complete length space, containing only quadruples of type~$\mathcal{E}_4$.
Show that $\spc{X}$ is isometric to a convex set in a Hilbert space.
\end{thm}

In the definition above, 
instead of  Euclidean space 
one can take  
hyperbolic space of curvature~$-1$.
In this case,
one obtains the definition of spaces with curvature bounded above or below by~$-1$.

To define spaces with curvature bounded above or below by $1$,
one has to take the unit 3-sphere 
and specify that only the quadruples of points such that each of the four triangles has perimeter 
less than $2\cdot\pi$ are checked.
The latter condition could be considered as a part of the {}\emph{spherical triangle inequality}.

{\sloppy

\section*{Acknowledgment}
We want to  thank 
David Berg,
Richard Bishop,
Yuri Burago,
Maxime Fortier Bourque,
Sergei Ivanov,
Michael Kapovich, 
Bernd Kirchheim, 
Bruce Kleiner,
Nikolai Kosovsky,
Greg Kuperberg,
Nina Lebedeva,
John Lott,
Alexander Lytchak,
Ricardo Mendes,
Dmitri Panov,
Kasra Rafi,
Stephan Stadler,
Wilderich Tuschmann,
and 
Sergio Zamora Barrera
for a number of discussions and suggestions.

}

We thank the mathematical institutions where we worked on this material, including
BIRS, 
EIMI,
MFO, 
the Henri Poincar\'{e} Institute,
the University of Cologne, and
the Max Planck Institute for Mathematics.

The first author was partially supported by  the Simons Foundation grant \#209053.
The second author was partially supported by a Discovery grant from NSERC and by the Simons Foundation grant \#390117.
The third author was partially supported by the
NSF grants DMS 1309340 and DMS 2005279
and the Simons Foundation grant \#584781.

%%!TEX root = the-preliminaries.tex
\chapter{Preliminaries}\label{prelim}

In this lecture we fix some conventions and recall the main definitions.
It may be used as a quick reference when reading the book.

To learn background in metric geometry, the reader  may consult the book of Dmitri Burago, Yuri Burago, and Sergei Ivanov~\cite{burago-burago-ivanov}
or the book by the third author \cite{petrunin-2023}.

\section{Metric spaces}
\label{sec:metric spaces}

The distance between two points $x$ and $y$ in a metric space $\spc{X}$ will be denoted by $\dist{x}{y}{}$ or $\dist{x}{y}{\spc{X}}$.
The latter notation is used if we need to emphasize 
that the distance is taken in the space~${\spc{X}}$.

The function 
\[\distfun_x\:y\mapsto \dist{x}{y}{}\]
is called the \index{distance function}\emph{distance function} from~$x$. 

\begin{itemize}
\item The \index{diameter}\emph{diameter} of a metric space $\spc{X}$ is defined as
\[\diam \spc{X}=\sup\set{\dist{x}{y}{\spc{X}}}{x,y\in \spc{X}}.\]

\item Given $R\in[0,\infty]$ and $x\in \spc{X}$, the sets
\begin{align*}
\oBall(x,R)&=\{y\in \spc{X}\mid \dist{x}{y}{}<R\},
\\
\cBall[x,R]&=\{y\in \spc{X}\mid \dist{x}{y}{}\le R\}
\end{align*}
are called, respectively, the  \index{open ball}\emph{open} and  the \index{closed ball}\emph{closed  balls}   of radius $R$ with center~$x$.
Again, if we need to emphasize that these balls are taken in the metric space $\spc{X}$,
we write 
\[\oBall(x,R)_{\spc{X}}\quad\text{and}\quad\cBall[x,R]_{\spc{X}}.\]
\end{itemize}

\section{Geodesics, triangles, and hinges}
\label{sec:geods}

\parbf{Geodesic.}
Let $\spc{X}$ be a metric space 
and $\II$\index{3@$\II$ (real interval)} be a real interval.
A~globally isometric map $\gamma\:\II\z\to \spc{X}$ is called a \index{geodesic}\emph{geodesic}%
\footnote{Various authors call it differently: {}\emph{shortest path}, {}\emph{minimizing geodesic}.}; 
in other words, $\gamma\:\II\z\to \spc{X}$ is a geodesic if
\[\dist{\gamma(s)}{\gamma(t)}{\spc{X}}=|s-t|\]
for any pair $s,t\in \II$.

We say that  $\gamma\:\II\to \spc{X}$ is a geodesic from point $p$ to point $q$ if 
$\II=[a,b]$ and $p=\gamma(a)$, $q=\gamma(b)$.
In this case the image of $\gamma$ is denoted by $[p q]$\index{2@$[x y]=[x y]_{\spc{X}}$ (geodesic)} and with an abuse of notation  we also call it a \index{geodesic}\emph{geodesic}.

We may write $[p q]_{\spc{X}}$ 
to emphasize that the geodesic $[p q]$ is in the space  ${\spc{X}}$.
We also use the following shortcut notation:
\begin{align*}
\left] p q \right[&=[pq]\setminus\{p,q\},
&
\left] p q \right]&=[pq]\setminus\{p\},
&
\left[ p q \right[&=[pq]\setminus\{q\}.
\end{align*}

In general, a geodesic between $p$ and $q$ need not exist and if it exists, it need not  be unique.  
However,  once we write $[p q]$ we mean that we have  made a choice of  a geodesic.

A metric space is called \index{geodesic}\emph{geodesic} if any pair of its points can be joined by a geodesic.

A \index{geodesic!path}\emph{geodesic path} is a geodesic with constant-speed parametrization by $[0,1]$.
%Given a geodesic $[p q]$, we denote by $\geodpath_{[pq]}$ the corresponding geodesic path; that is, $$\geodpath_{[pq]}(t)\z\df\geod_{[pq]}(t\cdot\dist[{{}}]{p}{q}{}).$$

A curve $\gamma\:\II\to \spc{X}$  is called a \index{geodesic!local geodesic}\emph{local geodesic} if for any $t\in\II$ there is a neighborhood $U$ of $t$ in $\II$ such that the restriction $\gamma|_U$ is a  geodesic.
A constant-speed parametrization of a local geodesic by the unit interval $[0,1]$ is called a \index{geodesic!local geodesic}\emph{local geodesic path}.

\parbf{Triangle.}
For a triple of points $p,q,r\in \spc{X}$, a choice of a triple of geodesics $([q r], [r p], [p q])$ will be called a \index{triangle}\emph{triangle}; we will use the short notation 
$\trig p q r=\trig p q r_{\spc{X}}=([q r], [r p], [p q])$\index{2@$\trig xyz=\trig xyz_{\spc{X}}$ (triangle)}.

Again, given a triple $p,q,r\in \spc{X}$ there may be no triangle 
$\trig p q r$ simply because one of the pairs of these points cannot be joined by a geodesic.
Also, many different triangles with these vertices may exist, any of which can be denoted by $\trig p q r$.
However, if we write $\trig p q r$, it means that we have made a choice of such a triangle; 
that is, we have  fixed a choice of the geodesics $[q r]$, $[r p]$, and $[p q]$.

The value 
\[\dist{p}{q}{}+\dist{q}{r}{}+\dist{r}{p}{}\]
will be called the {}\emph{perimeter of the triangle} $\trig p q r$.

\parbf{Hinge.}
Let $p,x,y\in \spc{X}$ be a triple of points such that $p$ is distinct from $x$ and~$y$.
A pair of geodesics $([p x],[p y])$ will be called  a \index{hinge}\emph{hinge} and will be denoted by 
$\hinge p x y=\hinge p x y_{\spc{X}}$\index{2@$\hinge yxz=\hinge yxz_{\spc{X}}$ (hinge)}.

\parbf{Convex set.}\label{def:convex-set}
A set $A$ in a metric space $\spc{X}$ is called 
\index{convex set}\emph{convex}
if for every two points $p,q\in A$, 
{}\emph{every} geodesic $[pq]$ in $\spc{X}$ 
lies in~$A$.

A set $A\subset\spc{X}$ is called 
\index{locally convex}\emph{locally convex}
if every point $a\in A$ admits an open neighborhood $\Omega\ni a$ in $\spc{X}$ 
such that any geodesic lying in $\Omega$ and with ends in $A$ lies completely in~$A$.

Note that any open set is locally convex by  definition.

\section{Length spaces}\label{sec:intrinsic}

A \emph{curve} is defined as a continuous map from a real interval to a space.
If the real interval is $[0,1]$, then the curve is called a \emph{path}.

\begin{thm}{Definition}
Let $\spc{X}$ be a metric space and
$\alpha\: \II\to \spc{X}$ be a curve.
We define the \index{length}\emph{length} of $\alpha$ as 
\[
\length \alpha \df \sup_{t_0\le t_1\le\ldots\le t_n}\sum_i \dist{\alpha(t_i)}{\alpha(t_{i-1})}{}.
\]
\end{thm}

Directly from the definition, it follows that if a path $\alpha\:[0,1]\to\spc{X}$ connects two points $x$ and $y$ 
(that is, if $\alpha(0)=x$ and $\alpha(1)=y$), then 
\[\length\alpha\ge \dist{x}{y}{}.\]

Let $A$ be a subset of a metric space $\spc{X}$.
Given two points $x,y\in A$,
consider the value
\[\dist{x}{y}{A}=\inf_{\alpha}\{\length\alpha\},\]
where the infimum is taken for all paths $\alpha$ from $x$ to $y$ in $A$.\footnote{Note that while this notation slightly conflicts with the previously defined notation for distance on a general metric space, we will usually work with ambient length spaces where the meaning will be unambiguous.}

If $\dist{x}{y}{A}$ takes a finite value for each pair $x,y\in A$,
then $\dist{x}{y}{A}$ defines a metric on  $A$;
this metric will be called the \index{induced length metric}\emph{induced length metric} on $A$.

If for any $\eps>0$ and any pair of points $x$ and $y$ in a metric space $\spc{X}$, there is a path $\alpha$ connecting $x$ to $y$ such that
\[\length\alpha< \dist{x}{y}{}+\eps,\]
then $\spc{X}$ is called a \index{length!space}\emph{length space} and the metric on $\spc{X}$ is called a \index{length!metric}\emph{length metric}.\label{page:length metric}

If $f\:\tilde{\spc{X}}\to\spc{X}$ is a covering,
then a length metric on $\spc{X}$ can be lifted to $\tilde{\spc{X}}$
by declaring 
\[\length_{\tilde{\spc{X}}}\gamma=\length_{\spc{X}}(f\circ\gamma)\]
for any curve $\gamma$ in $\tilde{\spc{X}}$.
The space $\tilde{\spc{X}}$ with this metric is called the  \index{metric cover}\emph{metric cover} of $\spc{X}$.

Note that any geodesic space is a length space.
As can be seen from the following example, the converse does not hold.

\begin{thm}{Example}
Let $\spc{X}$ be obtained by gluing a countable collection of disjoint intervals $\{\II_n\}$ of length $1+\tfrac1n$, where for each $\II_n$ the left end is glued to $p$ and the right end to~$q$.
Then $\spc{X}$ carries a natural complete length metric with respect to which $\dist{p}{q}{}=1$ but there is no geodesic connecting $p$ to~$q$.
\end{thm}

\begin{thm}{Exercise}\label{ex:no-geod}
Construct a complete length space for which no pair of distinct points can be joined by a geodesic.
\end{thm}

Let $\spc{X}$ be a metric space and $x,y\in\spc{X}$.

\begin{enumerate}[(i)]
\item A point $z\in \spc{X}$ is called a \index{midpoint}\emph{midpoint} between $x$ and $y$
if 
\[\dist{x}{z}{}=\dist{y}{z}{}=\tfrac12\cdot\dist[{{}}]{x}{y}{}.\]
\item Assume $\eps\ge 0$.
A point $z\in \spc{X}$ is called an {}\emph{$\eps$-midpoint} between $x$ and $y$
if 
\[\dist{x}{z}{},\quad\dist{y}{z}{}\le\tfrac12\cdot\dist[{{}}]{x}{y}{}+\eps.\]
\end{enumerate}

Note that a $0$-midpoint is the same as a midpoint.

\begin{thm}{Menger's lemma}\label{lem:mid>geod}
Let $\spc{X}$ be a complete metric space.
\begin{subthm}{lem:mid>length}
Assume that for any pair of points $x,y\in \spc{X}$  
 and any $\eps>0$
there is an $\eps$-midpoint~$z$.
Then $\spc{X}$ is a length space.
\end{subthm}

\begin{subthm}{lem:mid>geod:geod}
Assume that for any pair of points $x,y\in \spc{X}$, 
there is a midpoint~$z$.
Then $\spc{X}$ is a geodesic space.
\end{subthm}
\end{thm}

The second part of this lemma was proved by Karl Menger \cite[Section 6]{menger}.

\parit{Proof.}
We first prove \ref{SHORT.lem:mid>length}.
Let $x,y\in \spc{X}$ be a pair of points.

Set $\eps_n=\frac\eps{4^n}$, $\alpha(0)=x$ and $\alpha(1)=y$.

Let $\alpha(\tfrac12)$ be an $\eps_1$-midpoint between $\alpha(0)$ and $\alpha(1)$.
Further, let $\alpha(\frac14)$ 
and $\alpha(\frac34)$ be $\eps_2$-midpoints between the pairs $(\alpha(0),\alpha(\tfrac12))$ 
and $(\alpha(\tfrac12),\alpha(1))$ respectively.
Applying the above procedure recursively,
on the $n$-th step we define $\alpha(\tfrac{\kay}{2^n})$,
for every odd integer $\kay$ such that $0<\tfrac\kay{2^n}<1$, 
as an $\eps_{n}$-midpoint between the already defined
$\alpha(\tfrac{\kay-1}{2^n})$ and $\alpha(\tfrac{\kay+1}{2^n})$.

In this way we define $\alpha(t)$ for $t\in W$,
where $W$ denotes the set of dyadic rationals in $[0,1]$.
Since $\spc{X}$ is complete, the map $\alpha$ can be extended continuously to $[0,1]$.
Moreover,
\[\begin{aligned}
\length\alpha&\le \dist{x}{y}{}+\sum_{n=1}^\infty 2^{n-1}\cdot\eps_n\le
\\
&\le \dist{x}{y}{}+\tfrac\eps2.
\end{aligned}
\eqlbl{eq:eps-midpoint}
\]
Since $\eps>0$ is arbitrary, we get \ref{SHORT.lem:mid>length}.

To prove \ref{SHORT.lem:mid>geod:geod}, 
one should repeat the same argument 
taking midpoints instead of $\eps_n$-midpoints.
In this case \ref{eq:eps-midpoint} holds for $\eps_n=\eps=0$.
\qeds

A metric space $\spc{X}$ is called \index{proper space}\emph{proper} if all closed bounded sets in $\spc{X}$ are compact. 
This condition is equivalent to each of the following statements:
\begin{enumerate}
\item For some (and therefore any) point $p\in \spc{X}$ and any $R<\infty$, 
the closed ball $\cBall[p,R]\subset\spc{X}$ is compact. 
\item The function $\distfun_p\:\spc{X}\to\RR$ is proper for some (and therefore any) point $p\in \spc{X}$;
that is, for any compact set $K\subset \RR$, its inverse image 
$\set{x\in \spc{X}}{\dist{p}{x}{\spc{X}}\in K}$
is compact.
\end{enumerate}

Since in a compact space a sequence of $\tfrac1n$-midpoints $z_n$ contains a convergent subsequence, Menger's lemma immediately implies the following.

\begin{thm}{Proposition}\label{prop:length+proper=>geodesic}
A proper length space is geodesic.
\end{thm}

\begin{thm}{Hopf--Rinow theorem}\label{thm:Hopf-Rinow}
Any complete, locally compact length space is proper.
\end{thm}

\parit{Proof.}
Let $\spc{X}$ be a locally compact length space.
Given $x\in \spc{X}$, denote by $\rho(x)$ the supremum of all $R>0$ such that
the closed ball $\cBall[x,R]$ is compact.
Since $\spc{X}$ is locally compact, 
$$\rho(x)>0\ \ \text{for any}\ \ x\in \spc{X}.\eqlbl{eq:rho>0}$$
It is sufficient to show that $\rho(x)=\infty$ for some (and therefore any) point $x\in \spc{X}$.

Assume the contrary; that is, $\rho(x)<\infty$. We claim that

\begin{clm}{} $B=\cBall[x,\rho(x)]$ is compact for any~$x$.
\end{clm}

Indeed, $\spc{X}$ is a length space;
therefore for any $\eps>0$, 
the set $\cBall[x,\rho(x)-\eps]$ is a compact $\eps$-net in~$B$.
Since $B$ is closed and hence complete, it must be compact.
\claimqeds
Next we claim that
\begin{clm}{} $|\rho(x)-\rho(y)|\le \dist{x}{y}{\spc{X}}$ for any $x,y\in \spc{X}$;
in particular $\rho\:\spc{X}\to\RR$ is a continuous function.
\end{clm}

Indeed, 
assume the contrary; that is, $\rho(x)+|x-y|<\rho(y)$ for some $x,y\in \spc{X}$. 
Then 
$\cBall[x,\rho(x)+\eps]$ is a closed subset of $\cBall[y,\rho(y)]$ for some $\eps>0$, and compactness of $\cBall[y,\rho(y)]$ implies compactness of $\cBall[x,\rho(x)+\eps]$ --- a contradiction.
\claimqeds

Set $\eps=\min\set{\rho(y)}{y\in B}$; the minimum is defined since $B$ is compact.
From \ref{eq:rho>0}, we have $\eps>0$.

Choose a finite $\tfrac\eps{10}$-net $\{a_1,\ldots,a_n\}$ in $B$.
The union $W$ of the closed balls $\cBall[a_i,\eps]$ is compact.
Clearly 
$\cBall[x,\rho(x)+\frac\eps{10}]\subset W$, so $\cBall[x,\rho(x)+\frac\eps{10}]$ is compact ---
a contradiction.
\qeds

\begin{thm}{Exercise}\label{exercise from BH}
Construct a geodesic space that is locally compact,
but whose completion is neither geodesic nor locally compact.
\end{thm}

\section{Constructions}\label{sec:constructions}

\parbf{Product space.}
Given two metric spaces $\spc{U}$ and $\spc{V}$, the \index{product space}\emph{product space} 
$\spc{U}\times\spc{V}$ is defined as the set of all pairs $(u,v)$ where $u\in\spc{U}$ and $v\in \spc{V}$ 
with the metric defined by formula
\[\dist{(u_1,v_1)}{(u_2,v_2)}{\spc{U}\times\spc{V}}=\sqrt{\dist[2]{u_1}{u_2}{\spc{U}}+\dist[2]{v_1}{v_2}{\spc{V}}}.\]

\begin{thm}{Exercise}\label{ex:length-prod}
Show that the product of length spaces is a length space.
\end{thm}

\begin{thm}[!]{Exercise}\label{ex:geod-prod}
Show that projections of a geodesic path from the product space to its factors are geodesic paths.
\end{thm}

\parbf{Cone.}
The \index{cone}\emph{cone} $\spc{V}=\Cone\spc{U}$ over a metric space $\spc{U}$
is defined as the metric space whose underlying set consists of the
equivalence classes in
$[0,\infty)\times \spc{U}$ with the equivalence relation ``$\sim$'' given by $(0,p)\sim (0,q)$ for any points $p,q\in\spc{U}$,
and whose metric is given by the cosine rule
\[
\dist{(p,s)}{(q,t)}{\spc{V}} 
=
\sqrt{s^2+t^2-2\cdot s\cdot t\cdot \cos\alpha},
\]
where $\alpha= \min\{\pi, \dist{p}{q}{\spc{U}}\}$.

The point in the cone $\spc{V}$ formed by the equivalence class of $0\times\spc{U}$ is called the \index{tip of the cone}\emph{tip of the cone} and is denoted by $0$ or $0_{\spc{V}}$.
The distance $\dist{0}{v}{\spc{V}}$ is called the norm of $v$ and is denoted by $|v|$ or $|v|_{\spc{V}}$.

\begin{thm}[!]{Exercise}\label{ex:cone-geod}
Let $[pq]$ be a geodesic in $\Cone\spc{U}$;
assume it does not pass thru the tip.
Show that the projection of $[pq]$ to $\spc{U}$ (after reparametrization) is a geodesic of length less than $\pi$.
\end{thm}

\parbf{Suspension.}
The \index{suspension}\emph{suspension} $\spc{V}=\Susp\spc{U}$ over a metric space $\spc{U}$
is defined as the metric space whose underlying set consists of equivalence classes in
$[0,\pi]\times \spc{U}$ with the equivalence relation ``$\sim$'' given by $(0,p)\z\sim (0,q)$ and $(\pi,p)\z\sim (\pi,q)$ for any points $p,q\in\spc{U}$,
and whose metric is given by the  spherical cosine rule
\[
\cos\dist{(p,s)}{(q,t)}{\Susp\spc{U}} 
=
\cos s\cdot\cos t-\sin s\cdot\sin t\cdot\cos\alpha,
\]
where $\alpha= \min\{\pi, \dist{p}{q}{\spc{U}}\}$.

The points in $\spc{V}$ formed by the equivalence classes of $0\times\spc{U}$ and $\pi\times\spc{U}$ are called  the \emph{north} and the  \index{pole of suspension}\emph{south poles} of the suspension.

\begin{thm}{Exercise}\label{ex:product-cone}
Let $\spc{U}$ be a metric space.
Construct an isometry
\[\Cone[\Susp\spc{U}]\to\RR\times \Cone\spc{U}.\]

\end{thm}

\section{Model angles and triangles}\label{sec:mod-tri/angles}

Let $\spc{X}$ be a metric space and 
$p,q,r\in \spc{X}$. 
Let us define the \index{model triangle}\emph{model triangle} $\trig{\tilde p}{\tilde q}{\tilde r}$ 
(briefly, 
$\trig{\tilde p}{\tilde q}{\tilde r}=\modtrig(p q r)_{\EE^2}$%
\index{1@$\modtrig(xyz)$ (model triangle)}) to be a triangle in the plane $\EE^2$ with the same side lengths; that is,
\[\dist{\tilde p}{\tilde q}{}=\dist{p}{q}{},
\quad\dist{\tilde q}{\tilde r}{}=\dist{q}{r}{},
\quad\dist{\tilde r}{\tilde p}{}=\dist{r}{p}{}.\]

In the same way we can define the \index{hyperbolic model triangle}\emph{hyperbolic} and the \index{spherical model triangles}\emph{spherical model triangles} $\modtrig(p q r)_{\HH^2}$, $\modtrig(p q r)_{\mathbb{S}^2}$
in the hyperbolic plane $\HH^2$ and the unit sphere~$\mathbb{S}^2$.
In the latter case the model triangle is said to be defined if in addition
\[\dist{p}{q}{}+\dist{q}{r}{}+\dist{r}{p}{}< 2\cdot\pi.\]
In this case the model triangle again exists and is unique up to an isometry of~$\mathbb{S}^2$.

If 
$\trig{\tilde p}{\tilde q}{\tilde r}=\modtrig(p q r)_{\EE^2}$ 
and $\dist{p}{q}{},\dist{p}{r}{}>0$, 
the angle measure of 
$\trig{\tilde p}{\tilde q}{\tilde r}$ at $\tilde p$ 
will be called the \index{model angle}\emph{model angle} of the triple $p$, $q$, $r$ and will be denoted by
$\angk p q r_{\EE^2}$%
\index{1@$\angk yxz$ (model angle)}.
In the same way we define $\angk p q r_{\HH^2}$ and $\angk p q r_{\mathbb{S}^2}$;
in the latter case  we assume in addition that the model triangle $\modtrig(p q r)_{\mathbb{S}^2}$ is defined.

We may use the notation $\angk p q r$ if it is evident which of the model spaces $\HH^2$, $\EE^2$ or $\mathbb{S}^2$ is meant.

\begin{wrapfigure}{r}{25mm}
\vskip-0mm
\centering
\includegraphics{mppics/pic-730}
\end{wrapfigure}

\begin{thm}{Alexandrov's lemma}
\index{Alexandrov's lemma}
\label{lem:alex}  
Let $p,x,y,z$ be distinct points in a metric space such that $z\in \left]x y\right[$.
Then 
the following expressions for the Euclidean model angles have the same sign:
\begin{subthm}{lem-alex-difference}
$\angk x p y
-\angk x p z$,
\end{subthm} 

\begin{subthm}{lem-alex-angle}
$\angk z p x
+\angk z p y -\pi$.
\end{subthm}

Moreover,
\[\angk p x y \ge \angk p x z +  \angk p z y,\]
with equality if and only if the expressions in \ref{SHORT.lem-alex-difference} and \ref{SHORT.lem-alex-angle} vanish.

The same holds for the hyperbolic and spherical model angles, 
but in the latter case one has to assume in addition that
\[\dist{p}{z}{}+\dist{p}{y}{}+\dist{x}{y}{}< 2\cdot\pi.\]

\end{thm}

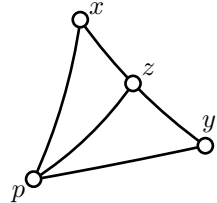
\parit{Proof.} 
Consider the model triangle $\trig{\tilde x}{\tilde p}{\tilde z}=\modtrig(x p z)$.
Take 
a point $\tilde y$ on the extension of 
$[\tilde x \tilde z]$ beyond $\tilde z$ so that $\dist{\tilde x}{\tilde y}{}=\dist{x}{y}{}$ (and therefore $\dist{\tilde x}{\tilde z}{}=\dist{x}{z}{}$). 

\begin{wrapfigure}{r}{33mm}
\vskip-0mm
\centering
\includegraphics{mppics/pic-740}
\end{wrapfigure}

Since increasing the opposite side in a plane triangle increases the corresponding angle, 
the following expressions have the same sign:
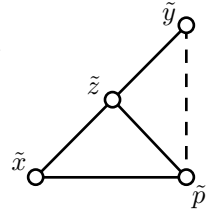
\begin{enumerate}[(i)]
\item $\mangle\hinge{\tilde x}{\tilde p}{\tilde y}-\angk{x}{p}{y}$,
\item $\dist{\tilde p}{\tilde y}{}-\dist{p}{y}{}$,
\item $\mangle\hinge{\tilde z}{\tilde p}{\tilde y}-\angk{z}{p}{y}$.
\end{enumerate}
Since 
\[\mangle\hinge{\tilde x}{\tilde p}{\tilde y}=\mangle\hinge{\tilde x}{\tilde p}{\tilde z}=\angk{x}{p}{z}\]
and
\[ \mangle\hinge{\tilde z}{\tilde p}{\tilde y}
=\pi-\mangle\hinge{\tilde z}{\tilde x}{\tilde p}
=\pi-\angk{z}{x}{p},\]
the first statement follows.

For the second statement, construct a model triangle $\trig{\tilde p}{\tilde z}{\tilde y'}\z=\modtrig(pzy)_{\EE^2}$ on the opposite side of $[\tilde p\tilde z]$ from $\trig{\tilde x}{\tilde p}{\tilde z}$.  
Note that 
\begin{align*}
\dist{\tilde x}{\tilde y'}{}
&\le \dist{\tilde x}{\tilde z}{} + \dist{\tilde z}{\tilde y'}{}=
\\
&=\dist{x}{z}{}+\dist{z}{y}{}=
\\
&=\dist{x}{y}{}.
\intertext{Therefore}
\angk{p}{x}{z} + \angk{p}{z}{y} 
&
= 
\mangle\hinge{\tilde p}{\tilde x}{\tilde z}+ \mangle\hinge{\tilde p}{\tilde z}{\tilde y'} 
=
\\
&
= 
\mangle\hinge{\tilde p}{\tilde x}{\tilde y'}
\le
\\
&\le  \angk p x y.
\end{align*}
Equality holds if and only  if $\dist{\tilde x}{\tilde y'}{}=\dist{x}{y}{}$, 
as required.
\qeds

\section{Angles and the first variation}\label{sec:angles}

Given a hinge $\hinge p x y$, we define its \index{angle}\emph{angle} as 
the limit\index{1@$\mangle\hinge y x z$ (angle)}
\[\mangle\hinge p x y
\df
\lim_{\bar x,\bar y\to p} \angk p{\bar x}{\bar y}_{\EE^2},\eqlbl{eq:def-angle}\]
where $\bar x\in\left]p x\right]$ and $\bar y\in\left]p y\right]$.
(The angle $\mangle\hinge p x y$ is defined if the limit exists.)

The value under the limit can be calculated from the  cosine law:
\[\cos\angk{p}{x}{y}_{\EE^2}
=
\frac{\dist[2]{p}{x}{}+\dist[2]{p}{y}{}-\dist[2]{x}{y}{}}{2\cdot \dist[{{}}]{p}{x}{}\cdot\dist[{{}}]{p}{y}{}}.\]

The following lemma implies that
changing $\angk p{\bar x}{\bar y}_{\EE^2}$ to $\angk p{\bar x}{\bar y}_{\mathbb{S}^2}$ or  $\angk p{\bar x}{\bar y}_{\HH^2}$ in \ref{eq:def-angle} makes no difference.

\begin{thm}{Lemma}\label{lem:k-K-angle}
For any three points $p,x,y$ in a metric space the following inequalities
\[
\begin{aligned}
|\angk p{x}{y}_{\mathbb{S}^2}-\angk p{x}{y}_{\EE^2}|
&\le 
100\cdot(\dist[{{}}]{p}{x}{}+\dist[{{}}]{p}{y}{})^2,
\\
|\angk p{x}{y}_{\HH^2}-\angk p{x}{y}_{\EE^2}|
&\le 
100\cdot(\dist[{{}}]{p}{x}{}+\dist[{{}}]{p}{y}{})^2
\end{aligned}
\eqlbl{eq:k-K}\]
hold whenever the left-hand side is defined.
\end{thm}

\parit{Proof.}
Note that 
\[\angk p{x}{y}_{\HH^2}\le\angk{p}{x}{y}_{\EE^2}\le \angk p{x}{y}_{\mathbb{S}^2}.\]
Therefore
\begin{align*}
0&\le \angk p{x}{y}_{\mathbb{S}^2}-\angk p{x}{y}_{\HH^2}\le
\\
&\le\quad \angk p{x}{y}_{\mathbb{S}^2}+\angk {x}p{y}_{\mathbb{S}^2}+\angk {y}p{x}_{\mathbb{S}^2}
-
\\
&\quad-\angk p{x}{y}_{\HH^2}-\angk {x}p{y}_{\HH^2}-\angk {y}p{x}_{\HH^2}
= 
\\
&=\area\modtrig(pxy)_{\mathbb{S}^2}+\area\modtrig(pxy)_{\HH^2}.
\end{align*}
The inequality \ref{eq:k-K} follows since 
\[
0
\le
\area\modtrig(pxy)_{\HH^2}
\le
\area\modtrig(pxy)_{\mathbb{S}^2}
\le
(\dist[{{}}]{p}{x}{}+\dist[{{}}]{p}{y}{})^2.
\]
\qedsf

\begin{thm}{Exercise}\label{ex:angle-on-shortest}
Let $\spc{X}$ be a complete length space and $p,x,y\in \spc{X}$.
Show that if $p\in \left] x y \right[$, then $\mangle\hinge pxy= \pi$.

\end{thm}

\begin{thm}{Triangle inequality for angles}
\label{claim:angle-3angle-inq}
Let  $[px_1]$, $[px_2]$ and $[px_3]$ be three geodesics in a metric space.
If all  the angles $\alpha_{i j}=\mangle\hinge p {x_i}{x_j}$ are defined, then they satisfy the triangle inequality:
\[\alpha_{13}\le \alpha_{12}+\alpha_{23}.\]

\end{thm}

\parit{Proof.} 
Since $\alpha_{13}\le\pi$, we may assume that $\alpha_{12}+\alpha_{23}< \pi$.

Let $\gamma_i$ be the unit-speed parametrization of $[px_i]$ from $p$ to $x_i$.
Given a small positive $\eps$, for all sufficiently small $t,\tau,s\in\RR_+$ we have
\begin{align*}
\dist{\gamma_1(t)}{\gamma_3(\tau)}{}
\le 
\ &\dist{\gamma_1(t)}{\gamma_2(s)}{}+\dist{\gamma_2(s)}{\gamma_3(\tau)}{}<\\
<
\ &\sqrt{t^2+s^2-2\cdot t\cdot  s\cdot \cos(\alpha_{12}+\eps)}\ +
\\
&+\sqrt{s^2+\tau^2-2\cdot s\cdot \tau\cdot \cos(\alpha_{23}+\eps)}\le
\end{align*}

\begin{wrapfigure}{o}{30 mm}
\vskip-10mm
\centering
\includegraphics{mppics/pic-615}
\vskip3mm
\end{wrapfigure}

Below we define 
$s(t,\tau)$ so that for 
$s=s(t,\tau)$, this chain of inequalities can be continued as follows:
\[\le
\ \sqrt{t^2+\tau^2-2\cdot t\cdot \tau\cdot \cos(\alpha_{12}+\alpha_{23}+2\cdot \eps)}.
\]
Thus for any $\eps>0$, 
\[\alpha_{13}\le \alpha_{12}+\alpha_{23}+2\cdot \eps.\]
Hence the result.

To define $s(t,\tau)$, consider three rays $\tilde\gamma_1$, $\tilde\gamma_2$, $\tilde\gamma_3$ on a Euclidean plane starting at one point, such that $\mangle(\tilde\gamma_1,\tilde\gamma_2)\z=\alpha_{12}+\eps$, $\mangle(\tilde\gamma_2,\tilde\gamma_3)\z=\alpha_{23}+\eps$ and $\mangle(\tilde\gamma_1,\tilde\gamma_3)\z=\alpha_{12}+\alpha_{23}+2\cdot \eps$.
We parametrize each ray by the distance from the starting point.
Given two positive numbers $t,\tau\in\RR_+$, let $s=s(t,\tau)$ be 
the number such that 
$\tilde\gamma_2(s)\in[\tilde\gamma_1(t)\ \tilde\gamma_3(\tau)]$.
Clearly $s\le\max\{t,\tau\}$, so $t,\tau,s$ may be taken sufficiently small.
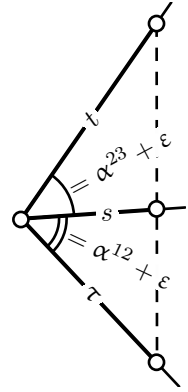
\qeds 

\begin{thm}{Exercise}\label{ex:adjacent-angles}
Prove that the sum of adjacent angles is at least~$\pi$.

More precisely: suppose two hinges $\hinge pxz$ and $\hinge pyz$ are \index{adjacent hinges}\emph{adjacent};
that is, they share side $[pz]$, and the union of two sides $[px]$ and $[py]$ forms a geodesic $[xy]$.
Show that
\[\mangle\hinge pxz+\mangle\hinge pyz\ge \pi\]
whenever  each angle on the left-hand side is defined.

Give an example showing that the inequality can be strict.
\end{thm}

\begin{thm}{First variation inequality}\label{lem:first-var}
Assume that for a  hinge $\hinge q p x$ 
the angle $\alpha=\mangle\hinge q p x$ is defined.
Then
\[\dist{p}{\gamma(t)}{}
\le
\dist{q}{p}{}-t\cdot \cos\alpha+o(t),\]
where $\gamma$ is the unit-speed parametrization of $[qx]$ from $q$ to $x$.
\end{thm}

\parit{Proof.} Take a sufficiently small $\eps>0$.
Denote by $\beta$ the unit-speed parametrization of $[qp]$ from $q$ to $p$.
For all sufficiently small $t>0$, we have 
\begin{align*}
\dist{\beta(t/\eps)}{\gamma(t)}{}
&\le 
\tfrac{t}{\eps}\cdot \sqrt{1+\eps^2 -2\cdot \eps\cdot \cos\alpha}+o(t)\le
\\
&\le \tfrac{t}{\eps} -t\cdot \cos\alpha + t\cdot \eps.
\end{align*}
Applying the triangle inequality, we get 
\begin{align*}
\dist{p}{\gamma(t)}{}
&\le \dist{p}{\beta(t/\eps)}{}+\dist{\beta(t/\eps)}{\gamma(t)}{}
\le 
\\
&\le
\dist{p}{q}{} -t\cdot \cos\alpha + t\cdot \eps
\end{align*}
for any fixed $\eps>0$ and all sufficiently small~$t$.
Hence the result.
\qeds

\section{Space of directions and tangent space}
\label{sec:tangent-space+directions}

Let $\spc{X}$ be a metric space with defined angles for all hinges.
Fix a point $p\in \spc{X}$. 

Consider the set $\mathfrak{S}_p$ 
of all nontrivial geodesics  that start at~$p$.
By \ref{claim:angle-3angle-inq}, the triangle inequality holds for the angle measure $\mangle$ on $\mathfrak{S}_p$,
so
 $(\mathfrak{S}_p,\mangle)$ 
forms a \index{semimetric space}\emph{semimetric space};
that is, $\mangle$ satisfies all the conditions of a metric on $\mathfrak{S}_p$, except that  the angle between distinct geodesics might vanish.

The metric space corresponding to  $(\mathfrak{S}_p,\mangle)$ is called the \index{space of geodesic directions}\emph{space of geodesic directions} at $p$, denoted by $\Sigma'_p$ or $\Sigma'_p\spc{X}$.
Elements of $\Sigma'_p$ are called \index{geodesic!directions}\emph{geodesic directions} at~$p$.
Each geodesic direction is formed by an equivalence class of geodesics in $\mathfrak{S}_p$
for the equivalence relation 
\[[px]\sim[py]\ \ \iff\ \ \mangle\hinge pxy=0.\]

The completion of $\Sigma'_p$ is called the 
\index{space of directions}\emph{space of directions} at $p$ and is denoted by $\Sigma_p$ or $\Sigma_p\spc{X}$.
Elements of $\Sigma_p$ are called \index{direction}\emph{directions} at~$p$.

The Euclidean cone $\Cone\Sigma_p$ over the space of directions $\Sigma_p$ is called the \index{tangent space}\emph{tangent space} at  $p$ and is denoted by $\T_p$ or $\T_p\spc{X}$.

The tangent space $\T_p$ could also be defined directly, without introducing the space of directions.
To do so, consider the set $\mathfrak{T}_p$ of all geodesics with constant-speed parametrizations starting at~$p$. 
Given $\alpha,\beta\in \mathfrak{T}_p$,
set 
\[\dist{\alpha}{\beta}{\mathfrak{T}_p}
=
\lim_{\eps\to0} 
\frac{\dist{\alpha(\eps)}{\beta(\eps)}{\spc{X}}}\eps
\eqlbl{eq:dist-in-T_p}\]
Since the angles in $\spc{X}$ are defined, 
\ref{eq:dist-in-T_p}
defines a semimetric on $\mathfrak{T}_p$.

The corresponding metric space admits a natural isometric identification with the cone $\T'_p=\Cone\Sigma'_p$.
The elements of $\T'_p$ are  equivalence classes for the relation 
\[\alpha\sim\beta\ \ \iff\ \ \dist{\alpha(t)}{\beta(t)}{\spc{X}}=o(t).\]
The completion of $\T'_p$ is therefore  naturally isometric to~$\T_p$.

Elements of $\T_p$ will be called 
\index{tangent vector}\emph{tangent vectors} 
at $p$, regardless of the fact that $\T_p$ is only a metric cone and need not be a vector space.
Elements of $\T'_p$ will be called 
\index{geodesic!tangent vector}\emph{geodesic tangent vectors}
at~$p$.

\section{Hausdorff convergence}

Let $\spc{X}$ be a metric space and $A\subset \spc{X}$.
We will denote by $\distfun_A(x)$ the distance from $A$ to a point $x$ in $\spc{X}$;
that is,
$$\distfun_A(x)\df\inf\set{\dist{a}{x}{\spc{X}}}{a\in A}.$$
We will use the convention that $\distfun_\emptyset(x)=\infty$ for any $x\in \spc{X}$.
%it's not an assumption but our choice of definition

\begin{thm}{Definition of Hausdorff convergence}\label{def:hausdorff-coverge}
Given a sequence of closed sets $(A_n)_{n=1}^\infty$ in a metric space $\spc{X}$, 
a closed set $A_\infty\subset \spc{X}$ is called the \emph{Hausdorff limit} of $(A_n)_{n=1}^\infty$,
briefly $A_n\to A_\infty$, if 
$$\distfun_{A_n}(x)\to\distfun_{A_\infty}(x)\ \ \text{as}\ \ n\to\infty$$
for every $x\in \spc{X}$.

In this case, the sequence of closed sets $(A_n)_{n=1}^\infty$ is said to \index{Hausdorff convergence}\emph{converge in the sense of Hausdorff}.
\end{thm}

\parbf{Examples.}
Let $D_n$ be the disc in the coordinate plane 
with center $(0,n)$ and radius~$n$.
Then $D_n$ converges to the upper half-plane as $n\to\infty$. 

Note that the sequence of one-point sets $\{(0,n)\}$ converges to the empty set.
Indeed, $\distfun_{\{(0,n)\}}(x)\to \infty=\distfun_\emptyset(x)$ for any~$x$.

\begin{thm}{Exercise}\label{ex:hausdorff-conv}
Let $A_n\to A_\infty$ as in Definition \ref{def:hausdorff-coverge}.  

Show that $A_\infty$ is the set of all points $p$ such that   $p_n\to p$ for some sequence of points  $p_n\in A_n$.

Does  the converse hold? That is, suppose $(A_n)_{n=1}^\infty, A_\infty$ are closed sets such that  $A_\infty$  is the set of all points $p$ such that   $p_n\to p$ for some sequence of points  $p_n\in A_n$.
Does this imply that $A_n\to A_\infty$?
\end{thm}

\begin{thm}{First selection theorem}
Let $\spc{X}$ be a proper metric space
and $(A_n)_{n=1}^\infty$ be a sequence of closed sets in~$\spc{X}$.
Then the sequence  $(A_n)_{n=1}^\infty$ has a convergent subsequence in the sense of Hausdorff.
\end{thm}

\parit{Proof.}
Since $\spc{X}$ is proper,
there is a countable dense set $\{x_1,x_2,\ldots\}$ in~$\spc{X}$.
We can assume that the sequence $d_n=\distfun_{A_n}(x_\kay)$ is bounded for each~$\kay$.
Otherwise there is a subsequence such that $\distfun_{A_n}(x)\to \infty$ as $n\to \infty$ for some (and therefore any) $x$.
Therefore, the corresponding subsequence of $A_n$ converges to the empty set.

Therefore, passing to a subsequence of $(A_n)_{n=1}^\infty$,
we can assume that $\distfun_{A_n}(x_\kay)$ converges as $n\to\infty$ for any fixed~$\kay$.

Note that for each $n$, the function $\distfun_{A_n}\:\spc{X}\to\RR$ is 1-Lipschitz and non-negative.
Therefore the sequence $\distfun_{A_n}$ converges pointwise to a 1-Lipschitz non-negative function $f\:\spc{X}\to\RR$.

Set $A_\infty=f^{-1}(0)$. Let us show that 
\[\distfun_{A_\infty}(y)\le f(y)\] 
for any~$y$.
Assume the contrary;
that is, 
$f(y)<R<\distfun_{A_\infty}(y)$
for some $y\in \spc{X}$ and $R>0$.
Then for any sufficiently large $n$, there is a point $z_n\in A_n$ such that
$\dist{y}{z_n}{}\le R$.
Since $\spc{X}$ is proper, we can pass to a partial limit $z_\infty$ of $z_n$ as $n\to\infty$.

It is clear that $f(z_\infty)=0$;
that is, $z_\infty\in A_\infty$. 
(Note that  this implies that $A_\infty\ne \emptyset$.)
On the other hand, 
\begin{align*}
\distfun_{A_\infty}(y)
\le
\dist{z_\infty}{y}{}
\le R
<
\distfun_{A_\infty}(y)
\end{align*}
 --- a contradiction.

Since $f$ is 1-Lipschitz, $\distfun_{A_\infty}(y)\ge f(y)$.
Therefore
\[\distfun_{A_\infty}(y)= f(y)\]
for any $y\in \spc{X}$; hence the result.
\qeds

\section{Gromov--Hausdorff convergence}

\begin{thm}{Definition}\label{def:comp-metr}
Let $\set{\spc{X}_\alpha}{\alpha\in\IndexSet}$ be a collection of metric spaces.
A metric $\rho$ on the disjoint union
$$\bm{X}=\bigsqcup_{\alpha\in\IndexSet} \spc{X}_\alpha$$
is called a  \index{compatible metric}\emph{compatible metric}
if the restriction of $\rho$ to every $\spc{X}_\alpha$ coincides with the original metric on~$\spc{X}_\alpha$.
\end{thm}

\begin{thm}{Definition}\label{def:GH}
Let $\spc{X}_1,\spc{X}_2,\ldots$
and $\spc{X}_\infty$ be proper metric spaces 
and $\rho$ be a compatible metric on their disjoint union~$\bm{X}$.
Assume that $(\bm{X},\rho)$ is proper,
$\spc{X}_n$ is an open set in
$(\bm{X},\rho)$ for each $n\ne\infty$, and 
$\spc{X}_n\to \spc{X}_\infty$ in $(\bm{X},\rho)$ as $n\to\infty$ in the sense of Hausdorff (see Definition~\ref{def:hausdorff-coverge}).

Then we say that $\rho$ defines a {}\emph{convergence}%
\footnote{Formally speaking, convergence in the topology induced by $\rho$ on~$\bm{X}$.} 
in the sense of \emph{Gromov--Hausdorff}%
\index{Gromov--Hausdorff convergence},
and write $\spc{X}_n\to \spc{X}_\infty$ or $\spc{X}_n\xrightarrow{\rho} \spc{X}_\infty$.
The space $\spc{X}_\infty$ is called the \emph{limit space} of the sequence $(\spc{X}_n)$ \emph{along} $\rho$.
\end{thm}

Usually Gromov--Hausdorff convergence is defined differently.
We prefer this definition since it induces convergence for a sequence of points $x_n\in\spc{X}_n$ (Exercise \ref{ex:hausdorff-conv}),
as well as 
weak convergence of measures $\mu_n$ on $\spc{X}_n$, 
and so on,
corresponding to convergence in the ambient space $(\bm{X},\rho)$.

Once we write $\spc{X}_n\to \spc{X}_\infty$, we mean that we have  made a choice of convergence.
Note that for a fixed sequence of metric spaces $(\spc{X}_n)$, it might be possible to construct different Gromov--Hausdorff convergences, say $\spc{X}_n\z{\xrightarrow{\rho}} \spc{X}_\infty$ and $\spc{X}_n\xrightarrow{\rho'} \spc{X}_\infty'$,   whose limit spaces $\spc{X}_\infty$ and $\spc{X}_\infty'$ need not be isometric to each other. 

For example, the constant sequence $\spc{X}_n\iso\RR_{\ge0}$
may converge to $\RR_{\ge0}$ and $\RR$.
The first convergence is evident and the second could be guessed from the diagram.

\begin{figure}[ht!]
\vskip-0mm
\centering
\includegraphics{mppics/pic-500}
\end{figure}

\begin{thm}{Second selection theorem}\label{thm:gromov-selection}
Let $\spc{X}_n$ be a sequence of proper metric spaces 
with marked points $x_n\in \spc{X}_n$.
Assume that for any fixed $R,\eps>0$, there is $N=N(R,\eps)\in\NN$ 
such that for each $n$
the ball $\cBall[x_n,R]_{\spc{X}_n}$ admits a finite $\eps$-net with at most $N$ points.
Then there is a subsequence of $\spc{X}_n$ admitting a Gromov--Hausdorff convergence 
such that the sequence of marked points $x_n\in\spc{X}_n$ converges.
\end{thm}

\parit{Proof.}
From the main assumption in the theorem,
in each space $\spc{X}_n$ 
there is a sequence of points $z_{i,n}\in\spc{X}_n$ such that the following condition holds for a fixed sequence of integers $M_1<M_2<\ldots$
\begin{itemize}
\item $\dist{z_{i,n}}{x_n}{\spc{X}_n}\le \kay+1$ if $i\le M_\kay$,
\item the points $z_{1,n},\ldots,z_{M_\kay,n}$ form an $\tfrac1\kay$-net in $\cBall[x_n,\kay]_{\spc{X}_n}$.
\end{itemize}

Passing to a subsequence, we can assume that the sequence \[\ell_n=\dist{z_{i,n}}{z_{j,n}}{\spc{X}_n}\] 
converges for any $i$ and~$j$.

Consider a countable set of points $\spc{W}=\{w_1,w_2,\ldots\}$
equipped with the semimetric defined by
\[\dist{w_i}{w_j}{\spc{W}}
=
\lim_{n\to\infty}\dist{z_{i,n}}{z_{j,n}}{\spc{X}_n}.\]
Let $\hat{\spc{W}}$ be the metric space corresponding to $\spc{W}$;
that is, points in $\hat{\spc{W}}$ are equivalence classes in $\spc{W}$
for the relation $\sim$, where $w_i\sim w_j$ if and only if $\dist{w_i}{w_j}{\spc{W}}=0$, 
and where 
\[\dist{[w_i]}{[w_j]}{\hat{\spc{W}}}\df\dist{w_i}{w_j}{\spc{W}}.\]
Denote by
$\spc{X}_\infty$ the completion of~$\hat{\spc{W}}$.

It remains to show that there is a Gromov--Hausdorff convergence 
$\spc{X}_n\to\spc{X}_\infty$ such that the sequence $x_n\in\spc{X}_n$ converges.
To prove this, we need to construct a metric $\rho$ on the disjoint union of \[\bm{X}=\spc{X}_\infty\sqcup\spc{X}_1\sqcup\spc{X}_2\sqcup\ldots\]  satisfying definitions \ref{def:comp-metr} and~\ref{def:GH}.
The metric $\rho$ can be constructed as the maximal compatible metric
such that 
\[\rho(z_{i,n},w_i)\le\tfrac1m\]
for any $n\ge N_m$ and $i<I_m$ for a suitable choice of two sequences 
$(I_m)$ and $(N_m)$ with $I_1=N_1=1$.
\qeds

\begin{thm}{Exercise}\label{ex:non-contracting-map}
Let $\spc{K}$  be a compact metric space and
\[f\:\spc{K}\z\to \spc{K}\] 
be a distance non-decreasing map.
Prove that $f$ is an isometry.
\end{thm}

\begin{thm}{Exercise}\label{ex:compact-proper-GH}
Let $\spc{X}_n$ be a sequence of metric spaces that admits 
two convergences $\spc{X}_n\xrightarrow{\rho}\spc{X}_\infty$ and $\spc{X}_n\xrightarrow{\rho'}\spc{X}_\infty'$ to nonempty spaces.
\begin{subthm}{ex:compact-proper-GH:a}
If $\spc{X}_\infty$ is compact, then $\spc{X}_\infty\iso\spc{X}_\infty'$.
\end{subthm}

\begin{subthm}{ex:compact-proper-GH:b}
If  $\spc{X}_\infty$ is proper and there is a sequence of points $x_n\in \spc{X}_n$ 
that  converges in both convergences, 
 then $\spc{X}_\infty\iso\spc{X}_\infty'$.
\end{subthm}
\end{thm}

\section{Notes}

Hausdorff convergence was essentially defined by Dimitrie Pompeiu in his thesis \cite[§ 22]{pompeiu1905};
it was rediscovered by Felix Hausdorff~\cite[VIII~§~6]{hausdorff} and a couple of years later by Wilhelm Blaschke in~\cite{blaschke}.
A refinement of this definition was introduced by Zden\v{e}k Frol\'{\i}k~\cite{frolik},
and later rediscovered by Robert Wijsman~\cite{wijsman}.
This refinement is intermediate between the original Hausdorff convergence and closed convergence, also introduced by Hausdorff in \cite{hausdorff},
and we call it Hausdorff convergence
instead of
Pompeiu--Hausdorff--Blaschke--Frol\'{\i}k--Wijsman convergence.

%%!TEX root = the-billiard.tex
\chapter{Gluing}\label{chapter:gluing}

In this lecture we define $\CAT(\kappa)$ spaces and prove Reshetnyak's gluing theorem.

Here ``$\CAT$'' is an acronym for Cartan, Alexandrov, and Toponogov.
It was coined by Mikhael Gromov in 1987.
Originally, Alexandrov called these spaces ``$\mathfrak{R}_\kappa$ domain'';
this term is still in use.

\section{The 4-point condition}

Given a quadruple of points $p,q,x,y$ in a metric space $\spc{X}$,
consider two model triangles in the plane 
$\trig{\tilde p}{\tilde x}{\tilde y}=\modtrig{}(pxy)_{\EE^2}$ 
and 
$\trig{\tilde q}{\tilde x}{\tilde y}=\modtrig{}(qxy)_{\EE^2}$ with common side $[\tilde x\tilde y]$.

\begin{wrapfigure}{r}{25mm}
\vskip-4mm
\centering
\includegraphics{mppics/pic-720}
\end{wrapfigure}

If the inequality
\[\dist{p}{q}{\spc{X}}\le \dist{\tilde p}{\tilde z}{\EE^2}+\dist{\tilde z}{\tilde q}{\EE^2}\]
holds for any point $\tilde z\in [\tilde x\tilde y]$, then we say that 
the quadruple $p,q,x,y$ satisfies \index{CAT(0) comparison@$\CAT(0)$ comparison}\emph{$\CAT(0)$ comparison}.
\label{page:CAT-comparison}

If we do the same for spherical model triangles  
$\trig{\tilde p}{\tilde x}{\tilde y}=\modtrig{}(pxy)_{\mathbb{S}^2}$ 
and 
$\trig{\tilde q}{\tilde x}{\tilde y}=\modtrig{}(qxy)_{\mathbb{S}^2}$,
then we arrive at the definition of $\CAT(1)$ comparison.
If one of the spherical model triangles is undefined\footnote{That is, if
\[\dist{p}{x}{}+\dist{p}{y}{}+\dist{x}{y}{}\ge 2\cdot\pi
\quad
\text{or}
\quad
\dist{q}{x}{}+\dist{q}{y}{}+\dist{x}{y}{}\ge 2\cdot\pi.\]},
then it is assumed that $\CAT(1)$ comparison automatically holds for this quadruple.

We can do the same for the model plane of given curvature $\kappa$;
that is, an appropriately rescaled sphere if $\kappa>0$,
the Euclidean plane if $\kappa\z=0$
and a rescaled Lobachevsky plane if $\kappa<0$.
In this way we arrive at the definition of $\CAT(\kappa)$ comparison for any real $\kappa$.
However, we will mostly consider  $\CAT(0)$ comparison and occasionally $\CAT(1)$ comparison;
so, if you see $\CAT(\kappa)$, you can safely assume that $\kappa$ is $0$ or~$1$.

If all quadruples in a metric space $\spc{X}$ satisfy $\CAT(\kappa)$ comparison, then we say that the space $\spc{X}$ is $\CAT(\kappa)$.
(Note that $\CAT(\kappa)$ is an adjective.)

In order to check $\CAT(\kappa)$ comparison for a given quadruple of points, it is sufficient to know the 6 distances between all pairs of points in it.
This observation implies the following.

\begin{thm}{Proposition}\label{prop:cat-limit}
Any Gromov--Hausdorff limit of a sequence of $\CAT(\kappa)$ spaces is $\CAT(\kappa)$. 
\end{thm}

In the proposition above, 
it does not matter which definition of convergence for metric spaces we use,
as long as any quadruple of points in the limit space can be arbitrarily well approximated by  quadruples in the sequence of metric spaces.
In particular, it works for the so-called \emph{ultralimits}; see for example \cite{petrunin-2023}.

\section{Geodesics}

The $\CAT$ comparison can be applied to any metric space,
but it is usually applied to geodesic spaces (or complete length spaces).
To simplify the presentation we will often assume in addition that the space is \index{proper space}\emph{proper}.
The latter means that any closed ball is compact; see Section~\ref{sec:intrinsic}.

\begin{thm}{Proposition}\label{ex:CAT-geodesic}
Let $\spc{U}$ be a complete geodesic $\CAT(0)$ space.
Then any two points in $\spc{U}$ are joined by a unique geodesic.
\end{thm}

\parit{Proof.} 
Suppose there are two geodesics between $x$ and $y$.
Then we can choose two points $p\ne q$ on these geodesics such that $\dist{x}{p}{}=\dist{x}{q}{}$ (and therefore $\dist{y}{p}{}=\dist{y}{q}{}$).

Observe that the model triangles $[\tilde p\tilde x\tilde y]=\modtrig(pxy)$ and $[\tilde q\tilde x\tilde y]\z=\modtrig(qxy)$ are degenerate and moreover $\tilde p=\tilde q$.
Applying $\CAT(0)$ comparison with $\tilde z=\tilde p=\tilde q$,
we get that $\dist{p}{q}{}=0$ --- a contradiction.
\qeds

\begin{wrapfigure}[3]{r}{25mm}
\vskip-8mm
\centering
\includegraphics{mppics/pic-750}
\end{wrapfigure}

\begin{thm}[!]{Exercise}\label{ex:noncreasing-CAT}
Given a hinge $\hinge p x y$ in a $\CAT(0)$ space $\spc{U}$ consider the function 
\[f\:(\dist{p}{\bar x}{},\dist{p}{\bar y}{})\mapsto \angk p{\bar x}{\bar y},\]
where $\bar x\in\left]p x\right]$ and $\bar y\in\left]p y\right]$.
Show that $f$ is nondecreasing in each argument.

Conclude that any hinge in a $\CAT(0)$ space has a well  defined angle.
\end{thm}

\begin{thm}[!]{Exercise}\label{ex:contractible}
Fix a point $p$ in a complete geodesic $\CAT(0)$ space~$\spc{U}$.
Given a point $x\in \spc{U}$, denote by $\gamma_x\:[0,1]\to \spc{U}$ the (necessarily unique) geodesic path from $p$ to $x$.

Show that the family of maps $h_t\: \spc{U}\to \spc{U}$ defined by 
\[h_t(x)= \gamma_x(t)\]
is a homotopy.
Conclude that $\spc{U}$ is contractible.
\end{thm}

The homotopy in the exercise is a special case of the so-called \index{geodesic!homotopy}\emph{geodesic homotopy}.
Namely, given two maps $h_0,h_1\:\spc{X}\to \spc{U}$ we define homotopy
\[h_t(x)= \gamma_x(t),\]
where $\gamma_x$ is the geodesic path from $h_0(x)$ to $h_1(x)$.

This construction will be useful in the following exercise.

\begin{thm}{Exercise}\label{ex:CAT-mnfld=>ext.geod}
Let $\spc{U}$ be a complete geodesic $\CAT(0)$ space.
Assume $\spc{U}$ is a topological manifold.
Show that any geodesic in $\spc{U}$ can be extended 
as a two-sided infinite geodesic.
\end{thm}

\begin{thm}{Exercise}\label{ex:geod-CBA}
Assume $\spc{U}$ is a proper geodesic $\CAT(0)$ space
 with extendable geodesics;
that is,
any geodesic extends to a geodesic $\RR\to \spc{U}$.
Show that the space of geodesic directions at any point in $\spc{U}$ is complete.

What if $\spc{U}$ is only complete?
\end{thm}

\begin{thm}{Advanced exercise}\label{ex:tan-CBA}
Let $\spc{U}$ be a proper length $\CAT(0)$ space.

\begin{subthm}{ex:tan-CBA:geodesic}
Show that $\spc{U}$ is geodesic.
\end{subthm}

\begin{subthm}{ex:tan-CBA:tan}
Show that for any $p\in\spc{U}$ the tangent space $\T_p\,\spc{U}$ is a geodesic $\CAT(0)$ space.
\end{subthm}

\end{thm}

\section{Thin triangles}

Let $\trig{\tilde x_1}{\tilde x_2}{\tilde x_3}\z=\modtrig({x_1}{x_2}{x_3})$
be a model triangle for a triangle 
$\trig{x_1}{x_2}{x_3}$ in a metric space.
The map that sends a point $\tilde z\in[\tilde x_i\tilde x_j]$ to the corresponding point $z\in[x_ix_j]$ for each side will be called \index{natural map}\emph{natural}.

\begin{thm}{Definition}\label{def:k-thin}
A triangle $\trig{x}{y}{z}$ in the metric space $\spc{U}$ 
is called \index{thin triangle}\emph{thin} if the natural map $\modtrig{}({x}{y}{z})_{\EE^2}\to \trig{x}{y}{z}$ is distance nonincreasing.

Analogously, a triangle $\trig{x}{y}{z}$ 
is called \index{spherically thin}\emph{spherically thin} if
the natural map from the spherical model triangle $\modtrig{}({x}{y}{z})_{\mathbb{S}^2}$ to $\trig{x}{y}{z}$ is distance nonincreasing.

\end{thm}

\begin{thm}{Proposition}\label{prop:thin=cat}
A geodesic space is $\CAT(0)$ 
($\CAT(1)$) 
if and only if 
all its triangles are thin (respectively, all its triangles of perimeter $<2\cdot\pi$ are spherically thin).
\end{thm}

\parit{Proof.} The if part follows from the triangle inequality and thinness of triangles $\trig pxy$ and $\trig qxy$, where $p$, $q$, $x$, and $y$ are as in the definition of the $\CAT(\kappa)$ comparison.

\parit{Only-if part.} 
Applying $\CAT(0)$ comparison to a quadruple $p,q,x,y$ with $q\in [xy]$ shows that any triangle satisfies \index{point-side comparison}\emph{point-side comparison}; that is, the distance from a vertex to a  point on the opposite side is no greater than the corresponding distance in the Euclidean model triangle.

Now consider a triangle $\trig{x}{y}{z}$ and let $p\in [xy]$ and $q\in [xz]$.
Let $\tilde p$, $\tilde q$ be the corresponding points on the sides of the model triangle $\modtrig({x}{y}{z})_{\EE^2}$.
Applying \ref{ex:noncreasing-CAT}, we get that
\[\angk {x} {y} {z}_{\EE^2} \ge \angk {x} p q _{\EE^2}.\]
Therefore $ \dist{\tilde p}{\tilde q}{\EE^2}\ge \dist{p}{q}{}$.

The $\CAT(1)$ argument is the same.
\qeds

Recall that a curve $\gamma\:\II\to \spc{U}$ is called a \index{geodesic!local geodesic}\emph{local geodesic} if for any $t\in\II$ there is a neighborhood $U$ of $t$ in $\II$ such that the restriction $\gamma|_U$ is a  geodesic.
Note that a local geodesic has to have unit speed.

\begin{thm}{Proposition}\label{cor:loc-geod-are-min}
Suppose $\spc{U}$ is a proper geodesic $\CAT(0)$ space.  
Then any local geodesic in $\spc{U}$ is a geodesic.

Analogously, if $\spc{U}$ is a proper geodesic $\CAT(1)$ space, then any local geodesic in $\spc{U}$ which is shorter than $\pi$ is a geodesic.
\end{thm}

\parit{Proof.}

We first consider the case when  $\spc{U}$ is $\CAT(0)$. 
Suppose $\gamma\:[0,\ell]\to\spc{U}$ is a local geodesic that is not a geodesic.
Choose $a$ to be the maximal value 
such that $\gamma$ is a geodesic on $[0,a]$. By assumption $a<\ell$.
Further, choose $b>a$ so that $\gamma$ is a geodesic on $[a,b]$.

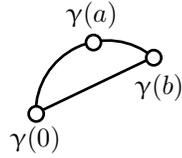
\begin{wrapfigure}{r}{25mm}
\vskip-5mm
\centering
\includegraphics{mppics/pic-800}
\end{wrapfigure}

Since the triangle $\trig{\gamma(0)}{\gamma(a)}{\gamma(b)}$ is thin and 
$\dist{\gamma(0)}{\gamma(b)}{}<b$, we have
\[\dist{\gamma(a-\eps)}{\gamma(a+\eps)}{}<2\cdot\eps\]
for all small~$\eps>0$.
That is, $\gamma$ is not length-minimizing on the interval $[a-\eps,a+\eps]$ for any $\eps>0$ ---
a contradiction.

The spherical case is proved in the same way.
\qeds

\begin{thm}[!]{Exercise}\label{ex:convex-distfun}
Let $\spc{U}$ be a complete geodesic space.
Show that $\spc{U}$ is $\CAT(0)$ if and only if the function $f=\tfrac12\cdot\distfun_p^2$ is {}\emph{1-convex} for any $p\in \spc{U}$;
that is, the function $t\mapsto f\circ \gamma(t) - \tfrac 12\cdot t^2$ is convex for any geodesic $\gamma$.
\end{thm}

\begin{thm}[!]{Exercise}\label{ex:convex-dist}
Suppose $\gamma_1,\gamma_2\:[0,1]\to \spc{U}$ are two geodesic paths in a complete geodesic $\CAT(0)$ space $\spc{U}$.
Show that
\[t\mapsto\dist{\gamma_1(t)}{\gamma_2(t)}{\spc{U}}\]
is a convex function.
\end{thm}

\begin{thm}{Exercise}\label{ex:displacement}
Let $\iota$ be an isometry of a $\CAT(0)$ space $\spc{U}$ to itself.
Show that its \index{displacement function}\emph{displacement function} $f\:\spc{U}\to\RR$ defined by
\[f(x)\df\dist{x}{\iota(x)}{\spc{U}}\]
is \index{convex function}\emph{convex};
that is, the function $t\mapsto f\circ \gamma(t)$ is convex for any geodesic~$\gamma$.
\end{thm}

\begin{thm}{Exercise}\label{ex:convex-nbhd}
Let $A$ be a closed convex set in a geodesic $\CAT(0)$ space $\spc{U}$;
that is, if $x,y\in A$, then $[xy]\subset A$.
Show that $\distfun_A$ is convex.

In particular, for any $r>0$ the set $A$ has a closed convex $r$-neighborhood
\[A_r=\set{x\in \spc{U}}{\distfun_Ax\le r}.\]

\end{thm}

\begin{thm}[!]{Exercise}\label{ex:closest-point}
Let  $\spc{U}$ be a proper geodesic $\CAT(0)$ space 
and $K\subset \spc{U}$ be a closed convex set.
Show that: 

\begin{subthm}{ex:closest-point:a}
For each point $p\in \spc{U}$ there is a unique point $p^*\in K$ that minimizes the distance $\dist{p}{p^*}{}$.
\end{subthm}

\begin{subthm}{ex:closest-point:b}
The nearest-point projection $p\mapsto p^*$ defined by \ref{SHORT.ex:closest-point:a} is short (that is, 1-Lipschitz). 
\end{subthm}

\end{thm}

Recall that a set $A$ in a metric space $\spc{U}$ is called \index{locally convex set}\emph{locally convex} if for any point $p\in A$ there is an open neighborhood $\Omega\ni p$ such that any geodesic in $\Omega$ with  ends in $A$ lies in~$A$.

\begin{thm}[!]{Exercise}\label{ex:locally-convex}
Let $\spc{U}$ be a proper geodesic $\CAT(0)$ space.
Show that any closed, connected, locally convex set in $\spc{U}$ is convex.
\end{thm}

\section{Inheritance lemma} \label{sec:thin-triangle}

\begin{thm}{Inheritance lemma}
\label{lem:inherit-angle} 
Assume that a triangle $\trig p x y$ 
in a metric space is \index{decomposed triangle}\emph{decomposed} 
into two triangles $\trig p x z$ and $\trig p y z$;
that is, $\trig p x z$ and $\trig p y z$ have a common side $[p z]$, and the sides $[x z]$ and $[z y]$ together form the side $[x y]$ of $\trig p x y$.

\begin{wrapfigure}{r}{25mm}
\vskip-0mm
\centering
\includegraphics{mppics/pic-810}
\end{wrapfigure}

If both triangles $\trig p x z$ and $\trig p y z$ are thin, 
then the triangle $\trig p x y$ is also thin.

Analogously, if $\trig p x y$ has perimeter $<2\cdot\pi$ and both triangles $\trig p x z$ and $\trig p y z$ are spherically thin, then triangle $\trig p x y$ is spherically thin.
\end{thm}

\parit{Proof.}
Construct  the model triangles $\trig{\dot p}{\dot x}{\dot z}\z=\modtrig(p x z)_{\EE^2}$ 
and $\trig {\dot p} {\dot y} {\dot z}\z=\modtrig(p y z)_{\EE^2}$ so that $\dot x$ and $\dot y$ lie on opposite sides of $[\dot p\dot z]$.

\begin{wrapfigure}{r}{33mm}
\vskip-0mm
\centering
\includegraphics{mppics/pic-821}
\vskip0mm
\end{wrapfigure}

Let us show that 
\[\angk{z}{p}{x}+\angk{z}{p}{y}
\ge
\pi.
\eqlbl{eq:<+<>=pi}\]
If not, then for some point $\dot w\in[\dot p\dot z]$, we have \[\dist{\dot x}{\dot w}{}+\dist{\dot w}{\dot y}{}
<
\dist{\dot x}{\dot z}{}+\dist{\dot z}{\dot y}{}=\dist{x}{y}{}.\]
Suppose $w\in[p z]$ corresponds to $\dot w$; that is, $\dist{z}{w}{}=\dist{\dot z}{\dot w}{}$.
Since $\trig p x z$ and $\trig p y z$ are thin, we have 
\[\dist{x}{w}{}+\dist{w}{y}{}<\dist{x}{y}{},\]
contradicting the triangle inequality. 

Denote by $\dot D$ the union of two solid triangles $\trig {\dot p}{\dot x}{\dot z}$ and $\trig {\dot p} {\dot y} {\dot z}$.
Further, denote by $\tilde D$ the solid triangle $\trig{\tilde  p}{\tilde  x}{\tilde  y}=\modtrig(p x y)_{\EE^2}$.
By \ref{eq:<+<>=pi}, there is a \index{short map}\emph{short map}\footnote{In other words, distance-nonexpanding or 1-Lipschitz.}
$F\:\tilde D\to \dot D$ that sends 
\begin{align*}
\tilde p&\mapsto \dot p,
&
\tilde x&\mapsto \dot x,
&
\tilde z&\mapsto \dot z,
&
\tilde y&\mapsto \dot y.
\end{align*}

\begin{wrapfigure}{r}{50mm}
\vskip-2mm
\centering
\includegraphics{mppics/pic-830}
\vskip0mm
\end{wrapfigure}

Indeed, by Alexandrov's lemma (\ref{lem:alex}), 
there are nonoverlapping triangles 
\[\trig{\tilde p}{\tilde x}{\tilde z_x}\iso\trig {\dot p}{\dot x}{\dot z}\] 
and 
\[\trig{\tilde p}{\tilde y}{\tilde z_y}\iso\trig {\dot p}{\dot y}{\dot z}\]
inside the  triangle $\trig{\tilde p}{\tilde x}{\tilde y}$.

Connect  the points in each pair
$(\tilde z,\tilde z_x)$, 
$(\tilde z_x,\tilde z_y)$ 
and $(\tilde z_y,\tilde z)$ 
with arcs of circles centered at 
$\tilde y$, $\tilde p$, and $\tilde x$ respectively. 
Define $F$ as follows:
\begin{itemize}

\item Map  $\Conv\trig{\tilde p}{\tilde x}{\tilde z_x}$ isometrically onto  $\Conv\trig {\dot p}{\dot x}{\dot z}$;
similarly map $\Conv \trig{\tilde p}{\tilde y}{\tilde z_y}$ onto $\Conv \trig {\dot p}{\dot y}{\dot z}$.

\item If $x$ is in one of the three circular sectors, say at distance $r$ from its center, set $F(x)$ to be the point on the corresponding segment 
$[p z]$, 
$[x z]$ 
or $[y z]$ whose distance from the left-hand endpoint of the segment is~$r$.

\item Finally, if $x$ lies in the remaining curvilinear triangle $\tilde z \tilde z_x \tilde z_y$, 
set $F(x) = z$. 
\end{itemize}
By construction, $F$ satisfies the conditions.

By assumption, the natural maps $\trig {\dot p} {\dot x} {\dot z}\to\trig p x z$ and $\trig {\dot p} {\dot y} {\dot z}\to\trig p y z$ are short.  
By composition, the natural map from $\trig{\tilde  p}{\tilde  x}{\tilde  y}$ to $\trig p x y$ is short, as claimed.

The spherical case is done along the same lines.
\qeds

\section{Reshetnyak's gluing}\label{sec:cba-gluing}

Suppose 
$\spc{U}_1$ and $\spc{U}_2$ are spaces
with isometric closed sets $A_1\subset\spc{U}_1$ and $A_2\subset\spc{U}_2$.
Let $\iota\:A_1\to A_2$ be an isometry.
Consider the space $\spc{W}$ of all equivalence classes in $\spc{U}_1\sqcup\spc{U}_2$ with the equivalence relation given by $a\sim\iota(a)$ for any $a\in A_1$.
Let us equip $\spc{W}$ with the minimal metric such that both natural maps $\spc{U}_i\to \spc{W}$ are short.

The space $\spc{W}$ is called the \index{gluing}\emph{gluing} of $\spc{U}_1$ and  $\spc{U}_2$ along~$\iota$.
(A more general definition is discussed in \cite[1D]{petrunin-2023}.)

If $\spc{U}_1$ and $\spc{U}_2$ are proper and geodesic and $A_i$ are convex in $\spc{U}_i$,
then it is easy to check that $\spc{W}$ is a proper geodesic space as well and its metric can be described in the following way:
\[\begin{aligned}
\dist{x}{y}{\spc{W}}&\df\dist{x}{y}{\spc{U}_i}
\\
&\quad\text{if}\quad x,y\in \spc{U}_i,\quad\text{and}
\\
\dist{x}{y}{\spc{W}}&\df\min\set{\dist{x}{a}{\spc{U}_1}+\dist{y}{\iota(a)}{\spc{U}_2}}{a\in A_1}
\\
&\quad\text{if}\quad x\in \spc{U}_1\quad\text{and}\quad y\in \spc{U}_2.
\end{aligned}
\eqlbl{eq:gluing}
\]
Abusing notation we denote by $x$ and $y$ both the points in $\spc{U}_1\sqcup\spc{U}_2$ and their equivalence classes in $\spc{U}_1\sqcup\spc{U}_2/{{\sim}}$.

We can (and will) identify $\spc{U}_i$ with its image in $\spc{W}$;
this way both subsets $A_i\subset \spc{U}_i$ will be identified and denoted further by~$A$.
Note that $A=\spc{U}_1\cap \spc{U}_2\subset \spc{W}$,
therefore $A$ is also a convex set in~$\spc{W}$.

{\sloppy

\begin{thm}{Reshetnyak gluing}\label{thm:gluing}
Suppose 
$\spc{U}_1$ and $\spc{U}_2$ are proper geodesic $\CAT(0)$ spaces
with isometric 
closed 
convex
sets $A_i\subset\spc{U}_i$, and $\iota\:A_1\z\to A_2$ is an isometry.
Then the gluing of $\spc{U}_1$ and  $\spc{U}_2$ along $\iota$ is a $\CAT(0)$ proper geodesic space.
\end{thm}

}

The condition that the spaces are proper will be needed only in order to use formula \ref{eq:gluing}.

\parit{Proof.} 
By construction of the gluing space, the statement can be reformulated in the following way:

\begin{thm}{Reformulation of \ref{thm:gluing}}
Let $\spc{W}$ be a 
proper geodesic space with two closed 
convex sets $\spc{U}_1,\spc{U}_2\subset\spc{W}$ such that
$\spc{U}_1\cup\spc{U}_2=\spc{W}$
and $\spc{U}_1$, $\spc{U}_2$ are $\CAT(0)$.
Then $\spc{W}$ is $\CAT(0)$.
\end{thm}

It suffices to show that any triangle $\trig {x}{y}{z}$ 
in $\spc{W}$ is thin.
This is obviously true if all three points $x$, $y$, $z$ lie in one of~$\spc{U}_i$.
Thus, without loss of generality, we may assume that $x\z\in\spc{U}_1$ and $y,z\z\in\spc{U}_2$.

\begin{wrapfigure}{o}{45mm}
\vskip-2mm
\centering
\includegraphics{mppics/pic-840}
\end{wrapfigure}

Choose points $a,b\in A\z=\spc{U}_1\z\cap\spc{U}_2$
that lie respectively on the sides $[xy], [xz]$.
Note that

\begin{itemize}
\item the triangle $\trig{x}{a}{b}$ lies in $\spc{U}_1$,
\item both triangles $\trig{y}{a}{b}$ and $\trig{y}{b}{z}$ lie in~$\spc{U}_2$.
\end{itemize}
In particular, each triangle $\trig{x}{a}{b}$,
$\trig{y}{a}{b}$, and $\trig{y}{b}{z}$ is thin.

Applying the inheritance lemma (\ref{lem:inherit-angle}) twice, 
we get that $\trig {x}{y}{b}$ 
and consequently $\trig {x}{y}{z}$ are thin.
\qeds

If one applies gluing to two copies of a single space $\spc{U}$ with a closed subset $A\subset \spc{U}$ and the identity map $\iota\:A\to A$, then the obtained space is called the \index{double}\emph{double} of $\spc{U}$ along~$A$.

\begin{thm}{Exercise}\label{ex:reshetnyak-doubling}
Suppose $\spc{W}$ is the double of a geodesic space $\spc{U}$ along its closed subset $A$.
Show that $\spc{W}$ is $\CAT(0)$ if and only if $\,\spc{U}$ is $\CAT(0)$, and $A$ is convex in $\spc{U}$.
\end{thm}

\section{Notes}

The motivation for the notion of $\CAT(\kappa)$ spaces comes from the fact that a Riemannian manifold is locally  $\CAT(\kappa)$ if and only if it has sectional curvature at most $\kappa$.
This easily follows from Rauch comparison for Jacobi fields and Proposition~\ref{prop:thin=cat}.

The gluing theorem (\ref{thm:gluing}) was proved by Yuri Reshetnyak~\cite{reshetnyak-1960}.
Applying the argument in \ref{ex:tan-CBA:geodesic}, it can be extended to all geodesic $\CAT(0)$ spaces (not necessarily proper).
It also admits a natural generalization to 
geodesic $\CAT(\kappa)$
spaces for arbitrary $\kappa$;
see the book of Martin Bridson and André Haefliger \cite{bridson-haefliger} and our book \cite{alexander-kapovitch-petrunin-2025} for details.

%%%%%%%%%%%%%%%%%%%%%%%%%%%%%%%%%%%%%%%%%%%%%%%%%%%%%%%%%%%%%%%%%%%%%%%%%%%%%%%%%%%%%%%%%%%%%%%%%

\chapter{Billiards}\label{chap:billiards}

Now we will look at an application of the previous lecture that is at once beautiful, deep, elementary, and unexpected---enjoy.

\section{Puff pastry}\label{sec:puff-pastry}

Here we introduce the notion of Reshetnyak puff pastry.
This construction will be applied to billiards (\ref{thm:collision}).

Let $\bm{A}=(A_1,\ldots,A_N)$ be an array of closed convex sets in the Euclidean space~$\EE^m$.
Consider an array of $N+1$ copies of~$\EE^m$.
Let us construct the space $\spc{R}$ by gluing successive pairs of spaces along  $A_1,\ldots,A_N$ respectively.

\begin{figure}[ht!]
\vskip-0mm
\centering
\includegraphics{mppics/pic-850}
\caption*{Puff pastry for $(A,B,A)$.}
\end{figure}
% superindexes in the pic should be changed to subindexes
The resulting space $\spc{R}$  will be called 
the
\index{puff pastry}\emph{Reshetnyak puff pastry} for array~$\bm{A}$.
The copies of $\EE^m$ in the puff pastry $\spc{R}$
will be called {}\emph{levels};
they will be denoted by $\spc{R}_0,\ldots,\spc{R}_N$.
The point in the $k$-th level $\spc{R}_k$
that corresponds to $x\in \EE^m$
will be denoted by~$x_k$.

Given $x\in \EE^m$, any point $x_k\in\spc{R}$ is called a {}\emph{lifting} of~$x$.
The map $x\mapsto x_k$ defines an isometry $\EE^m\to \spc{R}_k$;
in particular, we can talk about liftings of subsets in~$\EE^m$.

Note that: 
\begin{itemize}
\item The intersection $A_1\cap\ldots\cap A_N$ admits a unique lifting in~$\spc{R}$.
\item Moreover, $x_i=x_j$ for some $i<j$
if and only if 
\[x\in A_{i+1}\cap\ldots\cap A_j.\]
\item For any $k$, the restriction $\spc{R}_k\to \EE^m$
of the natural projection $x_k\mapsto x$ is an isometry.
\end{itemize}

\begin{thm}{Observation}\label{obs:puff pastry is CAT}
Any Reshetnyak puff pastry is a proper geodesic $\CAT(0)$ space.
\end{thm}

\parit{Proof.} Apply Reshetnyak gluing theorem (\ref{thm:gluing}) recursively for the convex sets in the array.
\qeds

\begin{thm}{Proposition}\label{prop:A-check-A}
Assume $(A_1,\ldots,A_N)$ and $(\check A_1,\ldots,\check A_N)$ are two arrays of closed convex sets in $\EE^m$
such that $ A_k\subset \check A_k$ for each~$k$.
Let $\spc{R}$ and $\check{\spc{R}}$ be the corresponding Reshetnyak  puff pastries.
Then the map $\spc{R}\to\check{\spc{R}}$
defined by $x_k\mapsto\check x_k$ is well-defined and short.

Moreover, if  
\[\dist{x_i}{y_j}{\spc{R}}=\dist{\check x_i}{\check y_j}{\check{\spc{R}}}\eqlbl{eq:dist=dist}\]
for some $x,y\in \EE^m$ and $i,j\in \{0,\ldots,N\}$, then the necessarily unique geodesic $[\check x_i \check y_j]_{\check{\spc{R}}}$ is the image of the necessarily unique geodesic $[x_i y_j]_{\spc{R}}$ under the map $x_i\mapsto \check x_i$.
\end{thm}

\parit{Proof.}
The first statement in the proposition 
follows from the construction of Reshetnyak  puff pastries.

By Observation~\ref{obs:puff pastry is CAT}, 
$\spc{R}$  and  $\check{\spc{R}}$ are proper geodesic $\CAT(0)$ spaces; 
hence $[x_i y_j]_{\spc{R}}$
and $[\check x_i \check y_j]_{\check{\spc{R}}}$ are unique.
By \ref{eq:dist=dist}, since the map $\spc{R}\to\check{\spc{R}}$ is short, 
the image of $[x_i y_j]_{\spc{R}}$
is a geodesic of $\check{\spc{R}}$ joining $\check x_i$ to~$\check y_j$.
Hence the second statement follows.
\qeds

%??? Sasha suggests to add the following def:
%Say that a puff-pastry is "simple" (or whatever) if there exists a geodesic connecting a point from the first with a point in the last  leave and which does not cut $A_i\cap A_j$  for no $i,j$.
%The theorem about pastry you prove is that if the angles are large there are no high simple puff pastries. 

\begin{thm}{Definition}
Consider a Reshetnyak puff pastry $\spc{R}$ with the levels 
$\spc{R}_0,\ldots,\spc{R}_N$.
We say that $\spc{R}$ is \index{end-to-end convex}\emph{end-to-end convex} 
if $\spc{R}_0\cup\spc{R}_N$, the union of its lower and upper levels,
forms a convex set in~$\spc{R}$;
that is, if $x,y\in \spc{R}_0\cup\spc{R}_N$, then $[xy]_{\spc{R}}\subset  \spc{R}_0\cup\spc{R}_N$.
\end{thm}

Note that if $\spc{R}$ is the Reshetnyak puff pastry for an array of convex sets $\bm{A}=(A_{1},\ldots, A_{N})$,
then $\spc{R}$ is end-to-end convex
if and only if the union of the lower and the upper levels
$\spc{R}_0\cup\spc{R}_N$ is isometric to the double of $\EE^m$ along the nonempty intersection $A_1\cap\ldots\cap A_N$.

\begin{thm}{Observation}\label{obs:end-to-end-convex}
Let $\check{\bm{A}}$ and $\bm{A}$ be arrays of convex bodies in~$\EE^m$.
Assume that array $\bm{A}$ is
obtained by inserting in $\check{\bm{A}}$ 
several copies of the bodies which were already listed in~$\check{\bm{A}}$.%
\footnote{For example, inserting $B$ and $A$ in $\check{\bm{A}}=(A,C,B,C,A)$, we obtain $\bm{A}\z=(A,B,C,A,B,C,A)$.}
Denote by $\check{\spc{R}}$ and $\spc{R}$ 
the Reshetnyak puff pastries for $\check{\bm{A}}$ and $\bm{A}$ respectively.

If $\check{\spc{R}}$ is end-to-end convex, then so is~$\spc{R}$.
\end{thm}

\parit{Proof.}
Without loss of generality, we may assume that $\bm{A}$ is 
obtained by inserting one element in $\check{\bm{A}}$,
say at the place number~$k$.

Note that $\check{\spc{R}}$ is isometric to the puff pastry 
for $\bm{A}$ with $A_k$ replaced by~$\EE^m$.
It remains to apply Proposition~\ref{prop:A-check-A}.
\qeds

{

\begin{wrapfigure}{o}{30mm}
\vskip-4mm
\centering
\includegraphics{mppics/pic-860}
\end{wrapfigure}

Let $X$ be a convex subset in a Euclidean space.
By a \index{dihedral angle}\emph{dihedral angle}, we understand an intersection of two half-spaces;
the intersection of corresponding hyperplanes is called the {}\emph{edge} of the angle.
We say that a dihedral angle $D$ 
supports
 $X$ at a point $p\in X$ 
if $D$ contains $X$ and the edge of $D$ contains~$p$.

}

\begin{thm}{Lemma}\label{lem:end-to-end-convex}
Let $A$ and $B$ be two convex sets in~$\EE^m$.
Assume that any dihedral angle supporting $A\cap B$ has angle measure at least~$\alpha$.
Then the Reshetnyak puff pastry for the array
\[(\underbrace{A,B,A,\ldots}_{\text{$\lceil\tfrac\pi\alpha\rceil+1$ times}}).\]
is end-to-end convex. 
\end{thm}

The proof of the lemma is based on its partial case,
which we formulate as the following sublemma.

\begin{thm}{Sublemma}\label{sublem:end-to-end-convex}
Let $\ddot A$ and $\ddot B$ be two  
half-planes in $\EE^2$, where $\ddot A\cap \ddot B$ is an angle with measure~$\alpha$.
Then the Reshetnyak puff pastry for the array \[(\underbrace{\ddot A,\ddot B,\ddot A,\ldots}_{\text{$\lceil\tfrac\pi\alpha\rceil+1$ times}})\]
is end-to-end convex. 
\end{thm}

\parit{Proof.}
Note that the puff pastry $\ddot{\spc{R}}$ is isometric to the cone over the space glued from the unit circles as shown on the diagram.

All the short arcs on the diagram have length $\alpha$;
the long arcs have length $\pi-\alpha$,
so making a circuit along any path will take~$2\cdot\pi$.

{

\begin{wrapfigure}{r}{35mm}
\vskip0mm
\centering
\includegraphics{mppics/pic-870}
\end{wrapfigure}

Applying Exercise \ref{ex:cone-geod}, we get that 
the end-to-end convexity of $\ddot{\spc{R}}$ follows if any geodesic shorter than $\pi$ with the ends on the inner and the outer circles lies completely in the union of these two circles.
The latter holds if the zigzag line in the picture has length at least~$\pi$.
This line is formed by $\lceil\tfrac\pi\alpha\rceil$ arcs with length $\alpha$ each.
Hence the sublemma.
\qeds

}

In the proof of \ref{lem:end-to-end-convex}, we will use the following exercise in convex geometry:
{

\begin{wrapfigure}{r}{26mm}
\vskip-3mm
\centering
\includegraphics{mppics/pic-880}
\vskip-4mm
\end{wrapfigure}

\begin{thm}[!]{Exercise}\label{ex:supporting-planes}
Let $A$ and $B$ be two closed convex sets in $\EE^m$.
Given two points $x,y\in \EE^m$, consider the function $f(p)\df\dist{x}{p}{}+\dist{y}{p}{}$.
Suppose point $z$ minimizes $f$ in $A\cap B$.

Show that there are half-spaces $\dot A$ and $\dot B$ such that
$\dot A\supset A$ and $\dot B\supset B$
and $z$ also minimizes $f$ in $\dot A\cap \dot B$.

\end{thm}

\parit{Proof of \ref{lem:end-to-end-convex}.}
Fix arbitrary $x,y\in \EE^m$.
Choose a point $z\in A\cap B$
for which the sum 
\[\dist{x}{z}{}+\dist{y}{z}{}\] 
is minimal.
To show the end-to-end convexity of  $\spc{R}$,
it is sufficient to prove the following:

}

\begin{clm}{}\label{clm:z in xy}
The geodesic $[x_0y_N]_\spc{R}$ contains $z_0=z_N\in \spc{R}$.
\end{clm}

Without loss of generality, we may assume that $z\in\partial A\cap\partial B$.
Indeed, assume $z$ lies in the interior of~$A$.
Since the puff pastry for the 1-array $(B)$ is end-to-end convex,
Proposition~\ref{prop:A-check-A} together with \ref{obs:end-to-end-convex}
imply \ref{clm:z in xy}.
The same way we can treat the case when $z$ lies in the interior of~$B$.

Suppose that $\dot A$ and $\dot B$ are provided by Exercise \ref{ex:supporting-planes}.
Note that $\EE^{m}$ admits
an
isometric splitting $\EE^{m-2}\times \EE^2$
such that
\begin{align*}
\dot A&=\EE^{m-2}\times \ddot A
\\
\dot B&=\EE^{m-2}\times \ddot B
\end{align*}
where $\ddot A$ and $\ddot B$ are half-planes in~$\EE^2$.

Let us replace each $A$ by $\dot A$ and each $B$ by $\dot B$
in the array, to get the array
\[(\underbrace{\dot A,\dot B,\dot A,\ldots}_{\text{$\lceil\tfrac\pi\alpha\rceil+1$ times}}).\]
The corresponding puff pastry $\dot{\spc{R}}$
splits as a product of $\EE^{m-2}$ and a puff pastry, 
call it $\ddot{\spc{R}}$,
glued from the copies of the plane $\EE^2$ for the array
\[(\underbrace{\ddot A,\ddot B,\ddot A,\ldots}_{\text{$\lceil\tfrac\pi\alpha\rceil+1$ times}}).\]

Note that the dihedral angle $\dot A\cap \dot B$ is at least~$\alpha$.
Therefore the angle measure of  $\ddot A\cap \ddot B$ is also at least $\alpha$.
According to Sublemma~\ref{sublem:end-to-end-convex} and Observation~\ref{obs:end-to-end-convex}, $\ddot{\spc{R}}$ is end-to-end convex.

Since $\dot{\spc{R}}\iso\EE^{m-2}\times\ddot{\spc{R}}$, 
the puff pastry $\dot{\spc{R}}$ is also end-to-end convex;
see \ref{ex:geod-prod}.

It follows that the geodesic $[\dot x_0\dot y_N]_{\dot{\spc{R}}}$ contains $\dot z_0=\dot z_N\in\dot{\spc{R}}$.
By Proposition~\ref{prop:A-check-A}, 
the image of $[\dot x_0\dot y_N]_{\dot{\spc{R}}}$
under the map $\dot x_k\mapsto x_k$
is the geodesic $[x_0 y_N]_{\spc{R}}$.
Hence \ref{clm:z in xy} and the lemma follow.
\qeds

\section{Wide corners}
\label{sec:wide-corners}

\begin{wrapfigure}{r}{26mm}
\vskip-4mm
\centering
\includegraphics{mppics/pic-885}
\end{wrapfigure}

We say that a closed convex set $A\subset \EE^m$ has {}\emph{$\eps$-wide corners}\label{page:wide corners} for given $\eps >0$
if together with each point $p$, 
the set $A$ contains a small right circular cone
with the tip at $p$ and aperture $\eps$;
that is, $\eps$ is the maximum angle between two generating lines of the cone.

For example, 
a plane polygon 
has $\eps$-wide corners
if all its interior angles are at least~$\eps$.

We will consider arrays of closed convex sets 
$A_1,\ldots,A_n$ in $\EE^m$
such that for any subset $F\subset\{1,\ldots,n\}$,
the intersection
$\bigcap_{i\in F}A_i$
has $\eps$-wide corners.
In this case, we may say briefly that \index{wide corners}\textit{all intersections of $A_i$ have $\eps$-wide corners}.

\begin{thm}{Exercise}\label{ex:compact-walls}
Assume $A_1,\ldots,A_n\subset\EE^m$ are compact, convex sets with a common interior point.
Show that all intersections of $A_i$ have $\eps$-wide corners for some~$\eps>0$.
\end{thm}

\begin{thm}[!]{Exercise}\label{ex:centrally-simmetric-walls}
Assume $A_1,\ldots,A_n\subset\EE^m$ are
convex sets with nonempty interiors that have a common center of symmetry.
Show that all intersections of $A_i$ have $\eps$-wide corners for some~$\eps>0$.
\end{thm}

The proof of the following proposition is based on lemma \ref{lem:end-to-end-convex}, which is essentially the case $n=2$ in the proposition.

\begin{thm}{Proposition}\label{prop:end-to-end-convex}
Given $\eps>0$ and a positive integer $n$, 
there is an array of integers $\bm{j}_\eps(n)=(j_1,\ldots,j_N)$
such that: 

\begin{subthm}{} For each $k$ we have $1\le j_k\le n$,
and each number $1,\ldots,n$ appears in $\bm{j}_\eps$ at least once.
\end{subthm}

\begin{subthm}{}
If $A_1,\ldots,A_n$ is a collection of closed convex sets in $\EE^m$ with a common point
and all their intersections have $\eps$-wide corners,  
then the puff pastry for the array
$(A_{j_1},\ldots,A_{j_N})$ is end-to-end convex.
\end{subthm}

Moreover, we can assume that $N\le (\lceil\tfrac\pi\eps\rceil+1)^n$.
\end{thm}

\parit{Proof.}
The array $\bm{j}_\eps(n)=(j_1,\ldots,j_N)$  is constructed recursively.
For $n=1$, we can take $\bm{j}_\eps(1)=(1)$.

Assume that $\bm{j}_\eps(n)$ is constructed.
Let us replace each occurrence of $n$ in $\bm{j}_\eps(n)$ by the alternating string 
\[\underbrace{n,n+1,n,\ldots}_{\text{$\lceil\tfrac\pi\eps\rceil+1$ times}}\]
Denote the obtained array by $\bm{j}_\eps(n+1)$.

By Lemma \ref{lem:end-to-end-convex},
the end-to-end convexity of the puff pastry for $\bm{j}_\eps(n\z+1)$
follows from the end-to-end convexity of the puff pastry for the array
where each string
\[\underbrace{A_n,A_{n+1},A_n,\ldots}_{\text{$\lceil\tfrac\pi\eps\rceil+1$ times}}\]
is replaced by  $Q=A_n\cap A_{n+1}$.
End-to-end convexity of the latter follows by the assumption on $\bm{j}_\eps(n)$, 
since all the intersections of $A_1,\ldots,A_{n-1},Q$
have $\eps$-wide corners.

The upper bound on $N$ follows directly from the construction.
\qeds

\section{Billiards}

Let $A_1,A_2,\ldots, A_n$ be a finite collection of closed convex bodies in~$\EE^m$.
Assume that for each $i$
the boundary $\partial A_i$ is a smooth hypersurface.

Consider the billiard table formed by the closure of the complement 
$$T=\overline{\EE^m\setminus \bigcup_{i} A_i}.$$
The sets $A_i$ will be called {}\emph{walls} of the table
and the billiards described above will be called {}\emph{billiards with convex walls}.

A billiard {}\emph{trajectory} on the table is a unit-speed polygonal line $\gamma$ that follows the standard law of billiards at the breakpoints on $\partial A_i$
--- in particular, the angle of reflection is equal to the angle of incidence.
The breakpoints of the trajectory will be called {}\emph{collisions}.
\textit{We assume that our trajectory (1) never meets a wall in a tangent direction and (2) meets at most one wall in a small neighborhood of every time $t$.
If $\gamma$ meets a wall tangentially,
or meets two or more walls at time $t$,
or collisions accumulate at $t$,
then we declare that $\gamma$ dies right before $t$ (that is, $\gamma$ is not defined at times  $\ge t$).}

Recall that the definition of sets with $\eps$-wide corners is given in \ref{sec:wide-corners}.

\begin{thm}{Collision theorem}\label{thm:collision}
Let $T\subset\EE^m$ be a billiard table with $n$ convex walls.
Assume that the walls of $T$ have a common interior point  and all their intersections have $\eps$-wide corners.
Then the number of collisions of any trajectory in  $T$  is bounded
by a number $N$ which depends only on $n$ and~$\eps$.
\end{thm}

As we will see from the proof, the value $N$ can be found explicitly;
say, $N=(\lceil\tfrac\pi\eps\rceil+1)^{n^2}$ will do.

\begin{thm}{Corollary}\label{cor:balls}
Consider $n$ homogeneous hard balls
moving freely and colliding
elastically in~$\RR^3$. 
Every ball moves
along a straight line with constant speed until two balls collide, and then
the new velocities of the two balls are determined by the
laws of classical mechanics. 
We assume that only two balls can collide at the same time.

Then the total number of collisions cannot exceed some number $N$ that  depends on the radii and masses of the balls.
If the balls are identical, then $N$ depends only on~$n$.
\end{thm}

The proof below admits a straightforward generalization to all dimensions.

\begin{thm}{Exercise}\label{cor:balls:dim=1}
Show that in the case of identical balls in the 1-dimensional space (that is, in $\RR$)
the total number of collisions cannot exceed $N=\tfrac{n\cdot(n-1)}2$.
\end{thm}

\parit{Proof of \ref{cor:balls} modulo \ref{thm:collision}.}
Denote by $a_i=(x_i,y_i,z_i) \in \RR^3$ the center of the $i$-th ball.
Consider the corresponding point in $\RR^{3\cdot n}$
\begin{align*}
\bm{a}&=(a_1, a_2 , \ldots , a_n ) =
\\
&=(x_1, y_1 , z_1 , x_2 , y_2 , z_2 , \ldots , x_n , y_n , z_n).
\end{align*}

The $i$-th and $j$-th balls intersect if 
$$|a_i - a_j | \le R_i+R_j,$$
where $R_i$ denotes the radius of the $i$-th ball.
These inequalities define $\tfrac{n\cdot(n-1)}{2}$ cylinders 
\[C_{i,j}=\set{(a_1, a_2 , \ldots , a_n )\in\RR^{3\cdot n}} {|a_i - a_j |\le R_i+R_j}.\]
The closure of the complement
\[T=\overline{\RR^{3\cdot n}\setminus \bigcup_{i< j} C_{i,j}}\] 
is the configuration space of our system. 
Its points correspond
to valid positions of the system of balls.

The evolution of the system is described by the motion of
the point $\bm{a}\in\RR^{3\cdot n}$.
It moves along a straight line at a
constant speed until it hits one of the cylinders $C_{i,j}$; 
this event corresponds
to a collision in the system of balls.

Consider the norm of $\bm{a}=(a_1,\ldots,a_n)\in \RR^{3\cdot n}$ defined by
\[\lVert \bm{a}\rVert
=
\sqrt{M_1\cdot|a_1|^2+\ldots+M_n\cdot |a_n|^2},\]
where $|a_i|=\sqrt{x_i^2+y_i^2+z_i^2}$ 
and $M_i$ denotes the mass of the $i$-th ball.
It is straightforward to check that for the metric defined by $\lVert {*}\rVert$, the collisions follow the standard law of billiards.

By construction, the number of collisions of hard balls that we need to estimate 
is the same as the number of collisions of the corresponding billiard trajectory on the table with $C_{i,j}$ as the walls.

Note that each cylinder $C_{i,j}$ is a convex set;
it has smooth boundary, 
and it is centrally symmetric around the origin.
By \ref{ex:centrally-simmetric-walls}, all the intersections of the walls have $\eps$-wide corners for some $\eps>0$ that depends on the radii $R_i$ and the masses~$M_i$.
It remains to apply the collision theorem (\ref{thm:collision}).
\qeds

Finally, we present the proof of the collision theorem (\ref{thm:collision})
based on the results developed in the previous section.

\parit{Proof of \ref{thm:collision}.}
Let us apply induction on~$n$.

\parit{Base: $n=1$.}
The number of collisions cannot exceed~1.  
Indeed, by the convexity of $A_1$,
if the trajectory is reflected once in $\partial A_1$,
then it cannot return to~$A_1$.

\parit{Step.}
Assume $\gamma$ is a trajectory that  meets the walls in the order $A_{i_1},\ldots,A_{i_N}$ for a large integer~$N$.

Consider the array 
\[\bm{A}^\gamma=(A_{i_1},\ldots,A_{i_N}).\]
The induction hypothesis implies:

\begin{clm}{}\label{clm:collision-induction hypothesis}
There is a positive integer $N_0$ such that any $N_0$ consecutive elements of $\bm{A}^\gamma$ contain each $A_i$ at least once.
\end{clm}

Let $\spc{R}^\gamma $ be  the  Reshetnyak puff pastry for~$\bm{A}^\gamma$.

Consider the lift of $\gamma$ to $\spc{R}^\gamma$,
defined by 
$\bar\gamma(t)=\gamma_k(t)\in \spc{R}^\gamma$
for any moment of time $t$ between the $k$-th and $(k+1)$-th collisions.  
Since $\gamma$ follows the standard law of billiards at breakpoints, the lift $\bar\gamma$ is locally a geodesic in~$\spc{R}^\gamma$.
By \ref{obs:puff pastry is CAT},
the puff pastry $\spc{R}^\gamma$ is a proper geodesic $\CAT(0)$ space.
Therefore $\bar\gamma$ is a geodesic.

Since $\gamma$ does not meet $A_1\cap\ldots\cap A_n$,
the lift $\bar\gamma$ does not lie in  $\spc{R}^\gamma_0\cup \spc{R}^\gamma_N$.
In particular, $\spc{R}^\gamma$ is not end-to-end convex.

Let 
\[\bm{B}=(A_{j_1},\ldots,A_{j_K})\]
be the array provided by Proposition~\ref{prop:end-to-end-convex};
so $\bm{B}$ contains each $A_i$ at least once
and the puff pastry $\spc{R}_{\bm{B}}$ for $\bm{B}$ is end-to-end convex.
If $N$ is sufficiently large, namely $N\ge K\cdot N_0$, then
 \ref{clm:collision-induction hypothesis}
implies that $\bm{A}^\gamma$ can be obtained
by inserting a finite number of $A_i$'s in~$\bm{B}$.

By \ref{obs:end-to-end-convex}, 
$\spc{R}^\gamma$ is end-to-end convex --- a contradiction.
\qeds

\begin{thm}{Exercise}\label{ex:collision}
Let $T\subset\EE^m$ be a billiard table with $n$ convex walls.
Assume that the walls of $T$ have a common point.
Show that any trajectory in $T$ has only a finite number of collisions in a finite time interval.

Construct an example of a billiard table on the plane with two convex walls
such that the walls have a common point,
but there is no upper bound on the number of collisions of its trajectories.
\end{thm}

\section{Notes}

Note that the interior points of the walls play a key role in the proof
despite the fact  that the trajectories never go inside the walls.

The collision theorem (\ref{thm:collision}) was proved by Dmitri Burago, Serge Ferleger and Alexey Kononenko \cite{burago-ferleger-kononenko-1997}.
Its corollary (\ref{cor:balls}) answers a question posed by Yakov Sinai \cite{galperin}.
Puff pastry is used to bound topological entropy of the billiard flow 
and to approximate the shortest billiard path that  touches given lines in a given order \cite{burago-ferleger-kononenko-1998,burago-grigoriev-slissenko}.
The lecture of Dmitri Burago \cite{burago-1998} gives a short survey on the subject.
Puff pastry was also used by the first author and Richard Bishop \cite{alexander-bishop-1998-Warped}
to find the upper curvature bound for warped products.

Joel Hass \cite{hass} constructed an example of a Riemannian metric on the closed 3-ball with negative curvature and concave boundary.
This example might decrease your appetite for generalizing the collision theorem --- while
locally such a 3-ball looks as good as the billiard table in the theorem, the number of collisions is obviously infinite.

It was shown by Dmitri Burago and Sergei Ivanov \cite{burago-ivanov} that the number of collisions that may occur between $n$ identical balls in $\RR^3$ grows at least exponentially in $n$; the 2-dimensional case remains open.

Roman Barinov and Sergei Ivanov used another method to prove analogous statements for normed spaces with strongly convex smooth norms \cite{barinov-ivanov}.

%%!TEX root = the-majorization.tex
\chapter{Majorization}\label{sec:resh-kirz}

This lecture provides a handy, intuitive tool that will make your life with $\CAT(0)$ spaces easier and pleasantly addictive.

\section{Formulation}

\begin{thm}{Definition}\label{def:majorize}
Let $\spc{X}$ be a metric space,
$\tilde \alpha$ be a simple closed curve of finite length  in $\EE^2$,
and $D\subset\EE^2$ be a closed region bounded by $\tilde \alpha$.
A length-nonincreasing map $F\:D\to\spc{X}$ is called \index{majorizing map}\emph{majorizing} if it is length-preserving on $\tilde \alpha$.

In this case, we say that $D$ \emph{majorizes} the curve $\alpha=F\circ\tilde \alpha$ under the map $F$.
\end{thm}

The following proposition is a consequence of the definition.

\begin{thm}{Proposition}
\label{prop:majorize-geodesic} 
Let  $\alpha$  be a closed curve in a metric space $\spc{X}$.
Suppose $D\subset\EE^2$ majorizes $\alpha$ under $F\: D \to \spc{X}$.  
Then any geodesic subarc of $\alpha$ is the image under $F$ of a subarc of $\partial D$ that is geodesic in the length metric of $D$.

In particular, if $D$ is convex, then the corresponding subarc is a geodesic in $\EE^2$.
\end{thm}

\parit{Proof.}
Given a geodesic subarc $\gamma\:[a,b]\to\spc{X}$ of $\alpha$,
let $\tilde \gamma = (F|_{\partial D})^{-1}\circ\gamma$,
\begin{align*}
r&=\dist{\gamma(a)}{\gamma(b)}{\spc{X}},
&
\tilde r&=\dist{\tilde \gamma(a)}{\tilde \gamma(b)}{D},
\\
s&=\length \gamma,
&
\tilde s&= \length \tilde \gamma.
\end{align*}
Then $\tilde r\ge r = s =\tilde s\ge\tilde r$.
Therefore $\tilde s=\tilde r$, and the statement follows.
\qeds

\begin{thm}{Corollary}\label{cor:maj-triangle}
Assume a convex region $D\subset \EE^2$ majorizes a triangle $\trig p x y$ in a metric space $\spc{X}$.
Then $D$ is a \emph{solid model triangle} of $\trig p x y$;
that is,
$D=\Conv\trig{\tilde p}{\tilde x}{\tilde y}$ for a model triangle $\trig{\tilde p}{\tilde x}{\tilde y}=\modtrig(p x y)$.
Moreover, the majorizing map sends  $\tilde p$, $\tilde x$ and $\tilde y$ respectively to $p$, $x$ and $y$.
\end{thm}

Now we come to the main theorem of this section.

\begin{thm}{Majorization theorem}\label{thm:majorization}
\label{thm:major}
Any closed rectifiable curve $\alpha$ in a geodesic $\CAT(0)$ space is majorized by a convex plane figure.
\end{thm}

\begin{thm}{Exercise}\label{ex:convex-dist:major}
Solve \ref{ex:convex-dist} using the majorization theorem.
\end{thm}

\section{Triangles}

The case when $\alpha$ is a triangle, say $\trig p x y$, is the base in the following proof, and it is nontrivial.
In this case, by Corollary~\ref{cor:maj-triangle}, the majorizing convex region is the solid model triangle.

\begin{thm}{Line-of-sight map} \label{def:sight}
Let $p$ be a point and $\alpha$ be a curve of finite length in a geodesic $\CAT(0)$ space~$\spc{U}$.
Let $\mathring\alpha:[0,1]\to\spc{U}$ be the constant-speed parametrization of~$\alpha$.  
If $\gamma_t\:[0,1]\to\spc{U}$ is a geodesic path from $p$ to $\mathring\alpha(t)$, we say that 
\[
[0,1]\times[0,1]\to\spc{U}\:(t,s)\mapsto\gamma_t(s)
\]
is the \index{line-of-sight map}\emph{line-of-sight map from $p$ to $\alpha$}.  
\end{thm}

Note that a line-of-sight map is essentially a geodesic homotopy from a constant map to a curve.

We will show that there is a majorizing map for $\trig p x y$ whose image $W$ is the image of the line-of-sight map for $[x y]$ from  $p$.
However, as one can see from the following example, the line-of-sight map might \textit{not} be majorizing.

\begin{wrapfigure}{r}{30 mm}
\vskip-0mm
\centering
\includegraphics{mppics/pic-951}
\end{wrapfigure}

\parbf{Example.} Let $\spc{Q}$ be a solid quadrangle $[p x z y]$ in $\EE^2$ formed by two congruent triangles, which is non-convex at $z$ (as in the picture).  
Equip $\spc{Q}$ with the length metric. 
By Reshetnyak gluing (\ref{thm:gluing}), $\spc{Q}$ is $\CAT(0)$.
For the triangle ${\trig p x y}_\spc{Q}$ in $\spc{Q}$ and its model triangle $\trig{\tilde p}{\tilde x}{\tilde y}$ in $\EE^2$,  
we have 
\[\dist{\tilde x}{\tilde y}{}=\dist{x}{y}{\spc{Q}}=\dist{x}{z}{}+\dist{z}{y}{}.\]
Then the map $F$ defined by matching the line-of-sight parameters satisfies $F(\tilde x)=x$ and $\dist{x}{F(\tilde w)}{}>\dist{\tilde x}{\tilde w}{}$ if $\tilde w$ is near the midpoint $\tilde z$ of $[\tilde x\tilde y]$ and lies on $[\tilde p\tilde z]$. 
Indeed,
\begin{align*}
\dist{\tilde x}{\tilde w}{}
&=\dist{\tilde x}{\tilde \gamma_\frac12(s)}{}
=\dist{x}{z}{}+o(1-s),
\\
\dist{x}{F(\tilde w)}{}
&=\dist{x}{\gamma_\frac12(s)}{}
=\dist{x}{z}{}+\eps\cdot(1-s)+o(1-s),
\end{align*}
where $\eps=-\dist{p}{z}{}\cdot\cos\hinge zpx>0$.
Thus $\dist{\tilde x}{\tilde w}{}<\dist{x}{F(\tilde w)}{}$ if $s$ is slightly smaller than $1$; therefore, $F$ is not majorizing.

\begin{wrapfigure}{r}{38 mm}
\vskip-0mm
\centering
\includegraphics{mppics/pic-811}
\end{wrapfigure}

\begin{thm}{Definition}\label{def:convex-devel}
Let $\tilde \gamma\:\II\to\EE^2$ be a curve.
Assume that for some point $\tilde p\in\EE^2$ the direction of $[\tilde p\,\tilde \gamma(t)]$ turns monotonically as $t$ grows.

The set formed by all geodesics from  $\tilde p$ to the points on $\tilde \gamma$ is called the \index{development!subgraph/supergraph} \emph{subgraph} of $\tilde \gamma$ with respect to $\tilde p$.

The set of all points $\tilde x\in\EE^2$ such that a geodesic $[\tilde p\tilde x]$ intersects $\tilde \gamma$ is called the \emph{supergraph} of $\tilde \gamma$ with respect to $\tilde p$.

The curve $\tilde \gamma$ is called \index{convex/concave curve with respect to a point}\emph{convex} (\emph{concave}) with respect to $\tilde p$ if the subgraph (respectively, supergraph) of $\tilde \gamma$ with respect to $\tilde p$ is convex.
\end{thm}

\begin{thm}{Lemma}\label{lem:majorize-subgraph}
Let $\beta$ be a plane curve from $x$ to $y$ that is concave with respect  to $p$, and let $D$  be the subgraph of $\beta$ with respect to $p$.
\begin{subthm}{curvilinear} 
Then $\beta$ forms a geodesic $[x y]_D$ in $D$ and therefore $\beta$, $[p x]$ and $[p y]$ form a triangle 
${\trig p x y}_D$ in the length metric of $D$.
\end{subthm}
\begin{subthm} {short-to-subgraph}
Let $\trig{\tilde p}{\tilde x}{\tilde y}$ be the model triangle for 
${\trig p x y}_D$.
Then there is a short map 
\[G\:\Conv\trig{\tilde p}{\tilde x}{\tilde y}\to D\]
such that $\tilde p\mapsto p$, $\tilde x\mapsto x$, $\tilde y\mapsto y$, and $G$ is length-preserving on each side of $\trig{\tilde p}{\tilde x}{\tilde y}$.
In particular, $\Conv\trig{\tilde p}{\tilde x}{\tilde y}$ majorizes the triangle $[p x y]_D$ in $D$ under~$G$.
\end{subthm}
\end{thm} 

The proof is based on repeated application of the argument in the proof of the inheritance lemma (\ref{lem:inherit-angle}).

\parit{Proof.}
We prove the lemma for a polygonal line $\beta$;
the general case then follows by approximation.
Namely, since $\beta$ is concave 
it can be approximated by polygonal lines that are concave with respect to $p$, 
with their lengths converging to $\length \beta$. 
Passing to a partial limit we will obtain the needed map $G$.  

Suppose $\beta=x_0x_1\ldots x_n$ is a polygonal line with $x_0=x$ and $x_n\z=y$.
Consider a sequence of polygonal lines $\beta_i=x_0x_1\ldots x_{i-1}y_i$ such that $\dist{p}{y_i}{}=\dist{p}{y}{}$ and
$\beta_i$ has the same length as $\beta$;
that is, 
\[\dist{x_{i-1}}{y_i}{}=\dist{x_{i-1}}{x_{i}}{}+\dist{x_{i}}{x_{i+1}}{}+\ldots+\dist{x_{n-1}}{x_n}{}.\]

\begin{figure}[!ht]
\vskip-0mm
\centering
\includegraphics{mppics/pic-955}
\end{figure}

Clearly $\beta_n=\beta$.
Sequentially applying Alexandrov's lemma (\ref{lem:alex}), we get that each of the polygonal lines $\beta_{n-1}, \beta_{n-2},\ldots,\beta_1$ is concave with respect to $p$.

Let $D_i$ be the subgraph of $\beta_i$ with respect to $p$.
Applying the argument in the inheritance lemma (\ref{lem:inherit-angle}) gives a short map $G_i\:D_{i}\z\to D_{i+1}$ that maps $y_{i}\mapsto y_{i+1}$ and does not move $p$ and $x$ (in fact,  $G_i$ is the identity everywhere except on $\Conv\trig{p}{x_{i-1}}{y_i}$).
Thus the composition 
\[G_{n-1}\circ\ldots\circ G_1\: D_1\to D_n\]
is short.
The result follows since $D_1\iso\Conv\trig{\tilde p}{\tilde x}{\tilde y}$.
\qeds

\begin{thm}{Lemma}\label{lem:devel}\label{def:devel}
Let $\spc{X}$ be a metric space, 
$\gamma\:\II\to \spc{X}$ be a $1$-Lipschitz curve,
$p\in \spc{X}$,
and $\tilde p\in\EE^2$.
Then there exists a unique up to rotation curve
$\tilde \gamma\: \II\to \EE^2$, parametrized by arc-length, 
such that
$\dist{\tilde p}{\tilde \gamma(t)}{}\z=\dist{p}{\gamma(t)}{}$ for all $t$
and the direction of
$[\tilde p\tilde \gamma(t)]$ monotonically turns around $\tilde p$ counterclockwise as $t$ increases.
\end{thm}

If $p$, $\tilde p$, $\gamma$, and $\tilde \gamma$ are as above,
then $\tilde \gamma$ is called the \index{development}\emph{development} of $\gamma$ with respect to $p$; 
the point $\tilde p$ is called the \index{development!basepoint of a development}\emph{basepoint} of the development.
This notion appears in \cite{alexandrov-1957}
and an earlier form of it can be found in \cite{liberman}.

\parit{Proof.}
Consider the functions $\rho$, $\theta\:\II\to\RR$ defined as 
\begin{align*}
\rho(t)
&=\dist{p}{\gamma(t)}{},
&
\theta(t)
&=
\int\limits_{t_0}^{t}\frac{\sqrt{1-(\rho')^2}}{\rho},
\end{align*}
where $t_0\in\II$ is a fixed number and $\int$ denotes Lebesgue integral.
Since $\gamma$ is $1$-Lipschitz, so is $\rho(t)$, and thus the function $\theta$ is defined and nondecreasing.

It is straightforward to check that $(\rho,\theta)$ uniquely describe $\tilde \gamma$ in polar coordinates on $\EE^2$ with center at $\tilde p$.
\qeds

\begin{thm}[!]{Exercise}\label{ex:devel-comp-CAT}
Show that a geodesic space $\spc{U}$ is $\CAT(0)$ if and only if the development of any geodesic with respect to any point is concave.
\end{thm}

\begin{thm}{Lemma}\label{lem:majorize-triangle}
Let $\trig{p}{x}{y}$ be a triangle in a geodesic $\CAT(0)$ space $\spc{U}$,
and let $\tilde \gamma$ be the development of $[x y]$ with respect to $p$, where $\tilde \gamma$ has basepoint $\tilde p$ and subgraph $D$.
Consider the map $H\:D\to\spc{U}$ that sends the point with parameter $(t,s)$ under the line-of-sight map for $\tilde \gamma$ with respect to $\tilde p$, to the point with the same parameter under the line-of-sight map $f$ for $[x y]$ with respect to $p$.
Then $H$ is  length-nonincreasing.
In particular, $D$ majorizes the triangle $\trig p x y$.
\end{thm}

\parit{Proof.}
Let $\gamma\:[0,T]\to \spc{U}$ be a unit-speed parametrization of $[xy]$; so, $T=\dist{x}{y}{}$.
Choose a partition 
\[0=t_0<t_1<\ldots<t_n=T,\]
and set $x_i=\gamma(t_i)$.
Construct a chain of model triangles  $\trig{\tilde p}{\tilde x_{i-1}}{\tilde x_i}\z=
\modtrig(p x_{i-1} {x_i})$, with $\tilde x_0=\tilde x$ and the direction of $[\tilde p\tilde x_i]$ turning counterclockwise as $i$ grows.
Let $D_n$ be the subgraph with respect to $\tilde p$ of the polygonal line $\tilde x_0\ldots \tilde x_n$.

Now we construct a map $H_n \: D_n\to\spc{U}$  that increases distances by at most $2\cdot\delta_n$, where $\delta_n=\max\{\,\dist{\tilde x_i}{\tilde x_{i-1}}{}\,\}$.
Suppose $w\in D_n$.
Then $w$ lies on or inside some triangle $\trig{\tilde p}{\tilde x_{i-1}}{\tilde x_i}$.
Define $H_n(w)$ by first mapping $w$ to a nearest point on $\trig{\tilde p}{\tilde x_{i-1}}{\tilde x_i}$ (choosing one if there are several), followed by the natural map to the triangle  $\trig {p}{x_{i-1}}{ x_i}$.

Since triangles in $\spc{U}$ are thin, the restriction of $H_n$ to each triangle $\trig{\tilde p}{\tilde x_{i-1}}{\tilde x_i}$ is short.
Then the triangle inequality implies that the restriction of $H_n$ to 
\[U_n=\bigcup_{1\le i\le n}\trig{\tilde p}{\tilde x_{i-1}}{\tilde x_i}\]
is short with respect to the length metric on $D_n$. 
Since nearest-point projection from $D_n$ to $U_n$ increases the $D_n$-distance between two points by at most $2\cdot\delta_n$, the map $H_n$ also increases the $D_n$-distance by at most $2\cdot\delta_n$. 

Consider converging sequences $v_n\to v$ and $w_n\to w$ such that $v_n,w_n\in D_n$ and therefore $v,w\in D$.
Note that 
\[\dist{H_n(v_n)}{H_n(w_n)}{} \le \dist{v_n}{w_n}{D_n} + 2\cdot\delta_n,\eqlbl{eq:|H(v)-H(w)|}\]
for each $n$.
Since $\delta_n\to 0$ and geodesics in $\spc{U}$ vary continuously with their endpoints (\ref{thm:alex-patch}), we have $H_n(v_n)\to 
H(v)$ and $H_n(w_n)\to H(w)$.
Therefore the left-hand side in \ref{eq:|H(v)-H(w)|} converges to $\dist{H(v)}{H(w)}{}$ and the right-hand side converges to $\dist{v}{w}{D}$, it follows that $H$ is short.
\qeds

\parit{Proof of \ref{thm:major} for triangles.}
Suppose $\alpha$ is a triangle, say $\trig p x y$.

Let $\tilde \gamma$ be the development of $[x y]$ with respect to $p$, where $\tilde \gamma$ has basepoint $\tilde p$ and subgraph $D$.
By \ref{ex:devel-comp-CAT}, $\tilde \gamma$ is concave.
By \ref{lem:majorize-subgraph},  there is a short map $G\:\Conv\modtrig(p x y)\to D$.
Further, by \ref{lem:majorize-triangle},  $D$ majorizes $\trig p x y$ under a majorizing map $H\:D\to\spc{U}$.
Clearly $H\circ G$ is a majorizing map for $\trig p x y$.
\qeds

\section{Polygons}\label{sec:polygons}

A closed piecewise geodesic curve in a metric space will be called a \index{polygon}\emph{polygon}.
Further, $[x_1\ldots x_n ]$ (for $n\ge 3$) will denote a polygon
with edges $[x_1x_2],\ldots,[x_{n-1}x_n]$, and $[x_n x_1]$.
For a subset $R$ of the ambient metric space,
we denote by $[x_1\ldots x_n ]_R$ a polygon in the length metric of $R$.

\begin{wrapfigure}{o}{40 mm}
\vskip-4mm
\centering
\includegraphics{mppics/pic-960}
\vskip0mm
\end{wrapfigure}

\parit{Proof of \ref{thm:major} for polygons.}
We claim that any closed $n$-gon $[x_1x_2 \ldots x_n]$ in a $\CAT(0)$ space $\spc{U}$ is majorized by a convex polygonal region
\[R_n=\Conv[\tilde x_1\tilde x_2\ldots\tilde x_n],\]
under a map $F_n$ that sends $\tilde x_i\mapsto x_i$ for each $i$.

The base case $n=3$ is proved above.
Assume the statement is true for $(n-1)$-gons, $n\ge 4$.  
Then  $[x_1 x_2 \ldots x_{n-1}]$  is majorized by a convex polygonal region
\[R_{n-1}=\Conv[\tilde x_1 \tilde x_2,\ldots, \tilde x_{n-1}],\]
in $\EE^2$ under a map $F_{n-1}$ satisfying $F_{n-1}(\tilde x_i)=x_i$ for all $i$.
Take $\dot x_n\in\EE^2$ such that $\trig{\tilde x_1}{\tilde x_{n-1}}{\dot x_n}=\modtrig(x_1 x_{n-1} x_n)$
and this triangle lies on the other side of $[\tilde x_1\tilde x_{n-1}]$ from $R_{n-1}$.
Let $\dot R\z=\Conv\trig{\tilde x_1}{\tilde x_{n-1}}{\dot x_n}$,
and $\dot F\:\dot R\to \spc{U}$ be a majorizing map for $\trig { x_1}{x_{n-1}}{ x_n}$ as provided above.

Set 
$R= R_{n-1}\cup \dot R$, where $R$ carries its length metric.
Since $F_{n-1}$ and $\dot F$ agree on $[\tilde x_1 \tilde x_{n-1}]$, we may define $F\:R\to\spc{U}$ by
\[
F(x)=
\begin{cases}
F_{n-1}(x),\quad & x\in R_{n-1},\\
\dot F(x),\quad & x\in \dot R.\\
\end{cases}
\]
Then $F$ is length-nonincreasing and is a majorizing map for $[x_1 x_2 \z\ldots x_n ]$ (as in Definition~\ref{def:majorize}).

If $R$ is a convex subset of $\EE^2$, we are done. 

If $R$ is not convex, the total internal angle of $R$ at $\tilde x_1$ or $ \tilde x_{n-1} $ or both is larger than $\pi$.
By relabeling we may suppose that this holds for~$\tilde x_{n-1}$.

The region $R$ is obtained by gluing $R_{n-1}$ to $\dot R$ by $[x_1x_{n-1}]$.
Thus, by Reshetnyak gluing (\ref{thm:gluing}), $R$ carrying its length metric is a $\CAT(0)$ space.
Moreover $[\tilde x_{n-2}\tilde x_{n-1}]\cup[\tilde x_{n-1} \dot x_n]$ is a geodesic of $R$.
Thus $[\tilde x_1 \tilde x_2 \ldots \tilde x_{n-2} \dot x_n]_R$ is a closed $(n-1)$-gon in $R$, to which the induction hypothesis applies. The resulting short map from a convex region in $ \EE^2$ to~$R$, followed by $F$,  is the desired majorizing map.
\qeds

For a polygon $[p_1\ldots p_n]$, the values $\theta_i=\pi-\mangle\hinge{p_i}{p_{i-1}}{p_{i+1}}$ for all $i\pmod n$ are called \emph{external angles} of the polygon.
The following exercise is a generalization of Fenchel's theorem.

\begin{thm}{Exercise}\label{ex:fenchel}
Show that the sum of external angles of any polygon in a complete length $\CAT(0)$ space cannot be smaller than $2\cdot\pi$. 
\end{thm}

The following exercise is a version of the Fáry--Milnor theorem for $\CAT(0)$ spaces.

\begin{thm}{Very advanced exercise}\label{ex:FM}
Suppose that a simple polygon $\beta$ in a complete length $\CAT(0)$ space does not bound an embedded disc.
Show that the sum of external angles of $\beta$ cannot be smaller than $4\cdot\pi$.

Give an example of such a polygon $\beta$ with the sum of external angles exactly $4\cdot\pi$.
\end{thm}
  
\begin{thm}{Arm lemma}\label{lem:arm+}
Let $P=[x_1\ldots x_n]$ be a polygon in a geodesic $\CAT(0)$ space $\spc{U}$.
Suppose $\tilde P=[\tilde x_1\ldots \tilde x_{n}]$ is a convex  polygon in $\EE^2$
such that 
\[
\dist{\tilde x_i}{\tilde x_{i-1}}{\EE^2}
=
\dist{x_i}{x_{i-1}}{\spc{U}}
\quad \text{and}\quad 
\mangle\hinge{x_i}{x_{i-1}}{x_{i+1}}\ge\mangle\hinge{\tilde x_i}{\tilde x_{i-1}}{\tilde x_{i+1}}
\eqlbl{eq:arm}
\]
for all $i$.
Then 
$\dist{x_1}{x_{n}}{\spc{U}}\ge \dist{\tilde x_1}{\tilde x_n}{\EE^2}$.
\end{thm}

If $\spc{U}$ is a Euclidean space, then the statement follows from the classical arm lemma;
see \cite{sabitov} and the references therein.

\begin{thm}{Exercise}\label{ex:arm-lemma}
Reduce \ref{lem:arm+} to the classical arm lemma.
\end{thm}

\section{General case}\label{seg:major-general}

If the space is proper, then the general case of the Majorization Theorem follows by applying polygonal case to inscribed polygonal lines and passing to the limit.

\section{Notes}

The argument in Section \ref{seg:major-general} makes essential use of the assumption that the space is proper.
Nevertheless, the theorem holds for any geodesic $\CAT(0)$ space \cite[9.56]{alexander-kapovitch-petrunin-2025}.

The majorization theorem can be generalized to $\CAT(\pm1)$ spaces;
in the $\CAT(1)$ case one has to assume that the closed curve has length at most $2\cdot\pi$.

This theorem was proved by Yuri Reshetnyak \cite{reshetnyak:major};
our proof uses a trick that we learned from the lectures of Werner Ballmann \cite{ballmann-1995}.
For complete spaces, another proof can be built  on the following closely related theorem.
It was discovered by Urs Lang and Viktor Schroeder \cite{lang-schroeder};
the third author proved it slightly earlier, but did not publish the proof for quite a while \cite{alexander-kapovitch-petrunin-2011,alexander-kapovitch-petrunin-2025}.

\begin{thm}{Generalized Kirszbraun's theorem}
\label{thm:kirsz+}
Let $\spc{U}$ be a complete length $\CAT(0)$ space, 
let $Q$ be an arbitrary subset of the Euclidean space~$\EE^m$.
Suppose $f\: Q\to\spc{U}$ is a short map.
Then $f\:Q\to\spc{U}$ can be extended to a short map 
$F\:\EE^m\to \spc{U}$.
\end{thm}

\begin{thm}{Open problem}
Consider  a closed rectifiable curve $\alpha$ in a $\CAT(0)$ space $\spc{U}$.
Note that if $\alpha$ is a geodesic triangle or it bounds an isometric copy of a convex plane figure in $\spc{U}$, then $\alpha$ has a unique (up to congruence) majorizing convex figure.

Does the converse hold?
\end{thm}

It is expected that the optimal isoperimetric inequality in the style of Frederick Almgren \cite{almgren} holds in \emph{Hadamard manifolds} (that is, $\CAT(0)$ Riemannian manifolds).
Namely, \textit{for any integers $n>m>0$ and real $r>0$, any $m$-dimensional surface with $m$-volume $a\cdot r^m$ in an $n$-dimensional Hadamard manifold can be filled with $m+1$-volume $v\cdot r^{m+1}$, where $a$ is the volume of the unit $m$-sphere, and $v$ is the volume of the unit ball in $\EE^{m+1}$}.
\begin{itemize}
 \item For $m=1$ it follows from the majorization theorem.
 \item For $m=3$ and $n=4$ it was proved by Christopher Croke \cite{croke}.
 \item For $m=2$ and $n=3$ it was proved later by Bruce Kleiner \cite{kleiner-1992}.
 \item For $m=2$ it was proved by Felix Schulze \cite{schulze}.
\end{itemize}
In all other cases, the problem remains open despite a number of unsuccessful attempts.

%%!TEX root = the-globalization.tex
\chapter{Globalization}\label{chapter:globalization}

This lecture gives a sufficient condition for locally $\CAT(0)$ spaces to be globally $\CAT(0)$.

\section{Locally CAT spaces}

We say that a space $\spc{U}$ is \index{locally $\CAT(\kappa)$ space}\emph{locally $\CAT(0)$} ({}\emph{locally $\CAT(1)$}) if
a small closed ball centered at any point $p$ in $\spc{U}$ is $\CAT(0)$ ($\CAT(1)$, respectively).

For example, the circle $\mathbb{S}^1=\RR/\ZZ$ is locally isometric to $\RR$, and so $\mathbb{S}^1$ is locally $\CAT(0)$.
On the other hand, $\mathbb{S}^1$ is not $\CAT(0)$, since closed local geodesics in $\mathbb{S}^1$ are not length-minimizing, so $\mathbb{S}^1$ does not satisfy \ref{cor:loc-geod-are-min}.

%\begin{thm}{Observation}
%Let $\spc{U}$ be a proper geodesic space.
%Then it is locally $\CAT(0)$ (or locally $\CAT(1)$) if and only ifeach point $p\in \spc{U}$ admits an open neighborhood $\Omega$ that is geodesic and such that any triangle in $\Omega$ is thin (or spherically thin, respectively).
%\end{thm}

\section{Space of local geodesic paths}\label{sec:geod-space}

In this section, we will study the behavior of local geodesics in locally $\CAT(\kappa)$  spaces.  
The results will be used in the proof of the globalization theorem (\ref{thm:hadamard-cartan}).

Recall that a \index{path}\emph{path} is a curve parametrized by the unit interval $[0,1]$,
and a \index{geodesic!local geodesic}\emph{local geodesic path} is a constant-speed parameterization of a local geodesic by $[0,1]$.
The space of paths in a metric space $\spc{U}$ comes with the natural metric
\[\dist{\alpha}{\beta}{}
=
\sup\set{\dist{\alpha(t)}{\beta(t)}{\spc{U}}}{t\in[0,1]}.
\eqlbl{eq:dist-between-paths}
\]

By \eqref{ex:convex-dist} and \eqref{ex:CAT-geodesic} geodesic paths in a complete geodesic $\CAT(0)$ space depend uniquely and continuously on their end points.
The same is true in $\CAT(1)$ spaces for geodesic paths shorter than $\pi$.

\begin{thm}{Proposition}\label{prop:geo-complete}
Let $\spc{U}$ be a proper geodesic, locally $\CAT(\kappa)$ space.

Assume $\gamma_n\:[0,1]\to\spc{U}$ is a sequence of local geodesic paths converging to a path $\gamma_\infty\:[0,1]\to\spc{U}$.
Then $\gamma_\infty$ is a local geodesic path.
Moreover 
\[\length\gamma_n\to\length\gamma_\infty\]
as $n\to\infty$.
\end{thm}

\parit{Proof.}
Assume $\kappa=0$.
Fix $t\in[0,1]$.  
Let $R>0$ be sufficiently small, so that $\cBall[\gamma_\infty(t),R]$ forms a proper geodesic $\CAT(0)$ space.
%"value" looks weird here. value usually means value of a function. V.

Assume that  a local geodesic $\sigma$  is shorter than $R/2$ and intersects the ball $\oBall(\gamma_\infty(t),R/2)$.
Then $\sigma$ cannot leave the ball $\cBall[\gamma_\infty(t),R]$.
By~\ref{cor:loc-geod-are-min}, $\sigma$ is a geodesic.  
In particular, for all sufficiently large $n$, any arc of $\gamma_n$ of length $R/2$ or less containing $\gamma_n(t)$ is a geodesic.

Since $\spc{B}=\cBall[\gamma_\infty(t),R]$ is a complete geodesic $\CAT(0)$ space, by \ref{ex:CAT-geodesic} and \ref{ex:convex-dist} geodesic segments in $\spc{B}$ depend uniquely and continuously on their endpoint pairs.
Thus there is a subinterval $\II$ of $[0,1]$,
that  contains a neighborhood of $t$ in $[0,1]$
and such that the arc $\gamma_n|_\II$ is minimizing for all large~$n$.
It follows that $\gamma_\infty|_\II$ is a geodesic,
and therefore $\gamma_\infty$ is a local geodesic.

The $\CAT(1)$ case is done in the same way, but one has to assume in addition that $R<\pi$.
\qeds

The following lemma allows a local geodesic path to be moved continuously so that its endpoints follow given trajectories.

\begin{thm}{Patchwork along a geodesic}
\label{lem:patch}
Let $\spc{U}$ be a proper geodesic, locally $\CAT(0)$ space, 
and $\gamma\:[0,1]\to\spc{U}$ be a locally geodesic path.

Then there is a proper geodesic  $\CAT(0)$ space   $\spc{N}$,
an open set $\hat\Omega\z\subset \spc{N}$,
and a geodesic path $\hat\gamma\:[0,1]\to\hat\Omega$,
such that there is an open locally distance-preserving map 
$\Phi\:\hat\Omega\looparrowright\spc{U}$ satisfying
$\Phi\circ\hat\gamma=\gamma$.

If $\length\gamma<\pi$,
then the same holds in the $\CAT(1)$ case.
Namely, we assume that $\spc{U}$ is a proper geodesic, locally $\CAT(1)$ space and construct a proper geodesic $\CAT(1)$ space $\spc{N}$ with the same property as above.
\end{thm}

\parit{Proof.} 
Fix $r>0$ so that for each $t\in[0,1]$,
the closed ball
$\cBall[\gamma(t),r]$ forms a proper geodesic $\CAT(0)$ space.

\begin{figure}[ht!]
\vskip-0mm
\centering
\includegraphics{mppics/pic-910}
\end{figure}

Choose a partition $0\z=t_0<t_1<\ldots<t_n\z=1$ such that
\[\oBall(\gamma(t_i),r)\supset \gamma([t_{i-1},t_i])\]
for each $i$.
Set $\spc{B}_i=\cBall[\gamma(t_i),r]$.
We can assume that $\spc{B}_{i-1}\cap \spc{B}_{i+1}\z\subset \spc{B}_{i}$ if both sets are defined.

Consider the disjoint union $\bigsqcup_i\spc{B}_i=\set{(i,x)}{x\in\spc{B}_i}$ with the minimal equivalence relation $\sim$ such that $(i,x)\sim(i-1,x)$ for all~$i$.
Let  $\spc{N}$ be the space obtained by gluing the $\spc{B}_i$ along~$\sim$.

Note that $A_i=\spc{B}_i\cap\spc{B}_{i-1}$ is convex in $\spc{B}_i$ and in $\spc{B}_{i-1}$.
Applying the Reshetnyak gluing theorem (\ref{thm:gluing}) $n$ times, 
we conclude that $\spc{N}$ is a proper geodesic $\CAT(0)$ space.

For $t\in[t_{i-1},t_i]$, define $\hat\gamma(t)$ as the equivalence class of $(i,\gamma(t))$ in~$\spc{N}$.
Let $\hat\Omega$ be the $\eps$-neighborhood of $\hat\gamma$ in $\spc{N}$, where $\eps>0$ is chosen so that $\oBall(\gamma(t),\eps)\subset\spc{B}_i$ for all $t\in[t_{i-1},t_i]$.

Define $\Phi\:\hat\Omega\to\spc{U}$
by sending the equivalence class of $(i,x)$ to~$x$.
It is straightforward to check that $\Phi$, 
$\hat\gamma$, and $\hat\Omega\subset\spc{N}$ satisfy the conclusion of  the lemma.

The $\CAT(1)$ case is proved in the same way.
\qeds

Recall that local geodesics are geodesics in any $\CAT(0)$ space; see \ref{cor:loc-geod-are-min}.
Using it with \ref{lem:patch} and the uniqueness of geodesics (\ref{ex:CAT-geodesic}), we get the following.

\begin{thm}{Corollary}\label{cor:discrete-paths}
If $\spc{U}$ is a proper geodesic, locally $\CAT(0)$ space, then for any pair of points $p,q\in\spc{U}$, the space of all local geodesic paths from $p$ to $q$ is discrete;
that is, for any local geodesic path $\gamma$ connecting $p$ to $q$, there is $\eps>0$ such that for any other local geodesic path $\delta$ from $p$ to $q$ we have
$\dist{\gamma(t)}{\delta(t)}{\spc{U}}>\eps$ for some $t\in[0,1]$.

Analogously, if $\spc{U}$ is a proper geodesic, locally $\CAT(1)$ space, then for any pair of points $p,q\in\spc{U}$,  the space of all local geodesic paths shorter than $\pi$ from $p$ to $q$ is discrete.
\end{thm}

\begin{thm}{Corollary}\label{cor:path-geod}
If $\spc{U}$ is a proper geodesic, locally $\CAT(0)$ space, then 
for any path $\alpha$ there is a choice of a local geodesic path $\gamma_\alpha$  connecting the ends of $\alpha$ such that the map $\alpha\mapsto\gamma_\alpha$ is continuous.
Moreover, if $\alpha$ is a local geodesic path, then $\gamma_\alpha=\alpha$.

Analogously, if $\spc{U}$ is a proper geodesic, locally $\CAT(1)$ space, then for any path $\alpha$ shorter than $\pi$, there is a choice of a local geodesic path $\gamma_\alpha$ shorter than $\pi$ connecting the ends of $\alpha$ such that the map $\alpha\mapsto\gamma_\alpha$ is continuous.
Moreover, if $\alpha$ is a local geodesic path, then $\gamma_\alpha=\alpha$.
\end{thm}

\parit{Proof of \ref{cor:path-geod}.} 
We only present the proof in the $\CAT(0)$ case;
the $\CAT(1)$ case is analogous.

Consider the maximal interval $\II\subset[0,1]$ containing $0$
such that there is a continuous one-parameter family of 
local geodesic paths $\gamma_t$ for $t\in \II$ connecting $\alpha(0)$ to $\alpha(t)$ (that is, $\gamma_t(0)=\alpha(0)$ and $\gamma_t(1)=\alpha(t)$ for any~$t$) such that $\gamma_0$ is the constant path at $\alpha(0)$.

By \ref{prop:geo-complete}, $\II$ is closed,
so we may assume $\II=[0,s]$ for some $s\in [0,1]$.

Applying  patchwork (\ref{lem:patch}) to  $\gamma_{s}$, 
we find that $\II$ is also open in $[0,1]$. 
Hence $\II=[0,1]$.
Set $\gamma_\alpha=\gamma_1$.

By construction,  if $\alpha$ is a local geodesic path, then $\gamma_\alpha=\alpha$. 

Moreover, from \ref{cor:discrete-paths},
the construction $\alpha\mapsto \gamma_\alpha$ produces close results for sufficiently close paths in the metric defined by \ref{eq:dist-between-paths};
that is, the map  $\alpha\mapsto \gamma_\alpha$ is continuous.
\qeds

Given a path $\alpha\:[0,1]\to\spc{U}$,
we denote by $\bar\alpha$ the same path traveled in the opposite direction;
that is,
\[\bar\alpha(t)=\alpha(1-t).\]
The \index{product of paths}\emph{product} of two paths  will be denoted with ``$*$'';
if two paths $\alpha$ and $\beta$ connect the same pair of points, then the product $\bar\alpha*\beta$ is a closed curve.

\begin{thm}[!]{Exercise}\label{ex:null-homotopic}
Assume $\spc{U}$ is a proper geodesic locally $\CAT(1)$ space. 
Let $\alpha\mapsto\gamma_\alpha$ be the construction provided by Corollary~\ref{cor:path-geod}.

Assume that $\alpha$ and $\beta$ are two paths connecting the same pair of points in $\spc{U}$, where 
each is shorter than $\pi$ 
and the product  
$\bar\alpha*\beta$ is null-homotopic in the class of closed curves shorter than $2\cdot\pi$.
Show that $\gamma_\alpha=\gamma_\beta$.
\end{thm}

\section{Globalization}\label{sec:Hadamard--Cartan}

\begin{thm}{Globalization theorem}
\label{thm:hadamard-cartan}
If a proper geodesic locally $\CAT(0)$ space is simply connected, then it 
is $\CAT(0)$.

Analogously, if $\spc{U}$ is a proper geodesic, locally $\CAT(1)$ space
such that any closed curve $\gamma\:\mathbb{S}^1\to \spc{U}$ shorter than $2\cdot\pi$ is null-homotopic in the class of closed curves shorter than $2\cdot\pi$, then $\spc{U}$ is $\CAT(1)$.
\end{thm}

The surface in the diagram is an example of a simply connected space that  is locally $\CAT(1)$ but not $\CAT(1)$;
\begin{figure}[ht!]
\vskip0mm
\centering
\includegraphics{mppics/pic-930}
\end{figure}
we assume that the marked curve is shorter than $2\cdot\pi$, but
to contract it, one has to increase its length to $2\cdot\pi$ or more.
In particular, this surface does not satisfy the assumptions of the globalization theorem.

The proof of the globalization theorem relies on the following theorem, 
which is essentially  \cite[Satz 9]{alexandrov-1957}.  

\begin{thm}{Patchwork globalization theorem}\label{thm:alex-patch}
A proper geodesic, locally $\CAT(0)$ space $\spc{U}$ is $\CAT(0)$
if and only if all pairs of points in $\spc{U}$  are joined by unique geodesics, and these geodesics depend continuously on their endpoint pairs.

Analogously, a proper geodesic, locally $\CAT(1)$ space $\spc{U}$ is $\CAT(1)$ 
if and only if all pairs of points in $\spc{U}$ at distance less than $\pi$ are joined by unique geodesics, and these geodesics depend continuously on their endpoint pairs.
\end{thm}

The proof uses a thin-triangle decomposition with the inheritance lemma (\ref{lem:inherit-angle}) and the line-of-sight map (\ref{def:sight}).

\parit{Proof.}
Note that the implication ``only if'' follows from \ref{ex:CAT-geodesic} and \ref{ex:convex-dist}; it remains to prove the ``if'' part.

Fix a triangle $\trig p x y$  in $\spc{U}$. 
We need to show that $\trig p x y$ is thin.

By the assumptions, the line-of-sight map  $(t,s)\mapsto\gamma_t(s)$ from $p$ to   $[x y]$ is uniquely defined and continuous.

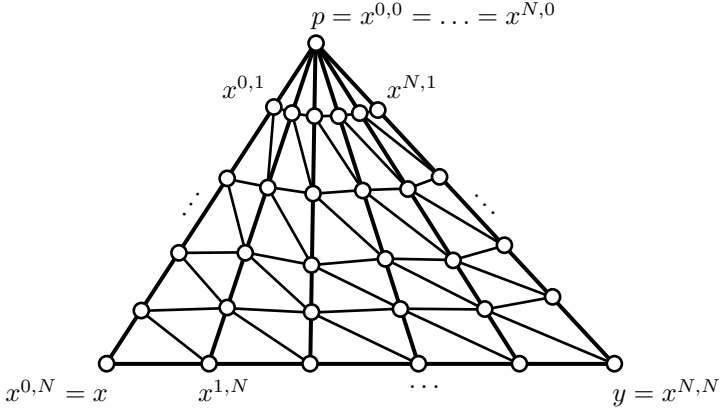
\begin{figure}[ht!]
\vskip0mm
\centering
\includegraphics{mppics/pic-950}
\end{figure}

Fix  a partition \[0\z=t_0\z<t_1\z<\z\ldots\z<t_N=1,\]
and set $x_{i,j}=\gamma_{t_i}(t_j)$.
Since the line-of-sight map is continuous and $\spc{U}$ is locally $\CAT(0)$, we may assume that the triangles 
\[\trig{x_{i,j}\,}{x_{i,j+1}\,}{x_{i+1,j+1}}\quad\text{and}\quad\trig{x_{i,j}\,}{x_{i+1,j}\,}{x_{i+1,j+1}}\]
are thin for each pair $i$,~$j$.

Now we show that the thin property propagates to $\trig p x y$ by repeated application of the inheritance lemma (\ref{lem:inherit-angle}):
\begin{itemize}
\item 
For fixed $i$, 
sequentially applying the lemma shows that the triangles 
$\trig{p\,}{x_{i,1}\,}{x_{i-1,2}}$,
$\trig{p\,}{x_{i,2}\,}{x_{i-1,2}}$,
$\trig{p\,}{x_{i,2}\,}{x_{i-1,3}}$,
and so on are thin. 
\end{itemize}
In particular, for each $i$, the long triangle $\trig{p\,}{x_{i,N}\,}{x_{i-1,N}}$ is thin.
\begin{itemize} 
\item 
By the same lemma the  triangles $\trig{p\,}{x_{0,N}\,}{x_{2,N}}$, $\trig{p\,}{x_{0,N}\,}{x_{3,N}}$, and so on, are thin.
\end{itemize}
In particular, $\trig p x y=\trig{p\,}{x_{0,N}\,}{x_{N,N}}$ is thin.
\qeds

\parit{Proof of the globalization theorem; $\CAT(0)$ case.}
Let $\spc{U}$ be a proper geodesic, locally $\CAT(0)$ space that is simply connected.
Given a path $\alpha$ in $\spc{U}$, 
denote by $\gamma_\alpha$ the local geodesic path provided by \ref{cor:path-geod}.
Since the map $\alpha\mapsto\gamma_\alpha$ is continuous, by \ref{cor:discrete-paths}
we have $\gamma_\alpha=\gamma_\beta$ for any pair of  paths $\alpha$ and $\beta$  homotopic relative to the ends.

Since $\spc{U}$ is simply connected, any pair of paths with common ends are homotopic.  In particular, if $\alpha$ and $\beta$ are local geodesics from $p$ to $q$, then $\alpha =\gamma_\alpha=\gamma_\beta=\beta$ by Corollary \ref{cor:path-geod}.
It follows that any two points $p,q\in\spc{U}$ are joined by a unique local geodesic that depends continuously on $(p,q)$.

Since $\spc{U}$ is geodesic, it remains to apply the patchwork globalization theorem (\ref{thm:alex-patch}).

\parit{$\CAT(1)$ case.}
The proof goes along the same lines, 
but one needs to use Exercise~\ref{ex:null-homotopic}. \qeds

\begin{thm}{Corollary}\label{cor:closed-geod-cat} 
Any compact geodesic, locally $\CAT(0)$ space that contains no closed local geodesics is $\CAT(0)$. 
 
Analogously, any compact geodesic, locally $\CAT(1)$ space that contains no closed local geodesics shorter than $2\cdot\pi$ is $\CAT(1)$.
\end{thm}

\parit{Proof.}
By the globalization theorem (\ref{thm:hadamard-cartan}), we need to show that the space is simply connected.
Assume the contrary. 
Fix a nontrivial homotopy class of closed curves.

Denote by $\ell$ the exact lower bound for the lengths of curves in the class.
Note that $\ell>0$;
otherwise, there would be a closed noncontractible curve in a $\CAT(0)$ neighborhood of some point, contradicting \ref{ex:contractible}.

Since the space is compact, the class contains a length-minimizing curve, 
which must be a closed local geodesic. 

The $\CAT(1)$ case is analogous, one only has to consider a homotopy in the class of closed curves shorter than $2\cdot\pi$.
\qeds

\begin{thm}{Exercise}\label{ex:geod-circle}
Prove that any compact geodesic, locally $\CAT(0)$ space $\spc{X}$ that is not $\CAT(0)$ contains a \index{geodesic!circle}\emph{geodesic circle};
that is, a simple closed curve $\gamma$ such that 
for any two points $p,q\in\gamma$, one of the arcs of $\gamma$ with endpoints $p$ and $q$ is a  geodesic.

Formulate and prove the analogous statement for $\CAT(1)$ spaces.
\end{thm}

\begin{thm}{Advanced exercise}\label{ex:branching-cover} 
Let $\spc{U}$ be a proper geodesic $\CAT(0)$ space.
Assume $\tilde{\spc{U}}\to \spc{U}$ is a metric double cover branching along a geodesic~$\gamma$.
(Formally speaking, $\tilde{\spc{U}}$ is the completion of a double cover of the complement $\spc{U}\setminus \gamma$.
For example, 3-dimensional Euclidean space admits a double cover branching along a line.)

Show that $\tilde{\spc{U}}$ is $\CAT(0)$.
\end{thm}

\section{Notes}

The lemma about patchwork along a geodesic (\ref{lem:patch}) and its proof were suggested to us by Alexander Lytchak.  
This statement was originally proved by the first author and Richard Bishop \cite{alexander-bishop-1990} using a different method.

In the globalization theorem (\ref{thm:hadamard-cartan}), properness can be weakened to completeness \cite[see][and the references therein]{alexander-kapovitch-petrunin-2025}.
The original formulation of the globalization theorem, or 
Hadamard--Cartan theorem, states that \textit{if $M$ is a complete Riemannian manifold with sectional curvature at most $0$,  
then the exponential map at any point $p\in M$ is a covering};
in particular, it implies that \textit{the universal cover of $M$ is diffeomorphic to the Euclidean space of the same dimension}.

In this generality, this theorem first appeared in the lectures of Elie Cartan~\cite{cartan}.
This theorem was proved for surfaces in Euclidean $3$-space 
by Hans von Mangoldt \cite{mangoldt}
and a few years later independently for 2-dimensional Riemannian manifolds by Jacques Hadamard \cite{hadamard}.

Formulations for metric spaces of different generality were proved by 
Herbert Busemann \cite{busemann-CBA},
Willi Rinow \cite{rinow},
Mikhael Gromov \cite[p.~119]{gromov-1987}. 
A detailed proof of Gromov's statement for proper spaces was given by Werner Ballmann \cite{ballmann-1995},
and for complete spaces by the first author and Richard Bishop \cite{alexander-bishop-1990}.

For proper $\CAT(1)$ spaces, the globalization theorem was proved by Brian Bowditch~\cite{bowditch-1995}.

{\sloppy

The globalization theorem holds for complete (not necessarily proper)  length spaces \cite{alexander-kapovitch-petrunin-2025}.

}

The patchwork globalization (\ref{thm:alex-patch}) was proved by Alexandrov \cite[Satz 9]{alexandrov-1957}.
For proper spaces one can remove the continuous dependence from the formulation; it follows from uniqueness.
For complete spaces, the latter is not true \cite[Chapter I, Exercise 3.14]{bridson-haefliger}.

%%!TEX root = the-polyhedral.tex
\chapter{Polyhedral spaces}\label{chap:poly}

This lecture introduces LEGO-like rules for assembling $\CAT(0)$ spaces out of Euclidean cubes and then uses them to construct exotic aspherical manifolds.

\section{Products, cones, and suspension}
\label{sec:Products and cones}

Products, cones, and suspension are defined in Section~\ref{sec:constructions}.

\begin{thm}{Proposition}\label{ex:product-CAT}
Let $\spc{U}$ and $\spc{V}$ be $\CAT(0)$ spaces.
Then the product space $\spc{U}\times\spc{V}$ is $\CAT(0)$.
\end{thm}

\parit{Proof.}
Fix a quadruple in $\spc{U}\times \spc{V}$:
\begin{align*}
p&=(p_1,p_2),
&
q&=(q_1,q_2), 
&
x&=(x_1,x_2),
&
y&=(y_1,y_2).
\end{align*}
For the quadruple $p_1,q_1,x_1,y_1$ in $\spc{U}$,
construct two model triangles $\trig{\tilde p_1}{\tilde x_1}{\tilde y_1}=\modtrig(p_1x_1y_1)_{\EE^2}$ 
and $\trig{\tilde q_1}{\tilde x_1}{\tilde y_1}=\modtrig(q_1x_1y_1)_{\EE^2}$.
Similarly, for the quadruple $p_2,q_2,x_2,y_2$ in $\spc{V}$
construct two model triangles $\trig{\tilde p_2}{\tilde x_2}{\tilde y_2}$ and $\trig{\tilde q_2}{\tilde x_2}{\tilde y_2}$.

Consider four points in $\EE^4=\EE^2\times\EE^2$ 
\begin{align*}
\tilde p&=(\tilde p_1,\tilde p_2),
&
\tilde q&=(\tilde q_1,\tilde q_2),
&
\tilde x&=(\tilde x_1,\tilde x_2),
&
\tilde y&=(\tilde y_1,\tilde y_2).
\end{align*}
Note that the triangles $\trig{\tilde p}{\tilde x}{\tilde y}$ and $\trig{\tilde q}{\tilde x}{\tilde y}$ in $\EE^4$ are isometric to the model triangles 
$\modtrig(pxy)_{\EE^2}$ and $\modtrig(qxy)_{\EE^2}$.
Note however that they need not lie in the same plane.

If $\tilde z=(\tilde z_1,\tilde z_2)\in [\tilde x\tilde y]$, then $\tilde z_1\in [\tilde x_1\tilde y_1]$ and $\tilde z_2\in [\tilde x_2\tilde y_2]$ and
\begin{align*}
\dist[2]{\tilde z}{\tilde p}{\EE^4}&=\dist[2]{\tilde z_1}{\tilde p_1}{\EE^2}+\dist[2]{\tilde z_2}{\tilde p_2}{\EE^2},
\\
\dist[2]{\tilde z}{\tilde q}{\EE^4}&=\dist[2]{\tilde z_1}{\tilde q_1}{\EE^2}+\dist[2]{\tilde z_2}{\tilde q_2}{\EE^2},
\\
\dist[2]{p}{q}{\spc{U}\times\spc{V}}&=\dist[2]{p_1}{q_1}{\spc{U}}+\dist[2]{p_2}{q_2}{\spc{V}}.
\end{align*}
Therefore $\CAT(0)$ comparison for the quadruples $p_1,q_1,x_1,y_1$ in $\spc{U}$
and 
$p_2,q_2,x_2,y_2$ in $\spc{V}$ implies that
\[\dist{p}{q}{\spc{U}\times\spc{V}}
\le
\dist{\tilde p}{\tilde z}{\EE^4}+\dist{\tilde z}{\tilde q}{\EE^4}.\]
This gives $\CAT(0)$ comparison for the quadruple $p,q,x,y$ in $\spc{U}\z\times \spc{V}$.
\qeds

\begin{thm}{Proposition}\label{ex:cone+susp}
Let $\spc{U}$ be a metric space.
Then $\Cone\spc{U}$ is $\CAT(0)$ if and only if $\spc{U}$ is $\CAT(1)$.
\end{thm}

\parit{Proof; if part.}
Given a point $x\in \Cone\spc{U}$, denote by $x'$ its projection to $\spc{U}$
and by $|x|$ the distance from $x$ to the tip of the cone;
if $x$ is the tip, then $|x|=0$ and we can take any point of $\spc{U}$ as~$x'$.

Let $p$, $q$, $x$, $y$
be a quadruple in $\Cone\spc{U}$.
Assume that the spherical model triangles $\trig{\tilde p'}{\tilde x'}{\tilde y'}_{\SSS^2}=\modtrig(p'x'y')_{\SSS^2}$ and $\trig{\tilde q'}{\tilde x'}{\tilde y'}_{\SSS^2}=\modtrig(q'x'y')_{\SSS^2}$ are defined.
Consider the following points in $\EE^3=\Cone\SSS^2$: 
\begin{align*}
\tilde p&=|p|\cdot\tilde p',
&
\tilde q&=|q|\cdot\tilde q',
&
\tilde x&=|x|\cdot\tilde x',
&
\tilde y&=|y|\cdot\tilde y'.
\end{align*}

Note that
$\trig{\tilde p}{\tilde x}{\tilde y}_{\EE^3}\iso\modtrig(pxy)_{\EE^2}$
and
$\trig{\tilde q}{\tilde x}{\tilde y}_{\EE^3}\iso\modtrig(qxy)_{\EE^2}$.
Further, note that if $\tilde z\in [\tilde x\tilde y]_{\EE^3}$, then
$\tilde z'=\tilde z/|\tilde z|$ lies on the geodesic $[\tilde x'\tilde y']_{\SSS^2}$.
Therefore, $\CAT(1)$ comparison for $\dist{p'}{q'}{}$ with $\tilde z'\in[\tilde x'\tilde y']_{\SSS^2}$ implies
$\CAT(0)$ comparison for $\dist{p}{q}{}$ with $\tilde z\in[\tilde x\tilde y]_{\EE^3}$.

If at least one of the model triangles $\modtrig(p'x'y')_{\SSS^2}$ and $\modtrig(q'x'y')_{\SSS^2}$ is undefined,
then the statement follows from the triangle inequalities 
\begin{align*}
|p'-x'|_{\spc{U}}+|q'-x'|_{\spc{U}}
&\ge |p'-q'|_{\spc{U}}
\\
|p'-y'|_{\spc{U}}+|q'-y'|_{\spc{U}}
&\ge |p'-q'|_{\spc{U}}
\end{align*}
The details are left as an exercise.

\parit{Only-if part.}
Suppose that $\tilde p'$, $\tilde q'$, $\tilde x'$, $\tilde y'$ are defined as above.
Assume all these points lie in a half-space of $\EE^3=\Cone\SSS^2$ with the origin at its boundary. 
Then we can choose positive values $a$, $b$, $c$, and $d$ such that the points $a\cdot \tilde p'$, $b\cdot \tilde q'$, $c\cdot \tilde x'$, $d\cdot \tilde y'$ lie in one plane.
Consider the corresponding points $a\cdot p'$, $b\cdot q'$, $c\cdot x'$, $d\cdot y'$ in $\Cone\spc{U}$.
$\CAT(0)$ comparison for these points implies $\CAT(1)$ comparison for the quadruple $ p'$, $q'$, $x'$, $ y'$ in $\spc{U}$.

It remains to consider the case when $\tilde p'$, $\tilde q'$, $\tilde x'$, $\tilde y'$ do not lie in a half-space.
Fix $\tilde z'\in [\tilde x' \tilde y']_{\SSS^2}$.
Observe that 
\begin{align*}\dist{\tilde p'}{\tilde x'}{\SSS^2}+\dist{\tilde q'}{\tilde x'}{\SSS^2}
 &
\le \dist{\tilde p'}{\tilde z'}{\SSS^2}+\dist{\tilde q'}{\tilde z'}{\SSS^2}
\intertext{or} 
\dist{\tilde p'}{\tilde y'}{\SSS^2}+\dist{\tilde q'}{\tilde y'}{\SSS^2}
&\le
\dist{\tilde p'}{\tilde z'}{\SSS^2}+\dist{\tilde q'}{\tilde z'}{\SSS^2}.\end{align*}
That is, in this case, $\CAT(1)$ comparison follows from the triangle inequality.
\qeds

Recall that suspension is a spherical analog of the cone construction; see Section~\ref{sec:constructions}.
The following statement is a direct analog of \ref{ex:cone+susp};
it can be proved along the same lines.

\begin{thm}{Proposition}\label{prop:susp}
Let $\spc{U}$ be a metric space,
and let $\spc{N}$ be a neighborhood of the north pole in $\Susp\spc{U}$ (possibly $\spc{N}=\Susp\spc{U}$).
Then $\spc{N}$ is $\CAT(1)$ if and only if so is $\spc{U}$.
\end{thm}

\section{Polyhedral spaces}

\begin{thm}{Definition}\label{def:poly}
A geodesic space $\spc{P}$ is called a (spherical) \index{polyhedral space}\emph{polyhedral space} 
if it admits a finite triangulation $\tau$ 
such that every simplex in $\tau$ is isometric to a simplex in a Euclidean space (or respectively a unit sphere) of appropriate dimension.

By a \index{triangulation of a polyhedral space}\emph{triangulation} of a polyhedral space, 
we will always understand a triangulation as above. 
\end{thm}

Note that according to the above definition,
all polyhedral spaces are compact.

The \index{dimension of a polyhedral space}\emph{dimension} of a polyhedral space $\spc{P}$
is defined as the maximal dimension of the simplices 
in one (and therefore any) triangulation of~$\spc{P}$.

\parbf{Links.}
Let $\tau$ be a triangulation of a polyhedral space $\spc{P}$,
and let $\sigma$ be a simplex in this triangulation.

The simplices that contain $\sigma$
form an abstract simplicial complex called the \index{link}\emph{link} of $\sigma$, 
denoted by $\Link_\sigma$.
If $m$ is the dimension of $\sigma$,
then the set of vertices of $\Link_\sigma$
is formed by the $(m+1)$-simplices that contain $\sigma$;
the set of its edges is formed by the $(m+2)$-simplices 
that contain $\sigma$; and so on.

The link $\Link_\sigma$
can be identified with the subcomplex of $\tau$ 
formed by all the simplices $\sigma'$ 
such that $\sigma\cap\sigma'=\emptyset$ 
but both $\sigma$ and $\sigma'$ are faces of a simplex of~$\tau$.

The points in $\Link_\sigma$ can be identified with the normal directions to $\sigma$ at a point in its interior.
The angle metric between directions makes $\Link_\sigma$ into a spherical polyhedral space.
We will always consider the link with this metric.

\parbf{Tangent space and space of directions.}
Recall that tangent space and space of directions are defined in Section~\ref{sec:tangent-space+directions}.

Let $\spc{P}$ be a polyhedral space (Euclidean or spherical) and $\tau$ be its triangulation.
If a point $p\in \spc{P}$ 
lies in the interior of a $k$-simplex $\sigma$ of $\tau$ 
then the tangent space $\T_p=\T_p\spc{P}$
is naturally isometric to
\[\EE^k\times(\Cone\Link_\sigma).\]

If $\spc{P}$ is an $m$-dimensional polyhedral space (Euclidean or spherical),
then for any $p\in \spc{P}$
the space of directions $\Sigma_p$ is a spherical polyhedral space
of dimension at most $m-1$. 

In particular, 
for any point $p$ in $\sigma$,
the isometry class of $\Link_\sigma$ together with $k=\dim\sigma$
determines the isometry class of $\Sigma_p$, 
 and the other way around: $\Sigma_p$ and $k$ determine the isometry class of $\Link_\sigma$.

A small neighborhood of $p$ is isometric to a neighborhood of the tip of $\Cone\Sigma_p$ (or $\Susp \Sigma_p$ if $\spc{P}$ is spehrical).
In fact, if this property holds at any point of a compact length space $\spc{P}$,
then $\spc{P}$ is a polyhedral space \cite{lebedeva-petrunin}.

\section{CAT criterion}

The following theorem provides a combinatorial description of polyhedral spaces with curvature bounded above.

\begin{thm}{Theorem}\label{thm:PL-CAT}
Let $\tau$ be a triangulation of a polyhedral space $\spc{P}$. 
The space $\spc{P}$ is locally $\CAT(0)$ if and only if the link of each simplex in $\tau$ has no closed local geodesics shorter than $2\cdot\pi$.

Analogously, let $\spc{P}$ be a spherical polyhedral space and $\tau$ be its triangulation. 
Then $\spc{P}$ is $\CAT(1)$ if and only if neither $\spc{P}$ nor the link of any simplex in $\tau$ contains a closed local geodesic shorter than $2\cdot\pi$.
\end{thm}

\parit{Proof.}
The ``only if'' part follows from \ref{cor:loc-geod-are-min}, \ref{prop:susp}, and \ref{ex:cone+susp}.

To prove the ``if'' part,
we apply induction on $\dim\spc{P}$.
The base case $\dim\spc{P}=0$ is evident.
Let us start with the $\CAT(1)$ case.

\parit{Step.}
Assume that the theorem is proved in the case $\dim\spc{P}<m$.
Suppose $\dim\spc{P}=m$.

Fix a point $p\in\spc{P}$.
A neighborhood of $p$ 
is isometric to a neighborhood of the north pole in the suspension over the space of directions~$\Sigma_p$.

Note that $\Sigma_p$ is a spherical polyhedral space, 
and its links are isometric to links of~$\spc{P}$. 
By the induction hypothesis, $\Sigma_p$ is $\CAT(1)$.
Thus, by \ref{prop:susp}, $\spc{P}$ is locally $\CAT(1)$.
Applying the second part of Corollary~\ref{cor:closed-geod-cat},
we get the statement.

Since the links of $\CAT(0)$ spaces are $\CAT(1)$, the $\CAT(0)$ case is done in exactly the same way except we need to use Proposition~\ref{ex:cone+susp} and the first part of Corollary~\ref{cor:closed-geod-cat} on the very last step when we finally have to talk about $\spc{P}$ (not about its links).
\qeds

\begin{thm}{Exercise}\label{ex:unique-geod=CAT}
Let $\spc{P}$
be a polyhedral space such that any two points can be connected by a unique geodesic.
Show that $\spc{P}$ is $\CAT(0)$.
\end{thm}

\begin{thm}{Advanced exercise}\label{ex:S3}
Construct a Euclidean polyhedral metric on $\mathbb{S}^3$
such that the total angle around each edge in its triangulation is at least $2\cdot \pi$.
\end{thm}

\section{Flag complexes}

\begin{thm}{Definition}\label{def:flag}
A simplicial complex $\mathcal{S}$ 
is called \index{flag complex}\emph{flag} if whenever $\{v_0,\z\ldots,v_k\}$
is a set of distinct vertices of $\mathcal{S}$
that are pairwise joined by edges, then the vertices $v_0,\ldots,v_k$
span a $k$-simplex in~$\mathcal{S}$.

If the above condition is satisfied for $k=2$, 
then we say that $\mathcal{S}$ satisfies 
the \index{no-triangle condition}\emph{no-triangle condition}.
\end{thm}

Note that every flag complex is determined by its one-skeleton.
Moreover, for any graph, its \index{clique}\emph{cliques} (that is, complete subgraphs) define a flag complex.
For that reason, flag complexes are also called \index{clique complex}\emph{clique complexes}.

From the definition of flag complex,
we get the following.

\begin{thm}{Observation}\label{obs:link-of-flag}
Any link of any simplex in a flag complex is flag.
\end{thm}

\begin{thm}[!]{Exercise}\label{ex:baricenric-flag}
Show that the barycentric subdivision of any simplicial complex is a flag complex.

Use the flag condition (see \ref{thm:flag} below)
to conclude that any finite simplicial complex is homeomorphic to a proper length $\CAT(1)$ space.

\end{thm}

\begin{thm}{Proposition}\label{prop:no-trig}
A simplicial complex $\mathcal{S}$ is flag if and only if 
$\mathcal{S}$ as well as the links of all its simplices
satisfy the no-triangle condition.
\end{thm}

\parit{Proof.}
By Observation~\ref{obs:link-of-flag}, the no-triangle condition holds 
for any flag complex and the links of all its simplices.

Now assume that a complex $\spc{S}$ and all its links satisfy 
the no-triangle condition.
It follows that $\spc{S}$ includes a 2-simplex for each triangle.
Applying the same observation for each edge we get that $\spc{S}$ 
includes a 3-simplex for any complete graph with 4 vertices.
Repeating this observation 
for triangles, 
4-simplices,
5-simplices,
and so on, we get that $\spc{S}$ is flag.
\qeds

\parbf{All-right triangulation.}
If each simplex of a triangulation $\tau$ of a spherical polyhedral space $\spc{P}$ is isometric to a spherical simplex with all edges $\tfrac\pi2$,
then $\tau$ is called an \index{all-right triangulation}\emph{all-right triangulation}.
(Note that such a simplex has all right angles.)
Similarly, we say that a simplicial complex 
is equipped with an \index{all-right spherical metric}\emph{all-right spherical metric}
if it is a length metric and each simplex is isometric 
to a spherical simplex with all edges $\tfrac\pi2$.

Spherical polyhedral $\CAT(1)$ spaces glued from right-angled simplices
admit the following characterization 
discovered by Mikhael Gromov \cite[p.~122]{gromov-1987}.

\begin{thm}{Flag condition}\label{thm:flag}
Assume that a spherical polyhedral space $\spc{P}$
admits an all-right triangulation~$\tau$.
Then $\spc{P}$ is $\CAT(1)$
if and only if $\tau$ is flag.
\end{thm}

\parit{Proof; only-if part.} 
Assume there are three vertices $v_1$, $v_2$, and $v_3$ of $\tau$
that are pairwise joined by edges 
but do not span a triangle.
Note that in this case 
\[
\mangle\hinge{v_1}{v_2}{v_3}=
\mangle\hinge{v_2}{v_3}{v_1}=
\mangle\hinge{v_3}{v_1}{v_2}=
\pi.
\]
Equivalently,
\begin{clm}{}\label{clm:3pi/2}
The product
of the geodesics $[v_1v_2]$, $[v_2v_3]$, and $[v_3v_1]$
forms a locally geodesic loop in~$\spc{P}$ of length $\tfrac32\cdot\pi$.
\end{clm}

Now assume that $\spc{P}$ is $\CAT(1)$.
Then by \ref{prop:susp},
$\Link_\sigma\spc{P}$ is $\CAT(1)$ for every simplex $\sigma$ 
in~$\tau$. 

Each of these links is an all-right spherical complex.
By \ref{cor:closed-geod-cat}, none of these links can contain a geodesic circle shorter than $2\cdot\pi$.

Therefore Proposition~\ref{prop:no-trig} and \ref{clm:3pi/2} 
imply the ``only if'' part.

\parit{If part.} 
By \ref{obs:link-of-flag} and \ref{cor:closed-geod-cat},
it is sufficient to show that any closed local geodesic $\gamma$ 
in a flag complex $\spc{S}$ with all-right metric has length at least $2\cdot\pi$.

Recall that the \index{star of vertex}\emph{closed star} of a vertex $v$ (briefly $\overline \Star_v$)
is formed by all the simplices containing~$v$. 
Similarly, $\Star_v$, the open star of $v$, is the union of all simplices containing $v$ with faces opposite $v$ removed.

Choose a vertex $v$ such that $\Star_v$ contains a point $\gamma(t_0)$ of $\gamma$.
Consider the maximal arc $\gamma_v$ of $\gamma$ 
that contains the point $\gamma(t_0)$
and runs in $\Star_v$.
Note that the distance $\dist{v}{\gamma_v(t)}{\spc{P}}$ behaves in exactly the same way 
as the distance from the north pole in~$\mathbb{S}^2$ to a geodesic in the northern hemisphere;
that is, there is a geodesic $\tilde\gamma_v$ in the northern hemisphere of $\mathbb{S}^2$ such that for any $t$ we have
\[\dist{v}{\gamma_v(t)}{\spc{P}}
=
\dist{n}{\tilde\gamma_v(t)}{\mathbb{S}^2},\]
where $n$ denotes the north pole of~$\mathbb{S}^2$.
In particular, 
\[\length\gamma_v=\pi;\]
that is, $\gamma$ spends time $\pi$ on every visit to $\Star_v$.

\begin{wrapfigure}{r}{45mm}
\vskip-0mm
\centering
\includegraphics{mppics/pic-1000}
\end{wrapfigure}

After leaving $\Star_v$,
the local geodesic $\gamma$ has to enter another simplex, 
say~$\sigma'$.
Since $\tau$ is flag, the simplex $\sigma'$
has a vertex $v'$ not joined to $v$ by an edge;
that is, 
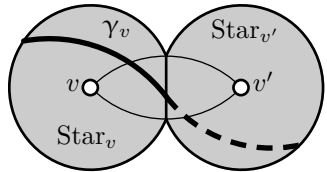
\[\Star_v\cap\Star_{v'}=\emptyset\]

The same argument as above shows that $\gamma$ spends time $\pi$ on every visit to $\Star_{v'}$.
Therefore the total length of $\gamma$ is at least~$2\cdot\pi$.
\qeds

\begin{thm}{Exercise}\label{ex:flag>=pi/2}
Assume that a spherical polyhedral space $\spc{P}$
admits a triangulation $\tau$ such that all edge lengths of all simplices are at least~$\tfrac\pi2$.
Show that $\spc{P}$ is $\CAT(1)$
if $\tau$ is flag.
\end{thm}

\begin{thm}{Exercise}\label{ex:polyhedron-glue}
Let $P$ be a convex polyhedron in $\EE^3$ with $n$ faces $F_1,\ldots, F_n$.
Suppose that each face of $P$ has only obtuse or right angles.
Let us take $2^n$ copies of $P$ indexed by all $n$-bit arrays of 0s and 1s.
Glue two copies of $P$ along $F_i$ if their arrays differ only in the $i$-th bit.
Show that the obtained space is a locally $\CAT(0)$ topological manifold.
\end{thm}

\parbf{The space of trees.}
The following construction is given by
Louis Billera,
Susan Holmes,
and Karen Vogtmann \cite{billera-holmes-vogtmann}.

Let $\spc{T}_n$ be the set of all metric trees with 
$n$ end vertices
labeled by $a_1,\ldots,a_n$.
To describe one tree in $\spc{T}_n$ we may fix a topological tree $t$ with these end vertices,
and all other vertices of degree 3,
and prescribe the lengths of $2\cdot n-3$ edges.
If the length of an edge vanishes, we assume that this edge degenerates;
such a tree can be also described using a different topological tree~$t'$.
The subset of $\spc{T}_n$ corresponding to the given topological tree $t$ can be identified with the octant
\[\set{(x_1,\ldots,x_{2\cdot n-3})\in\mathbb{R}^{2\cdot n-3}}{x_i\ge 0}.\]
Equip each such subset with the metric induced from $\mathbb{R}^{2\cdot n-3}$ and consider the length metric on $\spc{T}_n$ induced by these metrics.

\begin{thm}{Exercise}\label{ex:tree}
Show that $\spc{T}_n$ with the described metric is $\CAT(0)$.
\end{thm}

\section{Cubical complexes}

The definition of a cubical complex
mostly repeats the definition of a simplicial complex, 
with simplices replaced by cubes.

Formally, a \index{cubical complex}\emph{cubical complex} is defined as a subcomplex 
of the unit cube in the Euclidean space $\RR^N$ of large dimension;
that is, a collection of faces of the cube
such that together with each face it contains all its sub-faces.
Each cube face in this collection 
will be called a \index{cube}\emph{cube} of the cubical complex.

Note that according to this definition, 
any cubical complex is finite.

The union of all the cubes in a cubical complex $\spc{Q}$ will be called its \index{underlying space}\emph{underlying space}.
A homeomorphism from the underlying space of $\spc{Q}$ to a topological space $\spc{X}$ is called a \index{cubulation}\emph{cubulation of}~$\spc{X}$.

The underlying space of a cubical complex $\spc{Q}$ will be always considered with the length metric
induced from~$\RR^N$.
In particular, with this metric, 
each cube of $\spc{Q}$ is isometric to the unit cube of the corresponding dimension.

It is straightforward to construct a triangulation 
of the underlying space of $\spc{Q}$ 
such that each simplex is isometric to a Euclidean simplex.
In particular, the underlying space of $\spc{Q}$ is a Euclidean polyhedral space.

The link of a cube in a cubical complex is defined similarly to the link of a simplex in a simplicial complex.
It is a simplicial complex that admits a natural all-right triangulation --- each simplex corresponds to an adjacent cube.

\parbf{Cubical analog of a simplicial complex.}
Let $\spc{S}$ be a finite simplicial complex and $\{v_1,\ldots,v_N\}$ be the set of its vertices.

Consider $\RR^N$ with the standard basis $\{e_1,\ldots,e_N\}$.
Denote by $\square^N$ the standard unit cube in $\RR^N$;
that is, 
\[\square^N=\set{(x_1,\ldots,x_N)\in \RR^N}{0\le x_i\le 1\ \text{for each}\ i}.\]

Given a $k$-dimensional simplex $\<v_{i_0},\ldots,v_{i_k}\>$ in $\spc{S}$, let us mark the $(k\z+1)$-dimensional faces in $\square^N$ (there are $2^{N-k}$ of them) that are parallel to the coordinate $(k+1)$-plane spanned by $e_{i_0},\ldots,e_{i_k}$.

Note that the set of all marked faces of $\square^{N}$
forms a cubical complex;
it will be called 
the \index{cubical analog}\emph{cubical analog} of $\spc{S}$
and will be denoted as $\square_\spc{S}$.

\begin{thm}{Proposition}\label{prob:cubical-analog}
Let $\spc{S}$ be a finite connected simplicial complex
and $\spc{Q}=\square_{\spc{S}}$ be its cubical analog.
Then the underlying space of $\spc{Q}$ is connected 
and the link of any vertex of $\spc{Q}$
is isometric to ${\spc{S}}$
equipped with the all-right spherical metric.

In particular, if $\spc{S}$ is a flag complex,
then $\spc{Q}$ is locally $\CAT(0)$,
and therefore its universal cover $\tilde{\spc{Q}}$ is $\CAT(0)$.
\end{thm}

\parit{Proof.}
By the construction of $\square_{\spc{S}}$, the link of any vertex of $\spc{Q}$ is isometric to ${\spc{S}}$
equipped with the all-right spherical metric.
Further, $\square_{\spc{S}}$ is connected since it must contain the whole 1-dimensional skeleton of $\square^{N}$, which is connected.

If ${\spc{S}}$ is flag, then by the flag condition (\ref{thm:flag}), the link of any cube in $\spc{Q}$ is $\CAT(1)$.
Therefore, by the cone construction (\ref{ex:cone+susp}) $\spc{Q}$ is locally $\CAT(0)$.
It remains to apply the globalization theorem (\ref{thm:hadamard-cartan}).
\qeds

From Proposition \ref{prob:cubical-analog}, 
it follows that the cubical analog
of any flag complex is aspherical.
The following exercise states that the converse also holds; see \cite[5.4]{davis-2001}.

\begin{thm}{Exercise}\label{ex:flag-aspherical}
Show that a finite simplicial complex is flag 
if and only if its cubical analog is aspherical.
\end{thm}

\section{Exotic aspherical manifolds}

By \ref{ex:contractible}, any complete length $\CAT(0)$ space is contractible.
Therefore, by the globalization theorem (\ref{thm:hadamard-cartan}), all proper length, locally $\CAT(0)$ spaces 
are \index{aspherical}\emph{aspherical};
that is, they have contractible universal covers.
This observation will be used to construct examples of aspherical spaces. 

A proper topological space $\spc X$ is called
\index{simply connected space at infinity}\emph{simply connected at infinity} 
if for any compact set $K\subset\spc X$
there is a bigger compact set $K'\supset K$
such that $\spc X\setminus K'$ is path-connected 
and any loop which lies in $\spc X\setminus K'$
is null-homotopic in $\spc X\setminus K$.

Recall that path-connected spaces are not empty by definition.
In particular, compact spaces are not simply connected at infinity.

The following example was constructed by Michael Davis \cite{davis-1983}.

\begin{thm}{Proposition}\label{prop:aspherical}
For any $m\ge 4$, there is a closed aspherical $m$-dimensional manifold whose universal cover is not simply connected at infinity.

In particular, the universal cover of this manifold 
is not homeomorphic to the $m$-dimensional Euclidean space.
\end{thm}

The proof requires the following lemma.

\begin{thm}{Lemma}\label{lem:example-pi_infty}
Let $\spc{S}$ be a finite flag complex,
$\spc{Q}=\square_{\spc{S}}$ be its cubical analog
and $\tilde{\spc{Q}}$ be the universal cover of~$\spc{Q}$.

Assume $\tilde{\spc{Q}}$ is simply connected at infinity.
Then $\spc{S}$ is simply connected.
\end{thm}

\parit{Proof.}
Assume $\spc{S}$ is not simply connected. Equip $\spc{S}$ with an all-right spherical metric.
Choose a shortest noncontractible circle $\gamma\:\mathbb{S}^1\to\spc{S}$ formed by the edges of~$\spc{S}$.

Note that $\gamma$ forms a 1-dimensional subcomplex of $\spc{S}$ which is a closed local geodesic.
Denote by $G$ the subcomplex of $\spc{Q}$ which corresponds to~$\gamma$.

Fix a vertex $v\in G$;
let $G_v$ be the connected component of $v$ in~$G$.
Let $\tilde G$ be a connected component of the inverse image of $G_v$ in $\tilde{\spc{Q}}$
for the universal cover $\tilde{\spc{Q}}\to \spc{Q}$.
Fix a point $\tilde v\in\tilde G$ in the inverse image of~$v$.

\begin{wrapfigure}{r}{25mm}
\vskip-4mm
\centering
\includegraphics{mppics/pic-1100}
\vskip-4mm
\end{wrapfigure}
 
Note that 
\begin{clm}{}\label{tilde-G-convex}
$\tilde G$ is a convex set in~$\tilde{\spc{Q}}$.
\end{clm}

Indeed, according to Proposition \ref{prob:cubical-analog},
$\tilde{\spc{Q}}$ is $\CAT(0)$.
By Exercise \ref{ex:locally-convex},
it is sufficient to show that $\tilde G$ is locally convex in $\tilde{\spc{Q}}$,
or equivalently, $G$ is locally convex in~$\spc{Q}$.

Note that the latter can only fail if $\gamma$ contains two vertices, say $\xi$ 
and 
$\zeta$ in $\spc{S}$,
which are joined by an edge not in $\gamma$; 
denote this edge by~$e$.

Each edge of $\spc{S}$ has length~$\tfrac\pi2$.
Therefore each of the two circles formed by $e$ and an arc of $\gamma$
from $\xi$ to $\zeta$ is shorter than~$\gamma$.
Moreover,
at least one of them is noncontractible 
since $\gamma$ is 
noncontractible.
That is, 
$\gamma$ is not a shortest noncontractible circle --- a contradiction.
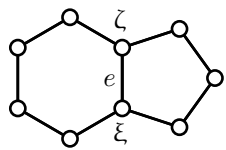
\claimqeds

Observe that $\tilde G$ is homeomorphic to the plane; indeed, $\tilde G$ is
a 2-dimensional manifold without boundary which
by the above is $\CAT(0)$ and hence is contractible.

Denote by $C_R$ the circle of radius $R$ in $\tilde G$ centered at~$\tilde v$.
All $C_R$ are homotopic to each other in $\tilde G\setminus\{\tilde v\}$ and therefore in $\tilde{\spc{Q}}\setminus \{\tilde v\}$.

Note that the map $\tilde{\spc{Q}}\setminus \{\tilde v\}\to \spc{S}$
which returns the direction of $[{\tilde v}{x}]$ for any $x\ne \tilde v$, maps $C_R$ to a circle homotopic to~$\gamma$.
Therefore $C_R$ is not contractible in $\tilde{\spc{Q}}\setminus \{\tilde v\}$.

If $R$ is large, the circle $C_R$ lies outside of any fixed compact set $K'$ in~$\tilde{\spc{Q}}$.
From above $C_R$ is not contractible in $\tilde{\spc{Q}}\setminus K$
if $K\supset \tilde v$.
It follows that $\tilde{\spc{Q}}$ is not simply connected at infinity --- a contradiction.
\qeds

\begin{thm}[!]{Exercise}\label{ex:example-pi_infty-new}
Show that under the assumptions of Lemma~\ref{lem:example-pi_infty},
for any vertex $v$ in $\spc{S}$
the complement $\spc{S}\setminus\{v\}$ is simply connected.
\end{thm}

The proof of \ref{prop:aspherical} will use the following topological fact, which we are not going to prove.
If $m\ge 5$, then it is proved by Michel Kervaire \cite{kervaire}.
For $m=4$, it follows from a construction of Barry Mazur~\cite{mazur}.

\begin{thm}{Fact}\label{thm:fact}
Given an integer $m\ge 4$ there is an $(m-1)$-dimensional smooth homology sphere $\Sigma^{m-1}$ that is not simply connected, and bounds a contractible smooth compact $m$-dimensional manifold~$\spc{W}$.
\end{thm}

\parit{Proof of \ref{prop:aspherical}.}
Choose $\spc{W}$ and $\Sigma^{m-1}$ provided by \ref{thm:fact}.
Pick any triangulation $\tau$ of $\spc{W}$ and let $\spc{S}$ be the resulting subcomplex that triangulates~$\Sigma$.

We can assume that $\spc{S}$ is flag; 
otherwise, pass to the barycentric subdivision 
of $\tau$ and apply Exercise~\ref{ex:baricenric-flag}.

Let $\spc{Q}=\square_{\spc{S}}$ be the cubical analog of~$\spc{S}$.

By Proposition~\ref{prob:cubical-analog},
$\spc{Q}$ is a homology manifold.
It follows that $\spc{Q}$ is a piecewise linear manifold 
with a finite number of singularities at its vertices.

Removing a small contractible neighborhood $V_v$ of each vertex $v$ in $\spc{Q}$,
we can obtain a piecewise linear manifold $\spc{N}$
whose boundary is formed by several copies of~$\Sigma$.

Let us glue a copy of $\spc{W}$ along its boundary to each copy of $\Sigma$ in the boundary of~$\spc{N}$.
This results in a closed manifold 
$\spc{M}$ with polyhedral metric which is homotopically equivalent to~$\spc{Q}$.

Indeed, since both $V_v$ and $\spc{W}$ are contractible, the identity map of their common boundary $\Sigma$ can be extended to a homotopy equivalence $V_v\to\spc{W}$ relative to the boundary.
Therefore the identity map on $\spc{N}$ extends to homotopy equivalences 
$f\: \spc Q\to \spc M$ and $g\:\spc M\to \spc Q$.

Finally, by Lemma~\ref{lem:example-pi_infty},
the universal cover $\tilde{\spc{Q}}$ of $\spc{Q}$
is not simply connected at infinity.

The same holds for 
the universal cover $\tilde{\spc{M}}$ of $\spc{M}$.
The latter follows since the constructed homotopy equivalences 
$f\: \spc Q\to \spc M$ and $g\:\spc M\z\to \spc Q$ 
lift to {}\emph{proper maps} 
$\tilde f \: \tilde{\spc{Q}}\to \tilde{\spc{M}}$
and $\tilde g \: \tilde{\spc{M}}\to \tilde{\spc{Q}}$;
that is, for any compact sets $A\subset \tilde{\spc{Q}}$ and $B\subset\tilde{\spc{M}}$, the inverse images $\tilde g^{-1}(A)$ and $\tilde f^{-1}(B)$ are compact.
\qeds

\section{Notes}

The following proposition was proved by
Fredric Ancel, 
Michael Davis,
and Craig Guilbault \cite{ancel-davis-guilbault};
it could be considered as a more geometric version of Proposition~\ref{prop:aspherical}.

\begin{thm}{Proposition}\label{prop:loc-CAT-mnfld}
Given $m\ge 5$, there is a Euclidean polyhedral space $\spc{P}$ such that:
\begin{subthm}{}
$\spc{P}$ is homeomorphic to a closed $m$-dimensional manifold.
\end{subthm}

\begin{subthm}{}
$\spc{P}$ is locally $\CAT(0)$.
\end{subthm}

\begin{subthm}{}
The universal cover of $\spc{P}$ is not simply connected at infinity.
\end{subthm}
\end{thm}

There are no such examples in dimensions 3 and 4.
In dimension 3 it was shown by
Dale Rolfsen \cite{rolfsen}
and the 4-dimensional case was done by
Alexander Lytchak,
Koichi Nagano,
and Stephan Stadler
\cite{lytchak-nagano-stadler}.

\parit{Proof.}
Apply Exercise~\ref{ex:example-pi_infty-new} to the barycentric subdivision of the simplicial complex $\spc{S}$ provided by Exercise~\ref{ex:funny-S}.
\qeds

\begin{thm}[!]{Exercise}\label{ex:funny-S}
Given an integer $m\ge 5$,
construct a finite $(m-1)$-dimensional simplicial complex $\spc{S}$ such that $\Cone\spc{S}$ is homeomorphic to $\EE^m$
and $\pi_1(\spc{S}\setminus\{v\})\ne0$ for some vertex $v$ in~$\spc{S}$.
\end{thm} 

Theorem \ref{thm:PL-CAT} gives a good-looking description of polyhedral $\CAT(\kappa)$ spaces,
but in fact, it is hard to check even in very simple cases.
For example, the description of those coverings of $\mathbb{S}^3$ branching at three 
great circles which are $\CAT(1)$ requires quite a bit of work \cite{charney-davis-1993} --- try to guess the answer before reading.

Another example is the braid space $\spc{B}_n$ that is the universal metric cover of $\CC^n$ infinitely branching in complex hyperplanes $z_i=z_j$.
So far it is not known if $\spc{B}_n$ is $\CAT(0)$ for any $n\ge 4$ \cite{panov-petrunin-2016}.
Understanding this space could help to study the braid group.
This circle of questions is closely related to the generalization of the flag condition (\ref{thm:flag}) to spherical simplices with few acute dihedral angles.

The construction used in the proof of Proposition~\ref{prop:aspherical} admits a number of modifications, several of which are discussed in a survey by Michael Davis \cite{davis-2001}.

A similar argument was used by Michael Davis, 
Tadeusz Ja\-nu\-szkie\-wicz,
and Jean-Fran\c{c}ois Lafont \cite{davis-januszkiewicz-lafont}.
They constructed a closed smooth 4-dimensional manifold $M$ with universal cover $\tilde M$ diffeomorphic to $\RR^4$, such that $M$ admits a polyhedral metric which is locally $\CAT(0)$, but does not admit a Riemannian metric with non-positive sectional curvature.
Another example of that type was constructed by Stephan Stadler \cite{stadler}.
There are no lower-dimensional examples of this type ---
the 2-dimensional case follows from the classification of surfaces,
and 
the 3-dimensional case follows from the geometrization conjecture.

It is noteworthy that any complete, simply connected Riemannian manifold with non-positive curvature is homeomorphic to the Euclidean space of the same dimension.
In fact, by the globalization theorem
(\ref{thm:hadamard-cartan}), 
the exponential map at a point of such a manifold is a homeomorphism.
In particular, there is no Riemannian analog of Proposition~\ref{prop:loc-CAT-mnfld}.

Recall that a triangulation of an $m$-dimensional manifold defines a piecewise linear structure if the link of every simplex $\Delta$ is homeomorphic to the sphere of dimension $m-1-\dim\Delta$.
According to Stone's theorem \cite{stone, davis-januszkiewicz}, the triangulation of $\spc{P}$ in Proposition~\ref{prop:loc-CAT-mnfld} 
cannot be made piecewise linear --- despite the fact that $\spc{P}$ is a manifold, its triangulation does not induce a piecewise linear structure.

The flag condition also leads to the so-called {}\emph{hyperbolization} procedure, a flexible tool for constructing aspherical spaces;
a survey on the subject is given by Ruth Charney and Michael Davis \cite{charney-davis-1995}.

The $\CAT(0)$ property of a cube complex admits interesting (and useful) geometric descriptions if one exchanges the $\ell^2$-metric for natural $\ell^1$ or $\ell^\infty$ metrics on each cube.

\begin{thm}{Theorem}
The following three conditions are equivalent.

\begin{subthm}{cube-2} A cube complex $Q$ equipped with the $\ell^2$-metric is $\CAT(0)$.
\end{subthm}

\begin{subthm}{cube-infty} A cube complex $Q$ equipped with the $\ell^\infty$-metric is \index{injective space}\emph{injective};
that is, given a metric space $\spc{X}$ with a subset $A$,
any short map $A\to Q$ can be extended to a short map $\spc{X}\to Q$.
\end{subthm}

\begin{subthm}{cube-1} A cube complex $Q$ equipped with the $\ell^1$-metric is {}\emph{median}.
 The latter means that for any three points $x,y,z$ there is a unique point $m$ (it is called the \index{median}\emph{median} of $x$, $y$, and $z$) and a choice of geodesics $[xy], [xz], [yz]$ such that $[xy]\ni m$, $[xz]\ni m$ and $[yz]\ni m$.
\end{subthm}
\end{thm}

A very readable paper on the subject was written by Brian Bowditch \cite{bowditch-2020}.
For more on $\CAT(0)$ cube complexes, see the draft by Montserrat Casals-Ruiz, Mark Hagen, and Ilya Kazachkov \cite{casals-hagen-kazachkov} and the references therein.

\begin{thm}{Exercise}\label{ex:cube-infty=>cube-2}
Prove the implication \ref{SHORT.cube-infty}$\Rightarrow$\ref{SHORT.cube-2} and/or
\ref{SHORT.cube-1}$\Rightarrow$\ref{SHORT.cube-2} in the theorem.
\end{thm}

All the topics discussed in this lecture link Alexandrov geometry with the fundamental group.
The theory of {}\emph{hyperbolic groups}, 
a branch of {}\emph{geometric group theory}, 
introduced by 
Mikhael Gromov \cite{gromov-1987},
could be considered as a further step in this direction.

%%!TEX root = the-subsets.tex
\chapter{Subsets}\label{chapter:shefel}

Here we give a partial answer to the following question:
\textit{Which subsets of Euclidean space, equipped with their induced length-metrics, are  $\CAT(0)$?}

\section{Motivating examples}\label{sec:solid-parabolonds}

Consider three subgraphs of different quadric surfaces:
\begin{align*}
A&=\set{(x,y,z)\in\EE^3}{z\le x^2+y^2},
\\
B&=\set{(x,y,z)\in\EE^3}{z\le -x^2-y^2},
\\
C&=\set{(x,y,z)\in\EE^3}{z\le x^2-y^2}.
\end{align*}

\begin{thm}{Question}\label{CAT(0)?}
Which of the sets $A$, $B$ and $C$, if equipped with the induced length metric, are $\CAT(0)$ and why?
\end{thm}

The answers are given below, but it is instructive to try to answer yourself before reading further.

\parit{$\bm{A}$.}
No, $A$ is not $\CAT(0)$.

Choose a triangle $[pqr]$ in the $(x,y)$-plane that surrounds the origin.
Lift it a bit up, so $[pqr]$ remains in $A$, and observe that this triangle is not thin in the induced length metric of $A$.
\qeds

Let us give another visual argument more in the style of  differential geometry.
 
\parit{Another way.} The boundary $\partial A$ is the paraboloid described by  $z\z=x^2+y^2$;  in particular it bounds an open convex set in $\EE^3$ whose complement is~$A$.
The nearest-point projection of $A\to\partial A$ is short (Exercise~\ref{ex:closest-point}).
It follows that $\partial A$ is a convex set in $A$ equipped with its induced length metric.

Therefore if $A$ is $\CAT(0)$, then so is~$\partial A$.
The latter is not true: $\partial A$ is a smooth convex surface, and has strictly positive curvature by the Gauss formula.
\qeds

\parit{$\bm{B}$.} Yes, $B$ is $\CAT(0)$.

Evidently $B$ is a closed convex set in~$\EE^3$.
Therefore the length metric on $B$ coincides with the Euclidean metric
and $\CAT(0)$ comparison holds.\qeds

\parit{$\bm{C}$.} Yes, $C$ is $\CAT(0)$,
but the proof is not as easy as before.

Set $f_t(x,y)=x^2-y^2 -2\cdot (x-t)^2$.
Consider the one-parameter family of sets 
\[V_t=\set{(x,y,z)\in\EE^3}{z\le f_t(x,y)}.\]

\begin{wrapfigure}{r}{25mm}
\vskip-10mm
\centering
\includegraphics{mppics/pic-1201}
\end{wrapfigure}

Each set $V_t$ is a solid paraboloid tangent to $\partial C$ along the parabola $y\mapsto(t,y,t^2-y^2)$.
The set $V_t$ is closed and convex for any $t$, and
\[C=\bigcup_t V_t.\]
Further note that the function $t\mapsto f_t(x,y)$ is concave for any fixed $x,y$.
Therefore
\begin{clm}{}\label{eq:VcVnV-0}
if $a<b<c$, then $V_b\supset V_a\cap V_c$.
\end{clm}

Given $t_1,\ldots, t_n\in \RR$ consider the  union
\[C'=V_{t_1}\cup\ldots\cup V_{t_n}.\]
The inclusion \ref{eq:VcVnV-0} makes it possible to apply Reshetnyak gluing theorem \ref{thm:gluing} recursively and show that $C'$ is $\CAT(0)$.

By approximation, the $\CAT(0)$ comparison holds for any 4 points in $C$;
hence $C$ is $\CAT(0)$.
More precisely, choose $x_1,x_2,x_3,x_4\in C$ and 6 geodesics $[x_ix_j]_C$ between them.
Choose $\eps>0$, shift each $[x_ix_j]_C$ down by $\eps$, and reconnect it to $x_i$ and $x_j$ by vertical $\eps$-segments.
Denote the obtained curve by $\gamma_{i,j}$;
note that 
\[\length\gamma_{i,j}=\dist{x_i}{x_j}{C}+2\cdot\eps.\]
We may assume that $C'$ contains all $\gamma_{i,j}$.
It follows that 
\[\dist{x_i}{x_j}{C}\le \dist{x_i}{x_j}{C'}\le \dist{x_i}{x_j}{C}+2\cdot\eps\]
Since $\eps>0$ is arbitrary and $C'$ is $\CAT(0)$, so is $C$.\qeds

\parbf{Remarks.}
The last argument will be reused in the proof of \ref{thm:set-with-smooth-bry:CBA}.
While the set $C$ is not convex, we shall see that it is \textit{two-convex} as defined in the next section.
In fact, two-convexity is closely related to the inheritance of an upper curvature bound by a subset.

\section{Two-convexity}

\begin{thm}{Definition}\label{def:two-convex}
We say that a subset $K\subset \EE^m$ is \index{two-convex set}\emph{two-convex}
if the following condition holds for any plane $W\subset\EE^m$:
If $\gamma$ is a simple closed curve in $W\cap K$ 
that is null-homotopic in $K$,  
then it is null-homotopic in $W\cap K$, and in particular the disc in $W$ bounded by $\gamma$ lies in~$K$.
\end{thm}

Note that two-convex sets do not have to be connected or simply connected. 
The following two propositions follow immediately from the definition.

\begin{thm}{Proposition}
Any subset in $\EE^2$ is two-convex.
\end{thm}

\begin{thm}{Proposition}\label{prop:two-hull}
The intersection of an arbitrary collection of two-convex sets in $\EE^m$ is two-convex.
\end{thm}

\begin{thm}{Proposition}\label{prop:two-hull-open}
The interior of any two-convex set in $\EE^m$ is a two-convex set.
\end{thm}

\parit{Proof.}
Fix a  two-convex set $K\subset \EE^m$ and a 2-plane $W$;
denote by \index{3@$\Int$ (interior)}$\Int K$ the interior of $K$.
Let $\gamma$ be a closed simple curve in $W\cap \Int K$ 
that  is contractible in the interior of~$K$.

Since $K$ is two-convex,
the plane disc $D$ bounded by $\gamma$ lies in~$K$.
The same holds for the translations of $D$ by small vectors.
Therefore $D$ lies in $\Int K$;
that is, $\Int K$ is two-convex.
\qeds

\begin{thm}{Definition}
Given a subset $K\subset \EE^m$, define its two-convex hull (briefly, $\Conv_2K$) as the intersection of all two-convex subsets containing~$K$.
\end{thm}

Note that by Proposition~\ref{prop:two-hull},
the two-convex hull of any set is two-convex.
Further, 
by \ref{prop:two-hull-open}, the
two-convex hull of an open set is open.

The next proposition describes  closed two-convex sets with smooth boundary.

\begin{thm}{Proposition}\label{prop:two-cove+smooth}
Let $K\subset \EE^m$ be a closed subset.

Assume that the boundary of $K$ is a smooth hypersurface~$S$.
Consider the unit normal vector field $\nu$ on $S$ that  points outside of~$K$.
Denote by $k_1\le \ldots\le k_{m-1}$ the principal curvature functions of $S$ with respect to $\nu$ (note that if $K$ is convex, then $k_1\ge 0$, and for connected $K$ the converse is also true).

Then $K$ is two-convex if and only if $k_2(p)\ge 0$ for any point $p\in S$.
Moreover, if $k_2(p)<0$ at some point $p$, then Definition~\ref{def:two-convex} fails for some curve $\gamma$ forming a triangle in an arbitrarily small neighborhood of~$p$.
\end{thm}

The proof uses a straightforward modification of the Morse theory for manifolds with boundary; the paper of Sergei Vakhrameev \cite{vakhrameev} contains all the necessary lemmas.

\begin{wrapfigure}{r}{25mm}
\vskip-4mm
\centering
\includegraphics{mppics/pic-1202}
\end{wrapfigure}

\parit{Proof; only-if part.}
Suppose $k_2(p)<0$ at some point $p\in S$,
consider the plane $W$ containing $p$ and spanned by the first two principal directions at~$p$.
Choose a small triangle in $W$ that surrounds $p$ and move it slightly in the direction of $\nu(p)$.
We get a triangle $\trig xyz$ that is null-homotopic in $K$,
but the solid triangle $\Delta=\Conv\{x,y,z\}$ bounded by $\trig xyz$ does not lie in $K$ completely.
Therefore $K$ is not two-convex.

\parit{If part.}
Recall that a smooth function $f\:\EE^m\to\RR$ is called \index{strongly convex function}\emph{strongly convex} if its Hessian is positive definite at each point.

Suppose $f\:\EE^m\to\RR$ is a smooth
 strongly
 convex function such that the restriction $f|_S$ is a Morse function.
Note that a generic smooth 
strongly 
convex function $f\:\EE^m\to\RR$ has this property.

For a critical point $p$ of $f|_S$, the outer normal vector $\nu(p)$ is parallel to the gradient $\nabla_pf$;
we say that $p$ is a 
\index{positive critical point}\emph{positive critical point}
if $\nu(p)$ and $\nabla_p f$ point in the same direction, 
and 
\index{negative critical point}\emph{negative} otherwise.
If $f$ is generic, then we can assume that the sign is defined for all critical points;
that is, $\nabla_pf\ne0$ for any critical point $p$ of~$f|_S$.

Since $k_2\ge 0$ and the function $f$ is  strongly
 convex, 
the negative critical points of $f|_S$
have index at most~1.

Given a real value $s$, set 
\[K_s=\set{x\in K}{f(x)<s}.\]
Assume  $\phi_0\:\DD\to K$ is a continuous map of the disc $\DD$
such that $\phi_0(\partial \DD)\subset K_s$.

Note that by the Morse lemma, 
there is a homotopy $\phi_t\:\DD\to K$ rel $\partial \DD$ such that 
$\phi_1(\DD)\subset K_s$.

Indeed, we can construct a homotopy $\phi_t\:\DD\to K$ that decreases the maximum of $f\circ\phi$ on $\DD$ until the maximum occurs at a critical point $p$ of~$f|_S$.
This point cannot be negative; otherwise, its index would be at least 2.
If this critical point is positive, then it is easy to decrease the maximum a little by pushing the disc from $S$ into $K$ in the direction of $-\nabla f_p$.

Consider a closed curve $\gamma\:\mathbb{S}^1\to K$ that is null-homotopic in~$K$.
Note that the distance function \[f_0(x)\z=\dist{\Conv\gamma}{x}{\EE^m}\] 
is convex on $\EE^m$ by Exercise \ref{ex:displacement}.
Therefore $f_0$ can be approximated by smooth strongly convex functions $f$ in general position.
From above, there is a disc in $K$ with boundary $\gamma$
that lies arbitrarily close to $\Conv\gamma$.
Since $K$ is closed, the statement follows.
\qeds

Note that the ``if'' part proves a somewhat stronger statement.
Namely,  any plane curve $\gamma$ (not necessary simple) which is  contractible in $K$ is also contractible in the intersection of $K$ with the plane of~$\gamma$.
The latter condition does not hold for the complement 
of two planes in $\EE^4$, which is two-convex by Proposition~\ref{prop:two-hull};
see also Exercise~\ref{ex:two-planes} below.
By the following proposition, there are no such examples in~$\EE^3$.

\begin{thm}{Proposition}\label{prop:3d-strong-2-convexity}
Let $\Omega\subset\EE^3$ be an open two-convex subset.
Then for any plane $W\subset\EE^3$, 
any closed curve in $W\cap \Omega$ 
that is null-homotopic in $\Omega$ is also null-homotopic in $W\cap \Omega$.
\end{thm}

This statement is intuitively obvious, but the proof is not trivial;
it uses the following classical result.

\begin{thm}{Loop theorem}
Let $M$ be a 3-dimensional manifold with nonempty boundary $\partial M$.
Assume 
$f\: (\DD,\partial \DD)\to (M,\partial M)$
is a continuous map of pairs
such that the boundary curve $f|_{\partial \DD}$ is not null-homotopic in $\partial M$; here $\DD$ denotes the disc.
Then there is an embedding of pairs $h\: (\DD,\partial \DD)\to (M,\partial M)$ with the same property.
\end{thm}

The theorem is due to Christos Papakyriakopoulos \cite[a proof can be found in][]{hatcher}. 

\parit{Proof of \ref{prop:3d-strong-2-convexity}.}
Fix a closed plane curve $\gamma$ in $W\cap \Omega$ that  is null-homotopic in $\Omega$. Suppose $\gamma$ is not contractible in  $W\cap \Omega$.

Let $\phi\: \DD\to\Omega$ be a map of the disc with the boundary curve~$\gamma$.

Since $\Omega$ is open we can first change  $\phi$  slightly so that $\phi(x)\notin W$ for $1-\eps<|x|<1$ for some small $\eps>0$.
By further changing $\phi$ slightly we can assume that it is transversal to $W$  on  $\Int \DD$ and agrees with the previous map near~$\partial \DD$.

This means that $\phi^{-1} (W)\cap \Int \DD$ consists of finitely many simple closed curves which cut $\DD$ into several components. 
Consider one of the ``innermost'' components $c'$;
that is, $c'$ is a boundary curve of a disc $\DD'\subset \DD$,
$\phi(c')$ is a closed curve in $W$ and $\phi(\DD')$  completely lies in one of the two half-spaces  with boundary~$W$. 
Denote this half-space by~$H$.

If $\phi(c')$ is not contractible in $W\cap \Omega$, then applying the loop theorem to $M^3=H\cap \Omega$ we conclude that  there exists a {}\emph{simple} closed curve $\gamma'\subset \Omega\cap W$ which is not contractible in $\Omega\cap W$ but is contractible in $\Omega\cap H$. 
This contradicts two-convexity of~$\Omega$. 

Hence $\phi(c')$ is contractible in $W\cap \Omega$. Therefore $\phi$ can be changed in a small neighborhood $U$ of $\DD'$ so that the new map $\hat\phi$ maps $U$ to one side of~$W$. 
In particular, the set $\hat\phi^{-1}(W)$ consists of the same curves as $\phi^{-1} (W)$ with the exception of~$c'$.

Repeating this process several times we reduce the problem to the case where $\phi^{-1} (W)\cap \Int \DD=\emptyset$.
This means that $\phi(\DD)$ lies entirely in one of the half-spaces bounded by $W$.

Again applying the loop theorem, we obtain a simple closed curve in $W\cap \Omega$ which is not contractible in $W\cap \Omega$ but is contractible in~$\Omega$. 
This again contradicts two-convexity of~$\Omega$. 
Hence $\gamma$ is contractible in  $W\cap \Omega$ as claimed.
\qeds

\section{Sets with smooth boundary}\label{sec:smooth-bry}

In this section, we characterize the subsets with smooth boundary in $\EE^m$ that form $\CAT(0)$ spaces. 

\begin{thm}{Smooth two-convexity theorem}\label{thm:set-with-smooth-bry:CBA}
Let $K$ be a closed, simply connected subset in $\EE^m$ equipped with the induced length metric.
Assume $K$ is bounded by a smooth hypersurface.
Then 
$K$ is $\CAT(0)$ if and only if $K$ is two-convex.
\end{thm}

The proof below is based on the argument in \ref{CAT(0)?}$C$.

\parit{Proof.}
Denote by $S$ and by $\Omega$ the boundary and the interior of $K$ respectively. 
Since $K$ is connected and $S$ is smooth, $\Omega$ is also connected.

Denote by $k_1(p)\le\ldots\le k_{m-1}(p)$ the principal curvatures of $S$ at $p\in S$ with respect to the normal vector $\nu(p)$ pointing out of~$K$.
By Proposition~\ref{prop:two-cove+smooth}, $K$ is two-convex if and only if $k_2(p)\ge 0$ for any $p\in S$.

\parit{Only-if part.}
Assume $K$ is not two-convex.
Then by Proposition~\ref{prop:two-cove+smooth}, $k_2(p)<0$ at some point $p\in S$. Hence there is a triangle $\trig xyz$ in $K$ that is null-homotopic in $K$, but the solid triangle $\Delta=\Conv\{x,y,z\}$ does not lie in $K$ completely.
Evidently the triangle $\trig xyz$ is not thin in~$K$. 
Hence $K$ is not $\CAT(0)$.

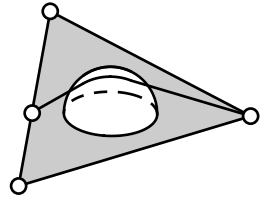
\parit{If part.}
Since $K$ is simply connected,
by the globalization theorem (\ref{thm:hadamard-cartan})
it suffices to show that any point $p\in K$ admits a $\CAT(0)$ neighborhood.

Evidently, such a neighborhood exists if $p\in\Int K$

Assume $p\in S$ and $k_2(p)>0$.
Fix a sufficiently small $\eps>0$ and set $K'=K\cap \cBall[p,\eps]$.
Let us show that 
\begin{clm}{}\label{K'-is-CAT}
$K'$ is $\CAT(0)$.
\end{clm}

Consider the coordinate system with the origin at $p$
and the principal directions and $\nu(p)$ as the coordinate directions.
For small $\eps>0$, the set $K'$ 
can be described as a subgraph
\[K'
=
\set
{(x_1,\ldots,x_m)\in \cBall[p,\eps]}
{x_m\le f(x_1,\ldots,x_{m-1})}.\]

Fix $s\in[-\eps,\eps]$.
Since $\eps$ is small and $k_2(p)>0$, the restriction 
$f|_{x_1=s}$
is concave in the $(m-2)$-dimensional cube defined by the inequalities $|x_i|<2\cdot\eps$ for $2\le i\le m-1$.

Fix a negative real value $\lambda<k_1(p)$.
Given $s\in (-\eps,\eps)$,
consider the set 
\[V_s
=
\set
{(x_1,\ldots,x_m)\in K'}
{x_m\le f(x_1,\ldots,x_{m-1})+\lambda\cdot (x_1-s)^2}.\]
Note that the function 
\[(x_1,\ldots, x_{m-1})\mapsto f(x_1,\ldots,x_{m-1})+\lambda\cdot (x_1-s)^2\]
is concave near the origin.
Since $\eps$ is small, we can assume that the $V_s$ are convex subsets of~$\EE^m$.

Further note that 
\[K'=\bigcup_{s\in[-\eps,\eps]}V_s.\]

Also, the same argument as in \ref{CAT(0)?}$C$ shows that
\begin{clm}{}\label{eq:VcVnV}
If $a<b<c$, then $V_b\supset V_a\cap V_c$.

\end{clm}

Given an array of values $s_1<\ldots<s_k$ in $[-\eps,\eps]$,
set $V_i=V_{s_i}$ and
consider the unions 
\[W_i=V_1\cup\ldots\cup V_i\]
equipped with the induced length metric.

Note that the array $(s_1,\ldots,s_k)$ can be chosen in such a way that
$W_k$ is arbitrarily close to~$K'$ in the sense of Hausdorff.

Arguing as in \ref{CAT(0)?}$C$, we get that the following claim implies \ref{K'-is-CAT}.
\begin{clm}{}
All $W_i$ are $\CAT(0)$.
\end{clm}

This claim is proved by induction.
Base: $W_1=V_1$ is $\CAT(0)$ as a convex subset in~$\EE^m$.

\parit{Induction step:} Assume that $W_i$ is $\CAT(0)$.
According to \ref{eq:VcVnV}, 
\[V_{i+1}\cap W_i=V_{i+1}\cap V_i.\]
Moreover, this is a convex set in $\EE^m$ 
and therefore it is a convex set in $W_i$ and in $V_{i+1}$.
By the Reshetnyak gluing theorem (\ref{thm:gluing}), $W_{i+1}$ is $\CAT(0)$.
Hence the claim follows.
\claimqeds

Summarizing, we have proved the following claim:

\begin{clm}{}\label{clm-strong2convex}
$K'$ is $\CAT(0)$ if $K$ is 
\index{strongly two-convex set}\emph{strongly two-convex},
that is, $k_2(p)>0$ at any point $p\in S$.
\end{clm}

It remains to show that $p$ admits a $\CAT(0)$ neighborhood in the case $k_2(p)=0$.

Choose a coordinate system $(x_1,\ldots,x_m)$ as above,
so that the $(x_1,\ldots,x_{m-1})$-coordinate hyperplane is the tangent subspace to $S$ at~$p$.

Fix $\eps>0$ so that a neighborhood of $p$ in $S$ 
is the graph
\[x_m= f(x_1,\ldots,x_{m-1})\]
of a function $f$ defined on the open ball $B$ of radius $\eps$ centered at the origin in the $(x_1,\ldots,x_{m-1})$-hyperplane.
Fix a smooth positive 
strongly 
convex function $\phi\:B\to \RR_+$
such that $\phi(x)\to\infty$ as $x$ approaches the boundary of $B$.
Note that for $\delta>0$, the subgraph $K_\delta$ defined by the inequality
\[x_m\le f(x_1,\ldots,x_{m-1})-\delta\cdot\phi(x_1,\ldots,x_{m-1})\]
is strongly two-convex.
By \ref{clm-strong2convex}, $K_\delta$ is $\CAT(0)$.

Finally as $\delta\to0$, the closed $\eps$-neighborhoods of $p$ in $K_\delta$ 
converge to the closed $\eps$-neighborhood of $p$ in $K$.
By \ref{prop:cat-limit}, the $\eps$-neighborhood of $p$ is $\CAT(0)$.
\qeds

Note that two-convexity survives under affine transformations of a Euclidean space.
Therefore, as a consequence of the smooth two-convexity theorem (\ref{thm:set-with-smooth-bry:CBA}), we get the following.

\begin{thm}{Corollary}\label{cor:affine}
Let $K\subset \EE^m$ be a closed connected set that admits a smooth parametrization by a closed $m$-ball.
If $K$ equipped with the induced length metric is $\CAT(0)$, then the same holds for any affine transformation of $K$.
\end{thm}

Corollary~\ref{cor:shefel} will provide a more general statement in the 3-dimensional Euclidean space.

\section{Open plane sets}

In this section, we consider inheritance of upper curvature bounds by subsets of the Euclidean plane.

\begin{thm}{Theorem}\label{thm:bishop-plane}
Let $\Omega$ be an open simply connected subset of~$\EE^2$.
Equip $\Omega$ with the induced length metric and denote its completion 
by~$K$.
Then $K$ is $\CAT(0)$.
\end{thm}

If $\Omega$ is not open, then the induced intrinsic metric might be infinite for some pair of points.
In this case the proof implies that each metric component of $\Omega$ with the induced length metric is $\CAT(0)$.

\parit{Sketch of proof.}
It is sufficient to show that any triangle in $K$ is thin,
as defined in \ref{def:k-thin}.

Note that $K$ admits a length-preserving map to $\EE^2$ that extends the embedding $\Omega\hookrightarrow\EE^2$.
Therefore each triangle $\trig xyz$ in $K$ can be mapped to the plane in a length-preserving way.
Since $\Omega$ is simply connected, any open region, say $\Delta$, that is surrounded by the image of $\trig xyz$ lies completely in~$\Omega$.

\begin{wrapfigure}{r}{35mm}
\vskip-4mm
\centering
\includegraphics{mppics/pic-1203}
\end{wrapfigure}

Note that in each triangle $\trig xyz$ in $K$, the sides $[xy]$, $[yz]$ and $[zx]$ intersect each other along a geodesic starting at a common vertex, possibly a one-point geodesic.
In other words, every triangle in $K$ looks like the one in the diagram. 

Indeed, assuming the contrary, there will be a {}\emph{lune} in $K$ bounded by two minimizing geodesics with common ends but no other common points.
The image of this lune in the plane must have concave sides, since otherwise one could shorten the sides by pushing them into the interior.
Evidently, there is no plane lune with concave sides ---  a contradiction.

Note that it is sufficient to consider only simple triangles $\trig xyz$, 
that is, triangles whose sides $[xy]$, $[yz]$ and $[zx]$ intersect each other only at the common vertices.
If this is not the case, chopping the overlapping part of sides reduces to the case when the triangle is simple (this is formally stated in Exercise~\ref{ex:chopping-triangle}).

Again, the open region, say $\Delta$, bounded by the image of $\trig xyz$ has concave sides in the plane, since otherwise one could shorten the sides by pushing them into~$\Omega$.
It remains to solve Exercise~\ref{ex:concave-triangle}.
\qeds

\begin{thm}[!]{Exercise}\label{ex:chopping-triangle}
Assume that $[pq]$ is a common part of the two sides $[px]$ and $[py]$ of the triangle $[pxy]$.
Consider the triangle $[qxy]$ whose sides are formed by arcs of the sides of~$[pxy]$.
Show that if $[qxy]$ is thin, then so is~$[pxy]$.
\end{thm}

\begin{thm}[!]{Exercise}\label{ex:concave-triangle}
Assume $S$ is a closed plane 
region whose boundary is a plane triangle $T$ with concave sides.
Equip $S$ with the induced length metric.
Show that the triangle $T$ is thin in~$S$.
\end{thm}

\begin{thm}{Proposition}\label{prop:bishop-sphere}
Let $\Theta$ be an open connected subset of the unit sphere $\mathbb{S}^2$ that does not contain a closed hemisphere.
Equip $\Theta$ with the induced length metric.
Let $\tilde \Theta$ be a metric cover of $\Theta$ 
such that any closed curve in $\tilde \Theta$ shorter than $2\cdot\pi$ is contractible.

Show that the completion of $\tilde \Theta$ is $\CAT(1)$.
\end{thm}

This proposition will be used in the next section.
It is a spherical analog of Theorem \ref{thm:bishop-plane}
and can be proved along the same lines with an additional idea from the following exercise.

\begin{thm}{Exercise}\label{ex:bishop-sphere}
Let $K$ be a closed subset of the unit sphere $\mathbb{S}^2$.
Suppose $K$ is simply connected, bounded by a simple closed Lipschitz curve, and its interior does not contain a closed hemisphere.
Show that $K$ with the induced length metric is $\CAT(1)$.
\end{thm}

\section{Two-convexity theorem}

\begin{thm}{Theorem}\label{thm:shefel}
Let $\Omega$ be a connected open set in~$\EE^3$.
Equip $\Omega$ with the induced length metric
and denote by $\tilde K$ the completion of the universal metric cover of~$\Omega$.
Then $\tilde K$ is $\CAT(0)$ 
if and only if $\Omega$ is two-convex.
\end{thm}

The proof of this statement will be given in the following three sections.
First we prove its polyhedral analog, then we prove some properties of two-convex hulls in the 3-dimensional Euclidean space and only then do we prove the general statement.

The following exercise shows that the analogous statement does not hold in higher dimensions.

\begin{thm}{Exercise}\label{ex:two-planes}
Let $\Pi_1,\Pi_2$ be two planes in $\EE^4$ intersecting at a single point.
Let $\tilde K$ be the completion of the universal metric cover of $\EE^4\setminus(\Pi_1\cup\Pi_2)$.

Show that 
$\tilde K$ is $\CAT(0)$ if and only if $\Pi_1\perp\Pi_2$.

\end{thm}

Before coming to the proof of the two-convexity theorem, 
let us formulate a few corollaries.
The following corollary is a generalization of the smooth two-convexity theorem (\ref{thm:set-with-smooth-bry:CBA}) for 3-dimensional Euclidean space.

\begin{thm}{Corollary}\label{cor:shefel}
Let $K$ be a closed subset in $\EE^3$.
Assume there is a Lipschitz homeomorphism from a closed ball in $\EE^3$ to $K$.
Then $K$ with the induced length metric is $\CAT(0)$ if and only if the interior of $K$ is two-convex.
\end{thm}

\parit{Proof.}
Set $\Omega=\Int K$.
By the domain invariance theorem, $\Omega$ is homeomorphic to an open ball; in particular it is simply connected.

Apply the two-convexity theorem to~$\Omega$.
Note that the completion of $\Omega$ equipped with the induced length metric 
is isometric to $K$ with the induced length metric.
Hence the result.
\qeds

Note that the Lipschitz condition is used just once to show that the completion of $\Omega$ is isometric to $K$ with the induced length metric.
This property holds for a wider class of hypersurfaces;
for instance, some examples of Alexander horned balls have $\CAT(0)$ induced length metric.

Let $U$ be an open set in~$\RR^2$.
A continuous function $f\:U\to\RR$ is called 
\index{saddle function}\emph{saddle} 
if for any linear function $\ell\:\RR^2\to\RR$, the difference 
$f-\ell$
does not have local maxima or local minima in~$U$.
Equivalently, the open subgraph and epigraph of $f$
\begin{align*}
&\set{(x,y,z)\in\EE^3}{z<f(x,y),\ (x,y)\in U},
\\
&\set{(x,y,z)\in\EE^3}{z>f(x,y),\ (x,y)\in U}
\end{align*}
are two-convex. 

\begin{thm}{Theorem}\label{thm:shefel-graph}
Let $f\:\DD\to \RR$ be a Lipschitz function which is saddle in the interior of the closed unit disc~$\DD$. 
Then the graph
\[\Gamma=\set{(x,y,z)\in \EE^3}{z=f(x,y)},\] 
equipped with the induced length metric is $\CAT(0)$.
\end{thm}

\parit{Proof.}
Since the function $f$ is Lipschitz,
its graph $\Gamma$ with the induced length metric is bi-Lipschitz equivalent to $\DD$ with the Euclidean metric.

\begin{wrapfigure}{r}{21mm}
\vskip-4mm
\centering
\includegraphics{mppics/pic-1206}
\end{wrapfigure}

Consider the sequence of sets 
\[K_n
=
\set{(x,y,z)\in \EE^3}{z\lege f(x,y)\pm\tfrac1n,\ (x,y)\in \DD}.\]
Note that each $K_n$ meets the conditions in Corollary~\ref{cor:shefel}.
Therefore, $K_n$ equipped with the induced length metric is $\CAT(0)$.
It remains to note that $K_n\to \Gamma$ in the sense of Gromov--Hausdorff, and apply \ref{prop:cat-limit}.
\qeds

\section{Polyhedral case}

Now we are back to the proof of the two-convexity theorem (\ref{thm:shefel}).

Recall that a subset $P$ of $\EE^m$ is called \index{polyhedral set}\emph{polyhedral}
if it can be presented as a union of a finite number of simplices.
Similarly,
a \index{spherical polyhedral set}\emph{spherical polyhedral set}
is a union of a finite number of simplices in~$\mathbb{S}^m$.

Note that any polyhedral set admits a finite triangulation.
Therefore any polyhedral set equipped with the induced length metric
forms a Euclidean polyhedral space as defined in~\ref{def:poly}.
 
\begin{thm}{Lemma}\label{lem:poly-shefel}
The two-convexity theorem (\ref{thm:shefel}) holds if the set $\Omega$ is the interior of a polyhedral set.
\end{thm}

The statement might look obvious, but we will face a technical difficulty in the proof that is related to the following.
Let $P$ be a polyhedral set and $\Omega$ its interior,
both considered with the induced length metrics.
If the completion $K$ of $\Omega$
is isometric to~$P$, then the lemma follows easily from
\ref{thm:PL-CAT}.

\begin{wrapfigure}{r}{20mm}
\vskip-0mm
\centering
\includegraphics{mppics/pic-1207}
\end{wrapfigure}

However in general
we only have a locally distance-preserving map $K\to P$;
it does not have to be onto and it may not be injective. 
An example can be guessed from the picture.
Nevertheless, it is easy to see that $K$ is always a polyhedral space.

The proof uses the following two exercises.

\begin{thm}[!]{Exercise}\label{ex:hemisphere}
Show that any closed path of length $<2\cdot \pi$ in the unit sphere $\mathbb{S}^2$ lies in an open hemisphere.
\end{thm}

\begin{thm}[!]{Exercise}\label{ex:inner-support}
Assume $\Omega$ is an open subset in $\EE^3$
that is not two-convex.
Show that there is a plane $W$ such that the complement 
$W\setminus \Omega$ contains an isolated point and a small circle around this point in $W$ is contractible in $\Omega$.
\end{thm}

\parit{Proof of \ref{lem:poly-shefel}.}
The ``only if'' part can be proved in the same way as in the smooth two-convexity theorem (\ref{thm:set-with-smooth-bry:CBA}) with additional use of Exercise~\ref{ex:inner-support}.

\parit{If part.}
Assume that $\Omega$ is two-convex.
Denote by $\tilde\Omega$ the universal metric cover of~$\Omega$.
Let $\tilde K$ and $K$ be the completions of $\tilde\Omega$ and $\Omega$, respectively.

The main step is to show that $\tilde K$ is $\CAT(0)$. 

Note that $K$ is a polyhedral space and the covering $\tilde\Omega\to\Omega$ extends to a covering map $\tilde K\to K$ which might be branching at some vertices.%
\footnote{For example, if $K=\set{(x,y,z)\in\EE^3}{|z|\le |x|+|y|\le 1}$ and $p$ is the origin, then $\Sigma_p$,
the space of directions at $p$,
is not simply connected and $\tilde K\to K$ branches at~$p$.}

Fix a point $\tilde p\in \tilde K\setminus\tilde\Omega$; 
denote by $p$ the image of $\tilde p$ in~$K$.
Note that $\tilde K$ is a ramified cover of $K$ and hence is locally contractible.
Thus, any loop in $\tilde K$ is homotopic to a loop in $\tilde\Omega$.
Since $\tilde\Omega$ is simply connected, so is $\tilde K$.

Thus, by the globalization theorem (\ref{thm:hadamard-cartan}), it is sufficient to show that

\begin{clm}{}\label{eq:curv=<0}
a small neighborhood of $\tilde p$ in $\tilde K$ is $\CAT(0)$.
\end{clm}

Recall that $\Sigma_{\tilde p}=\Sigma_{\tilde p}\tilde K$ denotes the space of directions at~$\tilde p$.
Note that a small neighborhood of $\tilde p$ in $\tilde K$
is isometric to an open set in the cone over $\Sigma_{\tilde p}\tilde K$.
By \ref{ex:cone+susp}, \ref{eq:curv=<0} follows once we can show that

\begin{clm}{}\label{eq:curv=<1}
$\Sigma_{\tilde p}$ is $\CAT(1)$.
\end{clm}

By rescaling, we can assume that all faces of $K$ that do not contain $p$ lie at distance at least 2 from~$p$.
Let $\mathbb{S}^2$ be the unit sphere centered at $p$,
and let $\Theta=\mathbb{S}^2\cap\Omega$.
Note that $\Sigma_pK$ is isometric to the completion of $\Theta$
and $\Sigma_{\tilde p}\tilde K$ is the completion of the regular metric covering $\tilde\Theta$ of $\Theta$ induced by the universal metric cover $\tilde \Omega\to \Omega$.

By \ref{prop:bishop-sphere}, it remains to show the following:
\begin{clm}{}
Any closed curve in $\tilde\Theta$ shorter than $2\cdot\pi$ is contractible.
\end{clm}

Fix a closed curve $\tilde \gamma$ of length $<2\cdot\pi$ in~$\tilde\Theta$.
Its projection $\gamma$ in $\Theta\subset\mathbb{S}^2$ has the same length.
Therefore, by Exercise~\ref{ex:hemisphere}, $\gamma$ lies in an open hemisphere.
Then for a plane $\Pi$ passing close to $p$,
the central projection $\gamma'$ of $\gamma$ to $\Pi$ is defined and lies in~$\Omega$.
By construction of $\tilde\Theta$, the curve $\gamma$ and therefore $\gamma'$ are contractible in~$\Omega$.
From two-convexity of $\Omega$
and Proposition~\ref{prop:3d-strong-2-convexity},
the curve $\gamma'$ is contractible in $\Pi\cap \Omega$.

It follows that $\gamma$ is contractible in $\Theta$ 
and therefore $\tilde\gamma$ is contractible in~$\tilde\Theta$.
\qeds

\section{Two-convex hulls}

The following proposition 
describes a construction which produces the two-convex hull $\Conv_2 \Omega$ of an open set $\Omega\subset\EE^3$.
This construction is very close to the one given by Samuel Shefel~\cite{shefel-1964}.

\begin{thm}{Proposition}\label{prop:2-conv-construction}
Let $\Pi_1,\Pi_2\ldots$ be an everywhere dense
sequence of planes in~$\EE^3$.
Given an open set $\Omega$, consider 
the recursively defined sequence of open sets 
$\Omega=\Omega_0\subset\Omega_1\subset\ldots$
such that 
$\Omega_n$ is the union of $\Omega_{n-1}$ 
and all the bounded components of 
$\EE^3\setminus(\Pi_n\cup \Omega_{n-1})$.
Then 
\[\Conv_2\Omega=\bigcup_n\Omega_n.\]

\end{thm}

\parit{Proof.}
Set 
\[\Omega'=\bigcup_n\Omega_n.\eqlbl{eq:Omega'}\]
Note that $\Omega'$ is a union of open sets; in particular, $\Omega'$ is open.

Let us show that
\[\Conv_2\Omega\supset\Omega'.\eqlbl{eq:Conv2>Omega'}\]
Suppose we already know that $\Conv_2\Omega\supset\Omega_{n-1}$. 
Fix a bounded component $\mathfrak{C}$ of $\EE^3\setminus(\Pi_n\cup \Omega_{n-1})$.
It is sufficient to show that $\mathfrak{C}\subset\Conv_2\Omega$.

By \ref{prop:two-hull-open}, $\Conv_2\Omega$ is open.
Therefore, if $\mathfrak{C}\not\subset\Conv_2\Omega$,
then there is a point $p\in \mathfrak{C}\setminus\Conv_2\Omega$ lying at maximal distance from~$\Pi_n$.
Denote by $W_p$ the plane containing $p$ which is parallel to~$\Pi_n$.

Note that $p$ lies in a bounded component of $W_p\setminus \Conv_2\Omega$.
In particular $p$ can be surrounded by a simple closed curve $\gamma$ in $W_p\z\cap\Conv_2\Omega$.
Since $p$ lies at maximal distance from $\Pi_n$,
the curve $\gamma$ is null-homotopic in $\Conv_2\Omega$.
Therefore $p\in \Conv_2\Omega$ ---  a contradiction.

By induction, $\Conv_2\Omega\supset\Omega_n$ for each~$n$.
Therefore \ref{eq:Omega'} implies~\ref{eq:Conv2>Omega'}.

It remains to show that $\Omega'$ is two-convex.
Assume the contrary; 
that is, there is a plane $\Pi$ 
and a simple closed curve $\gamma\:\mathbb{S}^1\to \Pi\cap \Omega'$ 
which is null-homotopic in $\Omega'$,
but not null-homotopic in $\Pi\cap\Omega'$.

By approximation we can assume that $\Pi=\Pi_n$ for a large $n$, and that $\gamma$ lies in $\Omega_{n-1}$.
By the same argument as in the proof of Proposition~\ref{prop:3d-strong-2-convexity} using the loop theorem, we can assume that there is an embedding $\phi\: \DD\to \Omega'$ such that $\phi|_{\partial\DD}=\gamma$ and $\phi(D)$ lies entirely in one of the half-spaces bounded by~$\Pi$.
By the $n$-th step of the construction, the entire bounded domain $U$ bounded by $ \Pi_n$ and $\phi(D)$ is contained in $\Omega'$ and hence $\gamma$ is contractible in $\Pi\cap\Omega'$ --- a contradiction.
\qeds

\begin{thm}{Key lemma}\label{lem:key-shefel}
The two-convex hull of the interior of a polyhedral set in $\EE^3$
is also the interior of a polyhedral set.
\end{thm}

\parit{Proof.}
Fix a polyhedral set $P$ in~$\EE^3$.
Set $\Omega=\Int P$.
We may assume that $\Omega$ is dense in~$P$
(if not, redefine $P$ as the closure of $\Omega$).
Denote by $F_1,\ldots,F_m$ the facets of~$P$.
By subdividing $F_i$ if necessary, we may assume that all $F_i$ are convex polygons.

Set $\Omega'=\Conv_2\Omega$ and let $P'$ be the closure of~$\Omega'$.
Further, let $F'_i=F_i\setminus \Omega'$;
in other words,
$F'_i$ is the subset of the facet $F_i$ 
which remains on the boundary of~$P'$.

From the construction of the two-convex hull (\ref{prop:2-conv-construction}) we have:

\begin{clm}{}\label{clm:F'-convex}
$F'_i$ is a convex subset of~$F_i$.
\end{clm}

Further, since $\Omega'$ is two-convex we obtain the following:

\begin{clm}{}\label{clm:complement-of-F'-convex}
Each connected component of the complement $F_i\setminus F'_i$ is convex.
\end{clm}

{

\begin{wrapfigure}{r}{20mm}
\vskip-0mm
\centering
\includegraphics{mppics/pic-1209}
\end{wrapfigure}

Indeed, assume a connected component $A$ of $F_i\setminus F'_i$ fails to be convex.
Then $F'_i$ has an extreme point in the interior of $F_i$.
In other words, there is a supporting line $\ell$ to $F'_i$ touching $F'_i$ at a single point in the interior of~$F_i$.
Then one could rotate the plane of $F_i$ slightly around $\ell$ and move it parallelly to cut a ``cap'' from the complement of~$\Omega$.
The latter means that $\Omega$ is not two-convex
--- a contradiction.
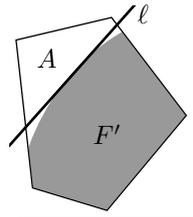
\claimqeds

}

From \ref{clm:F'-convex} and \ref{clm:complement-of-F'-convex}, we conclude
\begin{clm}{}\label{clm:F'} $F'_i$ is a convex polygon for each~$i$.
\end{clm}

Consider the complement 
$\EE^3\setminus \Omega$ 
equipped with the length metric.
By construction of the two-convex hull (\ref{prop:2-conv-construction}), 
the complement $L\z=\EE^3\setminus (\Omega'\cup P)$
is locally convex;
that is, any point of $L$ admits a convex neighborhood.

Summarizing: (1)
$\Omega'$ is a two-convex open set,
(2) the boundary $\partial\Omega'$ 
contains a finite number of polygons $F_i'$
and the remaining part $S$ of the boundary is locally concave.
It remains to do the following exercise.

\begin{thm}[!]{Advanced exercise}\label{ex:convex+saddle+broken=>PL}
Let $\bar S$ be the closure of $S$.
Denote by $\partial S$ the boundary of $\bar S$ in $\partial \Omega'$.

\begin{subthm}{ex:convex+saddle+broken=>PL:a}
Show that (1) and (2) imply that any point $p\in S$ lies in a line segment in $\bar S$ with ends in $\partial S$ or in a polygon in $\bar S$ with vertices in~$\partial S$.
\end{subthm}

\begin{subthm}{ex:convex+saddle+broken=>PL:b}
Use part \ref{SHORT.ex:convex+saddle+broken=>PL:a} to show that $\bar S$ (and therefore $\partial\Omega'$)
is piecewise linear. \qeds
\end{subthm}

\end{thm}

\section{Proof of Shefel's theorem}

\parit{Proof of \ref{thm:shefel}.}
The ``only if'' part can be proved in the same way as in the smooth two-convexity theorem (\ref{thm:set-with-smooth-bry:CBA}) with the additional use of Exercise~\ref{ex:inner-support}.

\parit{If part.}
Suppose $\Omega$ is two-convex. 
We need to show that $\tilde K$ is $\CAT(0)$.

Fix a quadruple of points $x_1,x_2,x_3,x_4\in\tilde \Omega$.
Let us show that
$\CAT(0)$ comparison holds for this quadruple.

Fix $\eps>0$.
Choose six polygonal lines in $\tilde \Omega$ connecting all pairs of points $x_1,x_2,x_3,x_4$, where the length of each polygonal line is at most $\eps$ bigger than
the distance between its ends in the length metric on~$\tilde \Omega$.
Denote by $X$ the union of these polygonal lines.
Choose a polyhedral set $P$ in $\Omega$ such that its interior $\Int P$ contains the projections of all six polygonal lines and discs which contract all the loops created by them (it is sufficient to take 3 discs).

Denote by $\Omega'$ the two-convex hull of the interior of~$P$.
According to the key lemma (\ref{lem:key-shefel}), $\Omega'$ is the interior of a polyhedral set.

Equip $\Omega'$ with the induced length metric.
Consider the universal metric cover $\tilde\Omega'$ of~$\Omega'$.
(The covering $\tilde\Omega'\to\Omega'$ might be nontrivial ---
even if $\Int P$ is simply connected, its two-convex hull $\Omega'$ might not be simply connected.)
Denote by $\tilde K'$ the completion of~$\tilde\Omega'$.

By Lemma~\ref{lem:poly-shefel}, $\tilde K'$ is $\CAT(0)$.

By construction of $\Int P$, the embedding $\Int P\hookrightarrow \Omega'$
admits a lift $\iota\:X\hookrightarrow \tilde K'$, and $\iota$ almost preserves the distances between the points $x_1,x_2,x_3,x_4$;
namely 
\begin{align*}
\dist{\iota(x_i)}{\iota(x_j)}{\tilde K'}\lessgtr \dist{x_i}{x_j}{\Int P}\pm\eps.
\end{align*}

Since $\eps>0$ is arbitrary and $\CAT(0)$ comparison holds in $\tilde K'$,
we get that $\CAT(0)$ comparison holds for $x_1,x_2,x_3,x_4$ in $\tilde \Omega$.

The statement follows since the quadruple $x_1,x_2,x_3,x_4\in\tilde\Omega$ is arbitrary.
\qeds

\begin{thm}{Exercise}\label{ex:CAT=>two-convex}\label{ex:two-convex-not-a-CAT}
Assume $K\subset \EE^{m}$ is a subgraph of a Lipschitz function $f\:\RR^{m-1}\to\RR$.
Equip $K$ with the induced length metric.
Show that if $K$ is $\CAT(0)$, then $K$ is two-convex.

Show that the converse does not hold.
\end{thm}

\section{Notes}

Under the name {}\emph{$(n-2)$-convex sets}, 
two-convex sets in $\EE^n$ were introduced by Mikhael Gromov \cite{gromov-1991}.
In addition to the inheritance of upper curvature bounds by two-convex sets discussed in this lecture, 
these sets appear as the maximal open sets with vanishing curvature in Riemannian manifolds with non-negative or non-positive sectional curvature [see Lemma 5.8 in \ncite{buyalo}, \ncite{andersson-howard}
and \ncite{panov-petrunin}].

Two-convex sets could be defined using homology instead of homotopy, as in Gromov's formulation of the Lefschetz theorem \cite[\S\textonehalf]{gromov-1991}.
Namely, we can say that $K$ is two-convex if the following condition holds: \textit{if a 1-dimensional cycle $z$ has support in the intersection of $K$ with a plane $W$ and bounds in $K$, then it bounds in $K\cap W$}.

The resulting definition is equivalent to the one used above.
But unlike our definition it can be generalized to define $k$-convex sets in $\EE^m$ for~$k>2$.
With this homological definition one can also avoid the use of the loop theorem, whose proof is quite involved.
We stick to the homotopy-based definition since it is slightly easier to visualize.

Both definitions work well for open sets; for general sets one should be able to give a similar definition using  an appropriate homotopy/ho\-mo\-logy theory.

Proposition~\ref{prop:two-cove+smooth} was proved by Mikhael Gromov \cite[\S\textonehalf]{gromov-1991}, but we added a few details.

Jointly with David Berg and Richard Bishop, the first author proved the following theorem, which gives the exact upper bound on Alexandrov's curvature for the Riemannian manifolds with boundary.
This theorem includes the smooth two-convexity theorem (\ref{thm:set-with-smooth-bry:CBA}) as a partial case.

\begin{thm}{Theorem}
Let $M$ be a Riemannian manifold with boundary~$\partial M$.
A direction tangent to the boundary will be called concave if there is a short geodesic in this direction which leaves the boundary and goes into the interior of~$M$.
A sectional direction (that is, a 2-plane) 
tangent to the boundary 
will be called {}\emph{concave} if all the directions in it are concave.

Denote by $\kappa$ an upper bound of sectional curvatures of $M$ and  
sectional curvatures of $\partial M$ in the concave sectional directions. 
Then $M$ is locally $\CAT(\kappa)$. 
\end{thm}

\begin{thm}{Corollary}
Let $M$ be a Riemannian manifold with boundary~$\partial M$. 
Assume that all the sectional curvatures of $M$ and $\partial M$ are bounded above by~$\kappa$.
Then $M$ is locally $\CAT(\kappa)$.
\end{thm}

The Theorem~\ref{thm:bishop-plane} was proved by Richard Bishop \cite{bishop}.
Generalizations of this result were proved by Alexander Lytchak, Stephan Stadler and Stefan Wenger; see \cite{lytchak-stadler} and \cite[Proposition 12.1]{lytchak-wenger}.

Samuel Shefel was intrigued by the survival of metric properties under affine transformation.
An instance of this can be seen in Corollary~\ref{cor:affine}.

Theorem~\ref{thm:shefel-graph} is from Shefel's original paper~\cite{shefel-1965}.
It is related to Alexandrov's theorem about ruled surfaces \cite{alexandrov-1957-ruled-surfaces}.

The following question remains open.

\begin{thm}{Shefel's conjecture}
Suppose $s\:\DD \to\EE^3$ is a Lipschitz embedding of the disc $\DD$,
and its image $D=s(\DD)$ is \index{saddle surface}\emph{saddle}; that is, for any plane $\Pi$ each connected component of $D\setminus \Pi$ meets the boundary of $D$.
Then $D$ equipped with the length-metric is locally $\CAT(0)$.
\end{thm}

A partial answer is given in Theorem \ref{thm:shefel-graph},
which is a result of Shefel~\cite{shefel-1964}.
He also solved an analogous question in $\EE^2$; see \cite{shefel-1965}.
Applying Crofton-type argument to his result, one can see that a saddle disc satisfies the isoperimetric inequality $a\le C\cdot \ell^2$
where $a$ is the area of a disc bounded by a curve of length $\ell$ and $C=\tfrac{1}{3\cdot\pi}$.
By a result of Alexander Lytchak and Stefan Wenger \cite{lytchak-wenger}, Shefel's conjecture would follow if the constant could be made optimal; that is, $C=\tfrac{1}{4\cdot\pi}$.
Several variations and generalizations of this conjecture are discussed by the third author and Stephan Stadler \cite{petrunin-stadler}.

%%!TEX root = the-kleiner.tex
\chapter{Barycenters}\label{chap:barycenter}

Barycenters provide yet another way to use your Euclidean intuition for $\CAT(0)$ spaces.
We will apply them to dimension theory.

\section{Definition}

Let us denote by $\triangle^k\subset \RR^{k+1}$\index{1@$\triangle^m$ (simplex)}
the \index{standard simplex}\emph{standard $k$-simplex}; 
that is, $\bm{\mu}=(\mu_0,\ldots,\mu_k)\in\triangle^k$ if $\mu_0+\ldots+\mu_k=1$ and $\mu_i\ge0$ for all $i$.

Consider a point array $\bm{p}=(p_0,\ldots,p_k)$ in the Euclidean space $\EE^n$.
Recall that
\[z=\mu_0\cdot p_0+\ldots+\mu_k\cdot p_k\]
is called the \index{barycenter}\emph{barycenter} of the point array $\bm{p}=(p_0,\ldots,p_k)$ with masses $\bm{\mu}\z=(\mu_0,\ldots, \mu_k)\in \triangle^k$.
Equivalently, 
\[z\df \argmin (\mu_0\cdot f_0+\ldots+\mu_k\cdot f_k),\eqlbl{eq:f-spx}\]
where $f_i=\tfrac12\cdot\distfun_{p_i}^2$ for each $i$, and $\argmin f$\index{3@$\argmin$ (point of minimum)} denotes a point of minimum of the function $f$.

The map $\spx{}\:\triangle^k\mapsto \EE^n$ defined by $\spx{}\:\bm{\mu}\mapsto z$ is called the \index{barycentric simplex} \emph{barycentric simplex} of the array $\bm{p}$.
If needed we may denote this map by $\spx{\bm{p}}$ or, more generally, $\spx{\bm{f}}$.
The latter means that we define the map via \ref{eq:f-spx} for an array of functions $\bm{f}=(f_0,\ldots,f_k)$.
Formally speaking, the definition \ref{eq:f-spx} makes sense for any array of functions in a metric space;
but in this case, the map might be undefined or non-uniquely defined.

Further, we will work with this definition in $\CAT(0)$ spaces.
We will essentially use only that on a geodesic $\CAT(0)$ space
functions of the type $f=\tfrac12\cdot\distfun_{p}^2$ are 1-convex; see \ref{ex:convex-distfun}.

\section{Barycentric simplex}

\begin{thm}{Theorem}\label{thm:barycenter}
Let $\spc{X}$ be a complete geodesic space 
and $\bm{f}\z=(f_0,\z\ldots,f_k)\:\spc{X}\to\RR^{k+1}$
be an array of non-negative 1-convex locally Lipschitz functions.
Then the barycentric simplex 
$\spx{\bm{f}}\:\triangle^k\to \spc{X}$ is a uniquely defined Lipschitz map.

In particular, we have that for any point array $\bm{p}\z=(p_0,\z\ldots,p_k)$ in a complete geodesic $\CAT(0)$ space
the barycentric simplex $\spx{\bm{p}}$ is a uniquely defined Lipschitz map.
\end{thm}

\begin{thm}{Lemma}\label{lem:argmin(convex)}
Suppose $\spc{X}$ is a complete geodesic space and $f\:\spc{X}\z\to\RR$ is a locally Lipschitz, $1$-convex function.
Then $\argmin f$ is uniquely defined.
\end{thm}

\parit{Proof.}
Note that
\begin{clm}{}\label{midpoint}
 if $z$ is a midpoint of a geodesic $[x y]$, then 
\[s\le f(z)
\le
\tfrac{1}{2}\cdot f(x)+\tfrac{1}{2}\cdot f(y)-\tfrac{1}{8}\cdot\dist[2]{x}{y}{},
\]
where $s$ is the infimum of $f$.
\end{clm}

\parit{Uniqueness.}
Assume that $x$ and $y$ are distinct minimum points of $f$.
Let $z$ be the midpoint of a geodesic $[x y]$.
From \ref{midpoint} we have
\[f(z)<f(x)=f(y)\] 
--- a contradiction. 

\parit{Existence.}
Fix a point $p\in \spc{X}$, and
let $\Lip\in\RR$ be a Lipschitz constant of $f$ in a neighborhood of $p$.

Choose a sequence of points $p_n\in \spc{X}$ such that $f(p_n)\to s$.
Applying \ref{midpoint} for $x\z=p_n$, $y\z=p_m$, we see that $p_n$ is a Cauchy sequence.
Thus the sequence $p_n$ converges to a minimum point of $f$.
\qeds

Assume $\phi$ is a convex function defined on a real interval $\II$.
For $t_0\in\II$, let us define the \index{right derivative}\emph{right} (respectively, \index{left derivative}\emph{left}) \index{3@$\phi^\pm$ (right/left derivative)}\emph{derivative} $\phi^+(t_0)$ ($\phi^-(t_0)$) at $t_0$ by
\[\phi^\pm(t_0)=\lim_{t\to t_0\pm} \frac{\phi(t)-\phi(t_0)}{|t-t_0|}.\]
Note that our sign convention for $\phi^-$ is not standard --- for $\phi(t)=t$ we have
$\phi^+(t)=1$ and $\phi^-(t)=-1$.

\parit{Proof of \ref{thm:barycenter}.}
Since each $f_i$ is $1$-convex, for any $\bm{\lambda}=(\lambda_0,\ldots,\lambda_k)\z\in\triangle^k$
the convex combination 
\[\left(\sum_i \lambda_i\cdot f_i\right)\:\spc{X}\to\RR\]
is also $1$-convex.
Therefore, according to \ref{lem:argmin(convex)}, the barycentric simplex 
%$\spx{\bm{f}}(\bm{\lambda})$ is defined for any $\bm{\lambda}\in\triangle^k$.
$\spx{\bm{f}}$ is uniquely defined on $\triangle^k$.
 
For $\bm{\lambda},\bm{\mu}\in\triangle^k$,
let 
\begin{align*}
f_{\bm{\lambda}}
&=\sum_i \lambda_i\cdot f_i,
&
f_{\bm{\mu}}
&=\sum_i \mu_i\cdot f_i,
\\
p
&=\spx{\bm{f}}(\bm{\lambda}),
&
q
&=\spx{\bm{f}}(\bm{\mu}),
\end{align*}
Choose a geodesic $\gamma$ from $p$ to $q$;
suppose $s=\dist{p}{q}{}$ and so $\gamma(0)=p$ and $\gamma(s)=q$.
Observe the following:
\begin{itemize}
\item The function $\phi(t)=f_{\bm{\lambda}}\circ\gamma(t)$ has a minimum at $0$.
Therefore $\phi^+(0)\ge 0$.

\item The function $\psi(t)=f_{\bm{\mu}}\circ\gamma(t)$ has a minimum at $s$.
Therefore $\psi^-(s)\ge 0$.
\end{itemize}
From $1$-convexity of $f_{\bm{\mu}}$, we have
\[\psi^+(0)+\psi^-(s)+s\le0.\]

Let $\Lip$ be a Lipschitz constant for all $f_i$ in a neighborhood $\Omega\ni p$.
Then 
\[\psi^+(0)
\le 
\phi^+(0)+\Lip\cdot\|\bm{\lambda}-\bm{\mu}\|_1,\]
where $\|\bm{\lambda}-\bm{\mu}\|_1=\sum_{i=0}^k|\lambda_i-\mu_i|$.
It follows that given $\bm{\lambda}\in\triangle^k$, there is a constant $\Lip$ such that
\begin{align*}
\dist{\spx{\bm{f}}(\bm{\lambda})}{\spx{\bm{f}}(\bm{\mu})}{}
&=
s
\le
\\
&\le 
\Lip\cdot\|\bm{\lambda}-\bm{\mu}\|_1
\end{align*}
for any $\bm{\mu}\in\triangle^k$.
In particular, there is $\eps>0$ such that if $\|\bm{\lambda}-\bm{\mu}\|_1<\eps,$ $\|\bm{\lambda}-\bm{\nu}\|_1 <\eps$, then $\spx{\bm{f}}(\bm{\mu})$, $\spx{\bm{f}}(\bm{\nu})\in\Omega$.
Thus the same argument as above implies 
\[\dist{\spx{\bm{f}}(\bm{\mu})}{\spx{\bm{f}}(\bm{\nu})}{}
\le \Lip\cdot\|\bm{\mu}-\bm{\nu}\|_1\]
for any $\bm{\mu}$ and $\bm{\nu}$ sufficiently close to $\bm{\lambda}$; that is, $\spx{\bm{f}}$ is locally Lipschitz.
Since $\triangle^k$ is compact, $\spx{\bm{f}}$ is Lipschitz.
\qeds

\begin{thm}{Exercise}\label{ex:finite-action-CAT}
Let $G$ be a subgroup of the group of isometries of a proper geodesic $\CAT(0)$ space.
Assume that
\begin{subthm}{ex:finite-action-CAT:finite}
$G$ is finite, or, more generally, 
\end{subthm}
\begin{subthm}{ex:finite-action-CAT:compact}
 $G$ is compact.
\end{subthm}
Show that the action of $G$ has a fixed point.
\end{thm}

\begin{thm}{Advanced exercise}\label{ex:bary-jensen}
Let $z$ be the barycenter of points $p_0,\ldots,p_k$ equipped with masses $\mu_0,\ldots,\mu_k$; we assume that $\mu_i\ge0$ for all $i$,
$\mu_0+\ldots+\mu_k=1$, and the points $p_0,\ldots,p_k$ are in a complete geodesic $\CAT(0)$ space $\spc{U}$.

\begin{subthm}{ex:bary-jensen:moment}
Show that
\[\sum_i \mu_i\cdot \dist[2]{q}{p_i}{}\ge \dist[2]{z}{q}{} +\sum_i \mu_i\cdot \dist[2]{z}{p_i}{}\]
for any $q\in \spc{U}$.
\end{subthm}

\begin{subthm}{ex:bary-jensen:moment+}
 Show that
\[\sum_i
\mu_i\cdot \dist[2]{z}{p_i}{}
\le
\sum_{i<j}
\mu_i\cdot\mu_j\cdot \dist[2]{p_i}{p_j}{}.\]
\end{subthm}

\begin{subthm}{ex:bary-jensen:jensen}
(Jensen inequality.)
Try to show that
\[h(z)\le \mu_0\cdot h(p_0)+\ldots+\mu_k\cdot h(p_k)\]
for any convex function $h\:\spc{U}\to\RR$.
\end{subthm}

\end{thm}

\section{Convexity of up-set}

\begin{thm}{Definition}\label{def:ordung}
For two real arrays $\bm{v}$, $\bm{w}\in \RR^{\kay+1}$,
$\bm{v}=(v_0,\z\ldots,v_\kay)$
and 
$\bm{w}=(w_0,\z\ldots,w_\kay)$,
we will write
$\bm{v}\succcurlyeq\bm{w}$ if $v_i\ge w_i$ for each $i$.
\end{thm}

Given a subset $Q\subset \RR^{\kay+1}$, 
denote by $\Up Q$ \label{PAGE.def:Up}
the smallest upper set containing $Q$;
that is,
\begin{align*}
\Up Q 
&=
\set{\bm{v}\in\RR^{\kay+1}}{\exists\, \bm{w}\in Q\ \text{such that}\ \bm{v}\succcurlyeq\bm{w}},
\end{align*}

\begin{thm}{Proposition}\label{thm:up-convex}
Let $\spc{X}$ be a complete geodesic space 
and $\bm{f}\z=(f_0,\z\ldots,f_\kay)\:\spc{X}\to\RR^{\kay+1}$
be an array of non-negative 1-convex locally Lipschitz functions.
Consider the set $W=\Up[\bm{f}(\spc{X})]\subset\RR^{k+1}$.
Then

\begin{subthm}{thm:up-convex:convex}
The set $W$ is convex.
\end{subthm}

\begin{subthm}{thm:up-convex:bry}
$\bm{f}[\spx{\bm{f}}(\triangle^k)]\subset\partial W$.
Moreover, $\bm{f}[\spx{\bm{f}}(\triangle^k)\setminus \spx{\bm{f}}(\partial\triangle^k)]$ is an open set in $\partial W$.
\end{subthm}

\begin{subthm}{thm:up-convex:bry+}
$W=\Up(\bm{f}[\spx{\bm{f}}(\triangle^k)])$; 
in other words, $\Up(\bm{f}[\spx{\bm{f}}(\triangle^k)])\supset\bm{f}(\spc{X})$.
In particular, $W$ is closed.
\end{subthm}

\end{thm}

\begin{wrapfigure}{r}{48mm}
\vskip0mm
\centering
\includegraphics{mppics/pic-1106}
\end{wrapfigure}

\parit{Proof.}
Let $V=\bm{f}(\spc{X})\subset\RR^{k+1}$; so $W\z=\Up V$.
Denote by $\bar V$ the closure of $V$.

\parit{\ref{SHORT.thm:up-convex:convex}.}
Convexity of all $f_i$ implies that
for any two points $p,q\in \spc{X}$ and $t\z\in[0,1]$ we have
\[(1-t)\cdot\bm{f}(p)+t\cdot \bm{f}(q)
\succcurlyeq
\bm{f}\circ\gamma(t),\]
where $\gamma$ denotes a geodesic path from $p$ to $q$. 
Therefore, $W$ is convex.

\parit{\ref{SHORT.thm:up-convex:bry}+\ref{SHORT.thm:up-convex:bry+}.}
Choose $p\in \spx{\bm{f}}(\triangle^k)$.
Note that if $\bm{f}(p)\succcurlyeq\bm{w}$ for some $\bm{w}\in W$, then $\bm{f}(p)\succcurlyeq\bm{w}\succcurlyeq\bm{f}(x)$ for some $x\in \spc{X}$ and hence $x=p$ by the uniqueness of $\spx{\bm{f}}$.
Therefore, $\bm{f}(p)=\bm{w}$.
It follows that $\bm{f}(p)\in\partial W$;
therefore $\bm{f}[\spx{\bm{f}}(\triangle^k)]$ lies in a convex hypersurface $\partial W$.

Choose $\bm{w}\in W$.
Observe that $\bm{w}\succcurlyeq\bm{v}$ for some $\bm{v}\in \bar V\cap \partial W$.
Note that $W$ is supported at $\bm{v}$ by a hyperplane 
\[\Pi=\set{(x_0,\ldots,x_k)\in \RR^{k+1}}{\mu_0\cdot x_0+\ldots+\mu_k\cdot x_k=\const}\]
for some $\bm{\mu}=(\mu_0,\ldots,\mu_k)\in\triangle^{\kay}$.
Let $p=\spx{\bm f}(\bm{\mu})$.
By \ref{lem:argmin(convex)}, $\bm{f}(p)\z=\bm{v}$;
in particular $\bm{v}\in V$.

Note that $p\in \spx{\bm{f}}(\triangle^k)\setminus\spx{\bm{f}}(\partial\triangle^k)$ if and only if 
$\bm{f}(p)$ is supported by a plane as above for some $\bm{\mu}\in\triangle^{\kay}$,
but not supported by any such plane with $\bm{\mu}\in\partial\triangle^{\kay}$.
This condition is open, therefore $\spx{\bm{f}}(\triangle^k)\z\setminus\spx{\bm{f}}(\partial\triangle^k)$ is an open set.

Since $\triangle^k$ is compact, $W=\Up(\bm{f}[\spx{\bm{f}}(\triangle^k)])$ is closed.
\qeds

\section{Nondegenerate simplex}

Given an array $\bm{f}=(f_0,\ldots,f_k)$,
we denote by $\bm{f}^{\without i}$ the subarray of $\bm{f}$ with $f_i$ removed;
that is, 
\[\bm{f}^{\without i\,}
\df
(f_0,\ldots,f_{i-1},f_{i+1},\ldots,f_k).\]
Note that $\spx{\bm{f}^{\without i}}$ coincides with the restriction of $\spx{\bm{f}}$ to the corresponding facet of $\triangle^k$.

If $\Im \spx{\bm{f}}$ is not covered by the union of $\Im \spx{\bm{f}^{\without i}}$ over all $i$, then we say that $\spx{\bm{f}}$ is nondegenerate.
In other words, $\spx{\bm{f}}$ is \index{nondegenerate simplex}\emph{nondegenerate} if 
\[\spx{\bm{f}}(\triangle^k)\setminus \spx{\bm{f}}(\partial \triangle^k)\ne\emptyset.\]

\begin{thm}{Exercise}\label{ex:barysimple}
Let $\spc{U}$ be a complete geodesic $\CAT(0)$ space.

Show that the image of the 1-dimensional barycentric simplex for a pair of points $p_0,p_1\in\spc{U}$ is the geodesic $[p_0\,p_1]$.

Construct a $\CAT(0)$ space with a three-point array $(p_0,p_1,p_2)$ such that its barycentric simplex is nondegenerate and noninjective.
\end{thm}

\begin{thm}[!]{Exercise}\label{lem:nondeg-test-with-balls}
Let $\bm{p}\z=(p_0,\ldots, p_\kay)$ be a point array
in a complete length $\CAT(0)$ space $\spc{U}$, 
and $B_i=\cBall[p_i,r_i]$ for some array of positive reals $(r_0,r_1,\ldots,r_\kay)$.

\begin{subthm}{lem:nondeg-test-with-ball:U}
Suppose $\bigcap_i B_i\ne \emptyset$.
Show that 
\[\Im\spx{\bm{p}}\subset \bigcup_i B_i.\]
\end{subthm}

\begin{subthm}{lem:nondeg-test-with-balls:nondeg}
Suppose $\bigcap_i B_i= \emptyset$, but $\bigcap_{i\ne j} B_i\ne \emptyset$ for any $j$.
Show that 
$\spx{\bm{p}}$ is nondegenerate.
\end{subthm}

\begin{subthm}{lem:nondeg-test-with-balls:nondeg+}
Suppose $\spx{\bm{p}}$ is nondegenerate.
Show that the condition in \ref{SHORT.lem:nondeg-test-with-balls:nondeg} holds for some
 array of positive reals $(r_0,\ldots,r_\kay)$.
\end{subthm}

\end{thm}

\section{bi-H\"older embedding}

\begin{thm}{Theorem}\label{thm:bihoelder}
Let $\spc{X}$ be a complete geodesic space 
and $\bm{f}\z=(f_0,\z\ldots,f_k)\:\spc{X}\to\RR^{k+1}$
be an array of 1-convex locally Lipschitz functions.
Then the set 
\[Z=\spx{\bm{f}}(\triangle^k)\setminus \spx{\bm{f}}(\partial \triangle^k)\]
is $C^{\frac12}$-bi-H\"older to an open domain in $\RR^k$.
\end{thm}

\parit{Proof.}
Let $\proj\:\RR^{k+1}\to\Pi$ be orthogonal projection to the hyperplane $x_0+\ldots+x_k=0$.
We will show that the restriction $\proj\circ \bm{f}|_Z$ is a bi-H\"older embedding.

Since the map $\proj\circ \bm{f}$ is Lipschitz,
it is sufficient to construct its right inverse and show that it is $C^{\frac12}$-continuous.

Given $\bm{v}=(v_0,\ldots,v_\kay)\in\Pi$, consider the function
$h_{\bm{v}}\: \spc{X}\to \RR$ defined by
\[h_{\bm{v}}(p)=\max_i\{f_i(p)-v_i\}.\]
Note that $h_{\bm{v}}$ is $1$-convex.
Let 
$$\map(\bm{v})\df\argmin h_{\bm{v}}.$$
According to Lemma~\ref{lem:argmin(convex)}, $\map(\bm{v})$ is uniquely defined.

If $\bm{v}=\proj \bm{f}(p)$, then 
\[f_i\circ \map(\bm{v})\le f_i(p)\]
for any $i$.
In particular, if $p\in \spx{\bm{f}}(\triangle^k)$, then $p=\map(\bm{v})$.
That is, $\map$ is a right inverse of the restriction $\bm{f}|_{\spx{\bm{f}}(\triangle^k)}$.

Given $\bm{v},\bm{w}\in\RR^{\kay+1}$,
set $p=\map (\bm{v})$ and $q=\map (\bm{w})$.
Since $h_{\bm{v}}$ and $h_{\bm{w}}$ are 1-convex, we have
\begin{align*}
h_{\bm{v}}(q)
&\ge 
h_{\bm{v}}(p)+\tfrac{1}{2}\cdot\dist[2]{p}{q}{},
&
h_{\bm{w}}(p)
&\ge 
h_{\bm{w}}(q)+\tfrac{1}{2}\cdot\dist[2]{p}{q}{}.
\end{align*}
Therefore,
\begin{align*}
\dist[2]{p}{q}{}
&\le 
2\cdot\sup_{x\in\spc{X}}\{ |h_{\bm{v}}(x)-h_{\bm{w}}(x)| \}
\le
\\
&\le 
2\cdot\max_{i}\{|v_i-w_i|\}.
\end{align*}
In particular,
$\map$ is $C^{\frac{1}{2}}$-continuous.

Finally, by \ref{thm:up-convex:bry}, $\bm{f}(Z)$ is a $k$-dimensional manifold --- hence the result.
\qeds

\section{Review of topological dimension}

Let $\spc{X}$ be a metric space and $\set{V_\beta}{\beta\in\mathcal{B}}$
 be an open cover of $\spc{X}$.
Let us recall two notions in general topology:
\begin{itemize}

\item The \index{order of a cover}\emph{order} of $\set{V_\beta}{\beta\in\mathcal{B}}$ is the supremum of all integers $n$ such that there is a collection of $n+1$ elements of $\{V_\beta\}$ with nonempty intersection.

\item An open cover $\set{W_\alpha}{\alpha\in\IndexSet}$ of $\spc{X}$ is called a \index{refinement of a cover}\emph{refinement} of $\set{V_\beta}{\beta\in\mathcal{B}}$ if for any $\alpha\in\IndexSet$ there is $\beta\in\mathcal{B}$ such that $W_\alpha\subset V_\beta$.

\end{itemize}

\begin{thm}{Definition}\label{def:TopDim}\index{dimension!topological dimension}\index{topological dimension}
Let $\spc{X}$ be a metric space. 
The topological dimension of $\spc{X}$ is defined to be the minimum of non-negative integers $n$
such that for any open cover of $\spc{X}$ there is a locally finite open refinement with order~$n$.
If no such $n$ exists, then the topological dimension of $\spc{X}$ is infinite.

The topological dimension of $\spc{X}$ will be denoted by $\TopDim\spc{X}$.
\end{thm}

The invariants satisfying statements \ref{dim-axiom-norm} and \ref{dim-axiom-sigma} are commonly called the \index{dimension}\emph{dimension};
for that reason, we call these statements axioms.

\begin{thm}{Normalization axiom}
\label{dim-axiom-norm} For any $m\in\ZZ_{\ge0}$,
\[\TopDim\EE^m=m.\]

\end{thm}

\begin{thm}{Cover axiom}\label{dim-axiom-sigma} 
If $\{A_n\}_{n=1}^\infty$ is a countable closed cover of $\spc{X}$, then
\begin{align*}
\TopDim \spc{X}&=\sup\nolimits_n\{\TopDim A_n\}.
\end{align*}

\end{thm}

\parbf{On product spaces.} 
The following inequality holds for arbitrary metric spaces
\begin{align*}
\TopDim (\spc{X}\times\spc{Y})
&\le 
\TopDim \spc{X}+ \TopDim\spc{Y}
\end{align*}
and this inequality might be strict \cite{pontyagin-surface}.

\medskip

\begin{thm}{Definition}
Let $\spc{X}$ be a metric space
and $F\:\spc{X}\to\RR^m$ be a continuous map.
A point $\bm{z}\in \Im F$ is called a \emph{stable value} of $F$
if there is $\eps>0$ such that $\bm{z}\in\Im F'$ 
for any \emph{$\eps$-close} to $F$ continuous map $F'\:\spc{X}\to\RR^m$,
that is, $|F'(x)-F(x)|<\eps$ for all $x\in \spc{X}$.
\end{thm}

The next theorem follows from \cite[theorems VI 1$\&$2]{hurewicz-wallman}.
(This theorem also holds for non-separable metric spaces \cite{nagata}, \cite[3.2.10]{engelking}). 

\begin{thm}{Stable value theorem}\label{thm:stable-value}
Let $\spc{X}$ be a separable metric space.
Then $\TopDim\spc{X}\ge m$ if and only if there is a map $F\:\spc{X}\to\RR^m$ with a stable value.
\end{thm}

\section{Dimension theorem}\label{sec:dim-CAT}

\begin{thm}{Theorem}\label{thm:dim-infty-CBA}
For any proper geodesic $\CAT(0)$ space $\spc{U}$, the following statements are equivalent:

\begin{subthm}{thm:dim-infty-CBA:TopDim}
\[\TopDim \spc{U}\ge m.\]
\end{subthm}

\begin{subthm}{thm:dim-infty-CBA:bary} 
For some $z\in \spc{U}$ there is an array of $m+1$ balls $B_i\z=\oBall(a_i,r_i)$
such that 
\[\bigcap_i B_i=\emptyset,
\quad\text{but}\quad
\bigcap_{i\ne j} B_i\ne \emptyset
\quad \text{for each $j$}.\]

\end{subthm}

\begin{subthm}{thm:dim-infty-CBA:mnfld} 
There is a $C^{\frac{1}{2}}$-embedding of an open set in $\RR^m$ into $\spc{U}$;
that is, $\map$ is bi-Hölder with exponent $\tfrac{1}{2}$.
\end{subthm}

\end{thm}

\begin{thm}{Lemma}\label{lem:approximation-cba}
Let $\spc{U}$ be a proper geodesic $\CAT(0)$ space
and $\rho\:\spc{U}\z\to\RR$ be a continuous positive function.
Then there is a locally finite countable simplicial complex $\spc{N}$,
a locally Lipschitz map $\map\:\spc{U}\z\to \spc{N}$, 
and a Lipschitz map $\Psi\:\spc{N}\to\spc{U}$ such that:

\begin{subthm}{lem:approximation-cba:displacement}
The displacement of the composition $\Psi\circ\map\:\spc{U}\to\spc{U}$ is bounded by $\rho$;
that is,
\[\dist{x}{\Psi\circ\map(x)}{}<\rho(x)\] 
for any $x\in\spc{U}$.
\end{subthm}

\begin{subthm}{lem:approximation-cba:im}
If $\TopDim\spc{U}\le m$, 
then the $\Psi$-image of $\spc{N}$ 
coincides with the image of its $m$-skeleton.
\end{subthm}

\end{thm}

\parit{Proof.}
Choose a locally finite countable covering $\set{\Omega_\alpha}{\alpha\in\IndexSet}$ of $\spc{U}$ such that $\Omega_\alpha\subset \oBall(x,\tfrac{1}{3}\cdot\rho(x))$ for any $x\in \Omega_\alpha$. 

Denote by $\spc{N}$ the \index{nerve}\emph{nerve} of the covering $\{\Omega_\alpha\}$;
that is, $\spc{N}$ is an abstract simplicial complex with 
%set of vertices formed by $\IndexSet$,
vertex set $\IndexSet$,
such that a finite subset 
$\{\alpha_0,\z\ldots,\alpha_n\}\subset\IndexSet$
forms a simplex if and only if
\[\Omega_{\alpha_0}
\cap
\ldots\cap
\Omega_{\alpha_n}\ne\emptyset.\]

Choose a Lipschitz partition of unity 
$\phi_\alpha\:\spc{U}\to [0,1]$ subordinate to $\{\Omega_\alpha\}$.
Consider the map $\map\:\spc{U}\to \spc{N}$ such that the barycentric coordinate of $\map(p)$ is $\phi_\alpha(p)$.
Note that $\map$ is locally Lipschitz. 
Clearly, the $\map$-preimage of any open simplex in $\spc{N}$ lies in $\Omega_\alpha$ for some $\alpha\in\IndexSet$.

For each $\alpha\in\IndexSet$, 
choose $x_\alpha\in\Omega_\alpha$.
Let us extend the map $\alpha\mapsto x_\alpha$
to a map $\Psi\:\spc{N}\to\spc{U}$ that is barycentric on each simplex.
According to \ref{thm:barycenter}, this extension exists, 
and $\Psi$ is locally Lipschitz.

\parit{\ref{SHORT.lem:approximation-cba:displacement}.}
Fix $x\in\spc{U}$.
Denote by $\triangle$ the minimal simplex that contains $\map(x)$, 
and let $\alpha_0,\alpha_1,\ldots,\alpha_n$ be the vertices of $\triangle$.
%Denote by $\triangle$ the minimal simplex that contains $\map(x)$;
%and let $(\alpha_0,\alpha_1,\ldots,\alpha_n)$ be the vertices of $\triangle$.
Note that $\alpha$ is a vertex of $\triangle$ if and only if $\phi_{\alpha}(x)>0$.
Thus
\[\dist{x}{x_{\alpha_i}}{}<\tfrac{1}{3}\cdot\rho(x)\] 
for any $i$.
Therefore 
\[\diam\Psi(\triangle)
\le
\max_{i,j}\{\dist{x_{\alpha_i}}{x_{\alpha_j}}{}\}
<
\tfrac{2}{3}\cdot\rho(x).\]
In particular, 
\[\dist{x}{\Psi\circ\map(x)}{}\le\dist{x}{x_{\alpha_0}}{}+\diam \Psi(\triangle) <\rho(x).\]

\parit{\ref{SHORT.lem:approximation-cba:im}.}
Assume the contrary;
that is, $\Psi(\spc{N})$ is not contained in the $\Psi$-image of the $m$-skeleton of $\spc{N}$.
Then for some $\kay>m$,
there is a $\kay$-simplex $\triangle^\kay$ in $\spc{N}$
such that the barycentric simplex $\sigma=\Psi|_{\triangle^\kay}$ is nondegenerate; 
that is, 
$$W=\Psi(\triangle^\kay)\setminus\Psi(\partial\triangle^\kay)\ne \emptyset.
$$
By \ref{thm:bihoelder}, $\TopDim\spc{U}\ge \kay$ --- a contradiction.
\qeds

\parit{Proof of \ref{thm:dim-infty-CBA}; \ref{SHORT.thm:dim-infty-CBA:bary}$\Rightarrow$\ref{SHORT.thm:dim-infty-CBA:mnfld}$\Rightarrow$\ref{SHORT.thm:dim-infty-CBA:TopDim}.}
The implication \ref{SHORT.thm:dim-infty-CBA:bary}$\Rightarrow$\ref{SHORT.thm:dim-infty-CBA:mnfld} follows from Exercise~\ref{lem:nondeg-test-with-balls}
and Theorem~\ref{thm:bihoelder}, and \ref{SHORT.thm:dim-infty-CBA:mnfld}$\Rightarrow$\ref{SHORT.thm:dim-infty-CBA:TopDim} is trivial.
 
\parit{\ref{SHORT.thm:dim-infty-CBA:TopDim}$\Rightarrow$\ref{SHORT.thm:dim-infty-CBA:bary}.}
According to \ref{thm:stable-value}, 
there is a continuous map $F\:\spc{U}\to \RR^m$ with a stable value.

Fix $\eps>0$.
Since $F$ is continuous, there is a continuous positive function $\rho$ defined on $\spc{U}$ such that 
\[\dist{x}{y}{}<\rho(x)
\quad\Rightarrow\quad
|F(x)- F(y)|<\tfrac13\cdot\eps.\]
Apply \ref{lem:approximation-cba} to $\rho$.
For the resulting simplicial complex $\spc{N}$ 
 and the maps $\map\:\spc{U}\to \spc{N}$, $\Psi\:\spc{N}\to \spc{U}$, we have
\[|F\circ \Psi\circ\map(x)-F(x)|<\tfrac13\cdot\eps\] 
for any $x\in \spc{U}$.

Arguing by contradiction,
assume $\TopDim\spc{U}<m$.
By \ref{lem:approximation-cba:im},
the image of $F\circ\Psi\circ\map$ lies in the $F\circ\Psi$-image of the $(m-1)$-skeleton of $\spc{N}$;
In particular, it can be covered by a countable collection of Lipschitz images of $(m-1)$-simplexes.
Hence
$\bm{0}\in \RR^m$ is not a stable value of $F\circ\Psi\circ\map$.
Since $\eps>0$ is arbitrary, we get the result.
\qeds

The following exercise is a generalization of Helly's theorem; for closely related statements see \cite[Prop. 5.3]{kleiner-1999} and \cite{ivanov2014}.

\begin{thm}{Exercise}\label{ex:helly}
Let $K_1,\ldots,K_n$ be closed convex subsets in a proper length $\CAT(0)$ space $\spc{U}$.
Suppose that $\TopDim \spc{U}\le m$ for some $m<n$ and any $m+1$ subsets from $\{K_1,\ldots, K_n\}$ have a common point.
Show that all subsets $K_1,\ldots,K_n$ have a common point.
\end{thm}

\section{Hausdorff dimension}

\begin{thm}{Definition}
\label{def:HausDim}\index{dimension!Hausdorff dimension}\index{Hausdorff dimension}
Let $\spc{X}$ be a metric space. 
Its \emph{Hausdorff dimension} is defined as
\[\HausDim\spc{X}=\sup\set{\alpha\in\RR}{\HausMes_\alpha(\spc{X})>0},\]
where $\HausMes_\alpha$ denotes the $\alpha$-dimensional Hausdorff measure.
\end{thm}

The following theorem follows from \cite[theorems V 8 and VII 2]{hurewicz-wallman}.

\begin{thm}{Szpilrajn's theorem}\label{thm:szpilrajn} 
Let $\spc{X}$ be a separable metric space.
Assume $\TopDim\spc{X}\ge m$.
Then $\HausMes_m \spc{X}>0$.

In particular, 
$\TopDim\spc{X}\le\HausDim\spc{X}$.
\end{thm}

{\sloppy

Except for Szpilrajn's theorem, there are no other relations between topological and Hausdorff dimension for general separable spaces.
Moreover, the following exercise implies that the same holds for compact geodesic $\CAT(0)$ spaces of topological dimension at least 1.

}

\begin{thm}{Exercise}\label{ex:dim-top-haus-CAT}
Given $\alpha\ge 1$ construct a metric on the binary tree such that it has compact completion of Hausdorff dimension $\alpha$.

Conclude that for any integer $m\ge 1$ and real $\alpha\ge m$ there is a compact $\CAT(0)$ space with topological dimension $m$ and Hausdorff dimension $\alpha$. 
\end{thm}

\section{Notes}

The dimension theorem is the main application of barycenters in Alexandrov geometry.
It is due to Bruce Kleiner \cite{kleiner-1999};
a more precise version was proved by Alexander Lytchak \cite{lytchak:diff}.
Kleiner conjectured that the dimension theorem holds for arbitrary complete geodesic $\CAT(\kappa)$ spaces \cite{kleiner-1999}; see also \cite[p.~133]{gromov-1993}.
For separable spaces, the answer is ``yes'', and it follows from Kleiner's argument \cite[Corollary 14.13]{alexander-kapovitch-petrunin-2025}.

One may wonder if bi-H\"older condition \ref{thm:dim-infty-CBA:mnfld} can be improved to bi-Lipschitz;
this is unknown even for compact spaces.
However if a compact geodesic $\CAT(0)$ space $\spc{U}$ has finite topological dimension $m$,
then a slight modification of Kleiner's technique can be used to show that there is a bi-Lipschitz embedding of an $m$-cube into $\spc{U}$ \cite[Theorem 14.15]{alexander-kapovitch-petrunin-2025}.
In particular, there is a bi-Lipschitz embedding of an $n$-cube for any $n\le m$.
If $\TopDim\spc{U}=\infty$, then we expect existence of a bi-Lipschitz embedding of an $n$-cube for any integer $n\ge1$.
The statement is trivial for $n=1$; in this case any geodesic gives an isometric embedding.
For $n=2$, it follows since minimal and metric minimizing surfaces in $\spc{U}$ are $\CAT(0)$;
see \cite{petrunin-stadler} and the references therein.
Any such surface is locally bi-Lipschitz to the Euclidean plane.
For $n\ge 3$ the question remains open.
% ??? Maybe it is solved by Stadler and Wenger --- we need to check later.

Giuliano Basso studied barycenters in spaces with a conical bicombing~\cite{basso}, which in particular include all geodesic $\CAT(0)$ spaces.

\appendix
%%!TEX root = invitation-CAT.tex
\backmatter
\chapter*{Semisolutions}
\addcontentsline{toc}{chapter}{Semisolutions}
\chaptermark{Semisolutions}

\parbf{\ref{ex:two-components-of-M4}.} 
Let $\spc{X}$ be a 4-point metric space.

Fix a tetrahedron $\triangle$ in~$\RR^3$.
The vertices of $\triangle$, 
say $x_0$, $x_1$, $x_2$, $x_3$, can be identified with the points of~$\spc{X}$.

Note that there is a unique quadratic form $W$ on $\RR^3$
such that 
\[W(x_i-x_j)=\dist[2]{x_i}{x_j}{\spc{X}}\]
for all $i$ and~$j$.

By the triangle inequality, $W(v)\ge 0$ 
for any vector $v$ parallel to one of the faces of~$\triangle$.

Note that $\spc{X}$ is isometric to a 4-point subset in the 3-dimensional Euclidean space
if and only if $W(v)\ge 0$ for any vector $v$ in~$\RR^3$.

\begin{wrapfigure}{r}{33mm}
\vskip-0mm
\centering
\includegraphics{mppics/pic-710}
\end{wrapfigure}

Therefore, if $\spc{X}$ is not of type $\mathcal{E}_4$, then $W(v)<0$ for some vector~$v$.
From above, the vector $v$ must be transversal to each of the 4 faces of~$\triangle$.
Therefore if we project $\triangle$ along $v$ to a plane transversal to $v$ we see one of the two pictures on the right.

Note that the set of vectors $v$ such that $W(v)<0$ has two connected components;
the opposite vectors $v$ and $-v$ lie in different components.
If one moves $v$ continuously, keeping $W(v)<0$,
then the corresponding projection moves continuously and the projections of the 4 faces 
cannot degenerate. 
It follows that the combinatorics of the picture do not depend on the choice of~$v$. 
Hence $\mathcal{M}_4\setminus\mathcal{E}_4$ is not connected. 

It remains to show that if the combinatorics of the pictures for two spaces is the same, then one can continuously deform one space into the other.
This can be easily done by deforming $W$ and applying a permutation of $x_0$, $x_1$, $x_2$, $x_3$ if necessary.

\parit{Comment.}
This solution is inspired by \cite{petrunin-quest}.

\parbf{\ref{ex:convex-set}.}
The simplest proof we know requires the construction of tangent cones.
A closely related statement was proved by Wallace Wilson \cite[§ 12]{wilson}; see also \cite{blumenthal}.

\parbf{\ref{ex:no-geod}.}
Consider the unit ball $(B,\rho_0)$
in the space $c_0$ of all sequences converging to zero equipped with the sup-norm.

Consider another metric $\rho_1$ which is different from $\rho_0$ by the conformal factor
\[\phi(\bm{x})=2+\tfrac{1}2\cdot x_1+\tfrac{1}4\cdot x_2+\tfrac{1}8\cdot x_3+\ldots,\]
where $\bm{x}=(x_1,x_2\,\ldots)\in B$.
That is, if $\bm{x}(t)$, $t\in[0,\ell]$, is a curve parametrized by $\rho_0$-length 
then its $\rho_1$-length is 
\[\length_{\rho_1}\bm{x}=\int\limits_0^\ell\phi\circ\bm{x}.\]
Note that the metric $\rho_1$ is bi-Lipschitz to~$\rho_0$.

Assume $\bm{x}(t)$ and $\bm{x}'(t)$ are two curves parametrized by $\rho_0$-length that differ only in the $m$-th coordinate, denoted by $x_m(t)$ and $x_m'(t)$ respectively.
Note that if $x'_m(t)\le x_m(t)$ for any $t$ and 
the function $x'_m(t)$ is locally $1$-Lipschitz at all $t$ such that $x'_m(t)< x_m(t)$, then 
\[\length_{\rho_1}\bm{x}'\le \length_{\rho_1}\bm{x}.\]
Moreover this inequality is strict if $x'_m(t)< x_m(t)$ for some~$t$.

Fix a curve $\bm{x}(t)$, $t\in[0,\ell]$, parametrized by $\rho_0$-length.
We can choose $m$ large, so that $x_m(t)$ is sufficiently close to $0$ for any~$t$.
In particular, for some values $t$, we have $y_m(t)<x_m(t)$, where
\[y_m(t)=(1-\tfrac t\ell)\cdot x_m(0)
+\tfrac t\ell\cdot x_m(\ell)
-\tfrac 1{100}\cdot \min\{t,\ell-t\}.\]
Consider the curve $\bm{x}'(t)$ as above with
\[x'_m(t)=\min\{x_m(t),y_m(t)\}.\]
Note that $\bm{x}'(t)$ and $\bm{x}(t)$ have the same end points, and by the above
\[\length_{\rho_1}\bm{x}'<\length_{\rho_1}\bm{x}.\]
That is, for any curve $\bm{x}(t)$ in $(B,\rho_1)$, we can find a shorter curve $\bm{x}'(t)$ with the same end points.
In particular, $(B,\rho_1)$ has no geodesics.

\parit{Comment.}
This example is due to Fedor Nazarov \cite{nazarov}.

\begin{wrapfigure}{r}{20 mm}
\vskip-0mm
\centering
\includegraphics{mppics/pic-1}
\end{wrapfigure}

\parbf{\ref{exercise from BH}.}
Consider the following subset of $\RR^2$ equipped with the induced length metric
\[
\spc{X}
=
\bigl((0,1]\times\{0,1\}\bigr)
\cup
\bigl(\{1,\tfrac12,\tfrac13,\ldots\}\times[0,1]\bigr)
\]
Note that $\spc{X}$ is locally compact and geodesic.

Its completion $\bar{\spc{X}}$ is isometric to the closure of $\spc{X}$ equipped with the induced length metric;
$\bar{\spc{X}}$ is obtained from $\spc{X}$ by adding two points $p=(0,0)$ and $q=(0,1)$.

The point $p$ admits no compact neighborhood in $\bar{\spc{X}}$ 
and there is no geodesic connecting $p$ to $q$ in~$\bar{\spc{X}}$. 

\parit{Comment.}
This example is taken from the book by Martin Bridson and André Haefliger \cite[I.3.6(4)]{bridson-haefliger}.

\parbf{\ref{ex:length-prod}.}
Let $\spc{W}=\spc{U}\times \spc{V}$.
Choose two pairs of points $u_0,u_1\in \spc{U}$
and $v_0,v_1\in \spc{V}$.
Set $a=\dist{u_0}{u_1}{\spc{U}}$, $b=\dist{v_0}{v_1}{\spc{V}}$ and $c=\dist{(u_0,v_0)}{(u_1,v_1)}{\spc{W}}$.

Since $\spc{U}$ and $\spc{V}$ are length spaces,
given $\eps>0$,
we can choose curves $\alpha\:[0,1]\to \spc{U}$ from $u_0$ to $u_1$
and $\beta\:[0,1]\to \spc{V}$ from $v_0$ to $v_1$
such that 
\[\length \alpha<a+\eps
\quad\text{and}\quad
\length \beta<b+\eps.\]
Reparametrizing the paths proportional to their lengths we can assume that
$\alpha$ is $(a+\eps)$-Lipschitz
and
$\beta$ is $(b+\eps)$-Lipschitz.
Therefore the path $t\mapsto (\alpha(t),\beta(t))$ is $(c+2\cdot \eps)$-Lipschitz.
In particular, its length cannot exceed $c+2\cdot \eps$ for any $\eps>0$.
Hence $\spc{W}$ meets the conditions in the definition of a length space.

\parbf{\ref{ex:geod-prod}.}
Let $\gamma\:t\mapsto (\alpha(t),\beta(t))$ be a geodesic path in $\spc{W}=\spc{U}\times \spc{V}$.
Show and use that 
\begin{align*}
\dist{\alpha(t_0)}{\alpha(t_1)}{\spc{U}}
&=
a\cdot \dist{\gamma(t_0)}{\gamma(t_1)}{\spc{W}}
\\
\dist{\beta(t_0)}{\beta(t_1)}{\spc{V}}
&=
b\cdot \dist{\gamma(t_0)}{\gamma(t_1)}{\spc{W}}
\end{align*}
for any $t_0$ and $t_1$ and some fixed values $a\ge 0$ and $b\ge 0$ such that $a^2+b^2=1$.
To see this, first show that both $\alpha$ and $\beta$ are distance minimizing. Use this to reduce the problem to the case where both $\spc{U}$ and $\spc{V}$ are intervals.

\parbf{\ref{ex:cone-geod}.}
Let $\gamma$ be a unit-speed parametrization of $[pq]$.
Show that after shifting the parametrization, we can assume that $|\gamma(t)|=\sqrt{a^2+t^2}$ for some constant $a$.

Let $\hat \gamma(t)$ be the projection of $\gamma(t)$ to $\spc{U}$.
Show and use that $t\z\mapsto \hat\gamma(a\cdot \tan t)$ is a geodesic in $\spc{U}$.
To see this, first observe that $\hat \gamma$ is distance minimizing. Use this to reduce the problem to the case when $\spc{U}$ is an interval.

\parbf{\ref{ex:product-cone}.} A point in $\RR\times \Cone \spc{U}$ can be described by a triple $(x,r,p)$, where $x\in \RR$, $r\in \RR_{\ge}$ and $p\in \spc{U}$.
Correspondingly, a point in $\Cone[\Susp\spc{U}]$ can be described by a triple $(\rho,\phi,p)$, where $\rho\in \RR_\ge$, $\phi\in [0,\pi]$ and $p\in \spc{U}$.

The map 
$\Cone[\Susp\spc{U}]\to\RR\times\Cone\spc{U}$ defined as
\[(\rho,\phi,p)\mapsto(\rho\cdot\cos\phi,\rho\cdot\sin\phi,p)\] 
is the needed isometry.

\parbf{\ref{ex:angle-on-shortest}.} Apply the triangle inequality and spell out the definitions.

\parbf{\ref{ex:adjacent-angles}.}
Use the triangle inequality for angles and \ref{ex:angle-on-shortest}.
An example with strict inequality can be obtained by gluing together two copies of a solid planar triangle with an obtuse angle along a side adjacent to the obtuse angle.

\parbf{\ref{ex:hausdorff-conv}.}
By definition of Hausdorff convergence
\[p\in A_\infty\quad\iff\quad\distfun_{A_n}(p)\to 0\quad\text{as}\quad n\to \infty.\] 
The latter is equivalent to the existence of a sequence $p_n\in A_n$ such that
$\dist{p_n}{p}{}\to0$ as $n\to \infty$;
or equivalently $p_n\to p$.
Hence the first statement follows.

The converse is false.
For example, consider the alternating sequence of two distinct closed sets $A,B,A,B,\ldots$;
note that it is not a convergent sequence in the sense of Hausdorff.
On the other hand, the set of all limit points is well-defined --- it is the intersection $A\cap B$.

\parit{Comment.} The set $\ushort{A}_\infty$ of all limits of sequences $p_n\in A_n$ is called the \emph{lower closed limit}
and the set $\bar{A}_\infty$ of all partial (also known as subsequential) limits of such sequences is called the \emph{upper closed limit}.
Clearly $\ushort{A}_\infty\subset \bar{A}_\infty$.
If $\ushort{A}_\infty=\bar{A}_\infty$, then it is called the \emph{closed limit} of~$A_n$.
All these convergences were introduced by Felix Hausdorff in~\cite{hausdorff}.

For the class of closed subsets of a proper metric space, closed limits coincide with limits in the sense of Hausdorff as we defined them.

\parbf{\ref{ex:non-contracting-map}.}
Given any pair of points $x_0,y_0\in \spc{K}$, 
consider two sequences $(x_n)$ and $(y_n)$
such that 
$x_{n+1}=f(x_n)$ and $y_{n+1}=f(y_n)$ for each~$n$.

Since $\spc{K}$ is compact, 
we can choose an increasing sequence of integers $n_i$
such that both sequences $(x_{n_i})_{i=1}^\infty$ and $(y_{n_i})_{i=1}^\infty$
converge.
In particular, both of these sequences are Cauchy;
that is,
\[
|x_{n_i}-x_{n_j}|_{\spc{K}}, |y_{n_i}-y_{n_j}|_{\spc{K}}\to 0
\ \ 
\text{as}
\ \ \min\{i,j\}\to\infty.
\]

Since $f$ is distance non-decreasing, we get
\[
|x_0-x_{|n_i-n_j|}|
\le 
|x_{n_i}-x_{n_j}|.
\]

It follows that
there is a sequence $m_i\to\infty$ such that
\[
x_{m_i}\to x_0\ \ \text{and}\ \ y_{m_i}\to y_0\ \ \text{as}\ \ i\to\infty.
\leqno({*})\]

Set \[\ell_n=|x_n-y_n|_{\spc{K}}.\]
Since $f$ is distance non-decreasing, $(\ell_n)$ is a non-decreasing sequence.

By $({*})$, $\ell_{m_i}\to\ell_0$ as $m_i\to\infty$.
It follows that $(\ell_n)$ is a constant sequence.

In particular 
\[|x_0-y_0|_{\spc{K}}=\ell_0=\ell_1=|f(x_0)-f(y_0)|_{\spc{K}}\]
for any pair $x_0$ and~$y_0$.
That is, $f$ is distance-preserving, in particular, injective.

From $({*})$, we also get that $f(\spc{K})$ is dense in $\spc{K}$.
Since $\spc{K}$ is compact, $f\:\spc{K}\to \spc{K}$ is surjective.
Hence the result.

\parit{Comment.}
This exercise is a basic introductory lemma on Gromov--Hausdorff distance (see for example \cite[7.3.30]{burago-burago-ivanov}).
The presented proof is not quite standard, it was found by Travis Morrison.

\parbf{\ref{ex:compact-proper-GH}.}
To prove part \ref{SHORT.ex:compact-proper-GH:a},
fix a countable dense set of points $\mathfrak{S}\subset\spc{X}_\infty$.
For each point $x\in \mathfrak{S}$, choose a sequence 
of points $x_n\in\spc{X}_n$ such that $x_n\xrightarrow{\rho} x$.

Applying the diagonal procedure, we can pass to a subsequence of $ \spc{X}_n$ such that each of the constructed sequences $\rho'$-converge;
that is, $x_n\xrightarrow{\rho'} x'$ for some $x'\in\spc{X}_\infty'$.

In this way we get a map $\mathfrak{S}\to\spc{X}_\infty'$ defined as $x\mapsto x'$.
Note that this map preserves distances and therefore can be extended to a distance-preserving map $\spc{X}_\infty\to\spc{X}_\infty'$.
Likewise we construct a distance-preserving map $\spc{X}_\infty'\to\spc{X}_\infty$.

It remains to apply Exercise~\ref{ex:non-contracting-map}.

The proof of part \ref{SHORT.ex:compact-proper-GH:b} is nearly identical,
but one has to apply Exercise~\ref{ex:non-contracting-map} to closed balls centered at the limits of $x_n$ in $\spc{X}_\infty$ and $\spc{X}_\infty'$. 

\parbf{\ref{ex:noncreasing-CAT}.}
Note that it is sufficient to show that 
$\angk p{\bar x}{y}\le \angk p{x}{y}$ for any $\bar x\in \left]px\right[$.
The latter follows from Alexandrov's lemma (\ref{lem:alex}) and the $\CAT(0)$ comparison for the quadruple $p, x, \bar x, y$.

The last statement follows from a monotone convergence argument.

\parbf{\ref{ex:contractible}.}
Observe that 
$\dist{h_t(x)}{h_t(y)}{}\le t\cdot \dist{x}{y}{}$
and therefore
\[\dist{h_{t_0}(x)}{h_{t_1}(y)}{}\le t_0\cdot \dist{x}{y}{}+|t_0-t_1|\cdot \dist{p}{y}{}.\]
The conclusion follows.

\parbf{\ref{ex:CAT-mnfld=>ext.geod}.}
Assume that a geodesic $[px]$ cannot be extended behind $x$.
Apply the homotopy from \ref{ex:contractible} to prove that $\spc{U}$ has vanishing local homology (and/or homotopy) groups at $x$.
Use that manifolds have some nontrivial local homologies (and/or homotopy) groups. 

\parbf{\ref{ex:geod-CBA}.}
Fix a sufficiently small $\eps>0$.
Recall that by Proposition~\ref{cor:loc-geod-are-min}, any local geodesic in $\spc{U}$ is a geodesic.

Consider a sequence of directions $\xi_n$ 
at $p$
of geodesics $[pq_n]$.
We can assume that the distances $\dist{p}{q_n}{\spc{U}}$ are equal to $\eps$ for all $n$;
here we use that the geodesics are extendable.

Since $\spc{U}$ is proper, we can pass to a converging subsequence of $q_n$;
suppose $q$ is its limit.
Show that the direction $\xi$ of $[pq]$ is the limit of directions $\xi_n$.

\begin{wrapfigure}{r}{45 mm}
\vskip-4mm
\centering
\includegraphics{mppics/pic-12}
\end{wrapfigure}

The unit disc in the plane with
a half-line attached to each point of its boundary is a complete $\CAT(0)$ length space (prove this!) with extendable geodesics.
However, the space of geodesic directions on the boundary of the disc is not complete --- there is no geodesic tangent to the boundary of the disc.
This provides a counterexample to the statement of the exercise if $\spc{U}$ is not assumed to be proper.

\parit{Comment.}
Properties of proper (locally) $\CAT(\kappa)$ spaces with extendable geodesics were studied by Alexander Lytchak and Koichi Nagano \cite{lytchak-nagano-2019}.

\parbf{\ref{ex:tan-CBA}}; \ref{SHORT.ex:tan-CBA:geodesic}.
Given two points $p$ and $q$ choose a sequence of $\eps$-midpoints for $\eps\to 0$.
Use the comparison to prove that this sequence is Cauchy.
Observe that its limit is a midpoint of $p$ and $q$ and apply Menger's lemma (\ref{lem:mid>geod}).

\parit{\ref{SHORT.ex:tan-CBA:tan}.}
By \ref{SHORT.ex:tan-CBA:geodesic}, it is sufficient to show that $\T_p$ is a length $\CAT(0)$ space.

Choose four geodesic vectors $v_1$, $v_2$, $v_3$, and $v_4$ in $\T_p$.
Let $\gamma_i$ be a geodesic that starts at $p$ in directions $v_i$.
Reparametrize each $\gamma_i$ with speed $|v_i|$ so that $\gamma_i(0)=p$.
Given small $t>0$, consider the quadruple of points $a_i=\gamma_i(t), i=1,2,3,4$.
Observe that $\CAT(0)$ comparison holds for the quadruple $(a_1$, $a_2$, $a_3$, $a_4)$ in $\spc{U}$ for all small $t>0$, and it implies $\CAT(0)$ comparison holds for $(v_1$, $v_2$, $v_3$, $v_4)$.
Since geodesic vectors are dense in $\T_p$, we get $\CAT(0)$ comparison in $\T_p$.

It remains to show that $\T_p$ is a length space.
By construction, $\T_p$ is complete.
By Menger's lemma (\ref{lem:mid>geod}) it is sufficient to construct
an $\eps$-midpoint between the given vectors $v_1,v_2\in\T_p$ and $\eps>0$.
Since geodesic vectors are dense in $\T_p$, we can assume that $v_1$, and $v_2$ are geodesic vectors.

Choose $t>0$ and define $a_i=\gamma_i(t)$ as above.
Let $q$ be the midpoint of $[a_1\,a_2]$.
Denote by $w$ the tangent vector of the geodesic from $p$ to $q$ with a constant-speed parametrization by $[0,t]$.
Observe that $w$ is an $\eps$-midpoint of $v_1$ and $v_2$, where $\eps\to0$ as $t\to 0$.

\parbf{\ref{ex:convex-distfun}}; \textit{Only-if part.}
It is sufficient to prove Jensen's inequality
\[h(\tfrac{t_0+t_1}2)
\le
\frac{h(t_0)+h(t_1)}2\]
where $h(t)\df f\circ \gamma(t)-\tfrac12\cdot t^2$.
Observe that the inequality holds in the Euclidean plane, and 
apply the thinness of the triangle $[p\,\gamma(t_0)\,\gamma(t_1)]_{\spc{U}}$.

\parit{If part.} Reverse the above argument.

\parbf{\ref{ex:convex-dist}.}
Observe that it is sufficient to show that 
\[\dist{\gamma_1(t)}{\gamma_2(t)}{\spc{U}}\le (1-t)\cdot\dist{\gamma_1(0)}{\gamma_2(0)}{\spc{U}}+t\cdot \dist{\gamma_1(1)}{\gamma_2(1)}{\spc{U}}.\]

\begin{figure}[htb!]
\vskip-0mm
\centering
\includegraphics{mppics/pic-1650}
\vskip0mm
\end{figure}

Let $\beta$ be the geodesic path from $\gamma_1(0)$ to $\gamma_2(1)$.
Observe that 
\begin{align*}
\dist{\gamma_1(t)}{\beta(t)}{\spc{U}}&\le t\cdot \dist{\gamma_1(1)}{\beta(1)}{\spc{U}},
\\
\dist{\beta(t)}{\gamma_2(t)}{\spc{U}}&\le (1-t)\cdot\dist{\beta(0)}{\gamma_2(0)}{\spc{U}},
\end{align*}
and apply the triangle inequality.

\parbf{\ref{ex:displacement}.}
Choose a geodesic $\gamma_1$.
Observe that $\gamma_2=\iota\circ\gamma_1$ is a geodesic as well and apply \ref{ex:convex-dist}.

\parbf{\ref{ex:convex-nbhd}.}
It is sufficient to show that 
\[\distfun_A\circ\gamma(t)\le (1-t)\cdot\distfun_A\circ\gamma(0)+t\cdot \distfun_A\circ \gamma(1).
\leqno({*})\]
for any geodesic path $\gamma$.
Note that given $\eps>0$, there are points $p,q\in A$ such that 
\[\dist{p}{\gamma(0)}{}<\distfun_A\circ\gamma(0)+\eps
\quad\text{and}\quad
\dist{q}{\gamma(1)}{}< \distfun_A\circ\gamma(1)+\eps.
\]
Let $\beta$ be a geodesic path from $p$ to $q$.
By \ref{ex:convex-dist},
\[\dist{\beta(t)}{\gamma(t)}{}\le (1-t)\cdot\dist{p}{\gamma(0)}{}+t\cdot \dist{q}{\gamma(1)}{}.\]
Since $A$ is convex, $\beta(t)\in A$ for any $t$.
Since $\eps>0$ is arbitrary, we get~$({*})$.

\parbf{\ref{ex:closest-point}.}
Since $\spc{U}$ is proper, the set $K\cap \cBall[p,R]$ is compact for any $R<\infty$.
The existence of at least one point $p^*$ that minimizes the distance from $p$ follows.

Assume $p^*$ is not uniquely defined;
that is, two distinct points in $K$, say $x$ and $y$, minimize the distance from~$p$.
Since $K$ is convex, the midpoint $z$ of $[xy]$ lies in~$K$.

Thinness of triangles implies that
\[\dist{p}{z}{}<\dist{p}{x}{}=\dist{p}{y}{}\]
--- a contradiction.

\begin{wrapfigure}{r}{37 mm}
\vskip-0mm
\centering
\includegraphics{mppics/pic-41}
\vskip-0mm
\end{wrapfigure}

It remains to show that the map $p\mapsto p^*$ is short, 
that is, 
\[\dist{p}{q}{}\ge \dist{p^*}{q^*}{}
\leqno({*})\]
for any $p,q\in \spc{U}$.

Assume $p\ne p^*$, $q\ne q^*$, $p^*\ne q^*$.
Construct the model triangles 
$\trig{\tilde p}{\tilde p^*}{\tilde q^*}$ and $\trig{\tilde p}{\tilde q}{\tilde q^*}$
of $\trig{p}{p^*}{q^*}$ and $\trig{p}{q}{q^*}$ so that 
the points $\tilde p^*$ and $\tilde q$ lie on the opposite sides from $[{\tilde p}{\tilde q^*}]$.

From thinness of the triangles $\trig{p}{p^*}{q^*}$ and $\trig{p}{q}{q^*}$, we get that
\[\mangle\hinge{\tilde p^*}{\tilde p}{\tilde q^*},\mangle\hinge{\tilde q^*}{\tilde p^*}{\tilde q}\ge \tfrac\pi2.\]
Hence 
\[\dist{\tilde p}{\tilde q}{}\ge \dist{\tilde p^*}{\tilde q^*}{}.\]
The latter is equivalent to $({*})$.

In the remaining cases: $({*})$ holds automatically if (1) $p^*= q^*$ or (2) $p= p^*$ and $q= q^*$.
If $p= p^*$, $q\ne q^*$ and $p^*\ne q^*$, then thinness of $\trig{p}{q}{q^*}$ implies that
\[
\mangle\hinge{\tilde q^*}{\tilde q}{\tilde p}\ge \tfrac\pi2,\]
and $({*})$ follows.

\parit{Comment.}
It is sufficient to assume that the space is complete length and $\CAT(0)$; see \cite{alexander-kapovitch-petrunin-2025}.

\parbf{\ref{ex:locally-convex}.}
Fix a closed, connected, locally convex set~$K$.
Apply \ref{ex:convex-nbhd} to show that $\distfun_{K}$ is convex in a neighborhood $\Omega\supset K$; that is, $\distfun_{K}$ is convex along any geodesic completely contained in~$\Omega$.

Since $K$ is locally convex,
it is locally path connected.
Since $K$ is connected, it is also path connected.

Fix two points $x,y\in K$. 
Let us connect $x$ to $y$ by a path $\alpha\:[0,1]\z\to K$.
Use \ref{prop:thin=cat} to show that the geodesics $[x\,\alpha(s)]$
are uniquely defined and depend continuously on~$s$.

\begin{wrapfigure}{r}{50 mm}
\vskip-0mm
\centering
\includegraphics{mppics/pic-9}
\end{wrapfigure}

If $[xy]=[x\,\alpha(1)]$ does not completely lie in $K$, then 
there is a value $s\in [0,1]$ such that $[x\,\alpha(s)]$ 
lies in $\Omega$,
but does not completely lie in~$K$.
Therefore $f=\distfun_{K}$ is convex 
along $[x\alpha(s)]$.

Note that $f(x)\z=f(\alpha(s))=0$ and $f\ge 0$, 
therefore $f(z)= 0$ for any $z\z\in [x\,\alpha(s)]$;
that is, $[x\,\alpha(s)]\z\subset K$ --- a contradiction.

\parit{Comment.}
The proof presented here is nearly identical to the original proof given by Heinrich Tietze of a less general statement~\cite{tietze}. 

\parbf{\ref{ex:reshetnyak-doubling}.}
The ``if'' part follows from the Reshetnyak gluing theorem (\ref{thm:gluing}).

Assume $\spc{W}$ is $\CAT(0)$.
Note that one copy of $\spc{U}$ embeds isometrically in $\spc{W}$.
Conclude that $\spc{U}$ is $\CAT(0)$.

Assume $A$ is not convex;
that is, $[xy]\not\subset A$ for some $x,y\in A$.
Observe that there are distinct geodesics from $x$ to $y$ in $\spc{W}$.
Arrive at a contradiction with the uniqueness of geodesics (\ref{ex:CAT-geodesic}).

\parbf{\ref{ex:supporting-planes}.}
Note that $f$ is convex as the sum of two convex functions.
By approximation, it is sufficient to consider the case when 
$A$ and $B$ have smooth boundary. 

If $[xy]\cap A\cap B\ne \emptyset$, then $z\in [xy]$ and $\dot A, \dot B$ can be chosen to be arbitrary half-spaces containing $A$ and $B$ respectively.

In the remaining case $[xy]\cap A\cap B=\emptyset$, 
we have $z\in\partial (A\cap B)$.
Consider the solid ellipsoid
\[C=\set{p\in\EE^m}{f(p)\le f(z)}.\]
Note that $C$ is compact, convex and has smooth boundary.

Suppose $z\in\partial A \cap \Int B$.
Then $A$ and $C$ touch at $z$ and we can set $\dot A$ to be the uniquely defined supporting half-space to $A$ at $z$ and $\dot B$ to be any half-space containing~$B$.
The case $z\in\partial B \cap \Int A$ is treated similarly.

Finally, suppose $z\in\partial A\cap\partial B$.
Then the set $\dot A$ (respectively, $\dot B$) is defined as the unique supporting half-space to $A$ (respectively, $B$) at $z$ containing $A$ (respectively, $B$).

Suppose $f(p)<f(z)$ for some $p\in \dot A\cap\dot B$.
Since $f$ is convex,
$f(\bar z)<f(z)$ for any $\bar z \in \left[pz\right[$.
Since $\left[pz\right[\cap A\cap B\ne\emptyset$,
the latter contradicts the fact that $z$ is a minimum point of $f$ on $A\cap B$.

\begin{wrapfigure}{r}{40 mm}
\vskip-0mm
\centering
\includegraphics{mppics/pic-43}
\vskip-0mm
\end{wrapfigure}

\parbf{\ref{ex:compact-walls}.}
Fix two open balls $B_1\z=\oBall(0,r_1)$ and $B_2=\oBall(0,r_2)$
such that 
\[B_1\subset A_i\subset B_2\]
for each wall~$A_i$.

Note that all the intersections of the walls have $\eps$-wide corners for
\[\eps=2\cdot \arcsin\tfrac{r_1}{r_2}.\]
The proof can be guessed from the picture.

\parbf{\ref{ex:centrally-simmetric-walls}.}
Observe that any closed unbounded centrally symmetric set in $\RR^m$ splits off an $\RR$ factor.
Therefore, any centrally symmetric closed convex set in $\RR^m$ is a product of a compact centrally symmetric convex set and a linear subspace of $\RR^m$.

It follows that there is $R<\infty$
such that if $X$ is an intersection of an arbitrary number of walls, then for any point $p\in X$ there is an isometry of $X$ 
that moves $p$ to a point in the ball $\oBall(0,R)$.

It remains to repeat the proof of Exercise~\ref{ex:compact-walls}.

\parbf{\ref{cor:balls:dim=1}.}
Imagine that each ball has zero radius; then we may think that the balls pass thru each other.
That is, every ball moves with constant speed along the line.
Let $x_i(t)$ be the coordinate of $i$-th ball at time~$t$.
Note that the graph $x_i$ in the $(t,x)$-plane is a straight line.
Every two lines have at most one intersection and each collision corresponds to one such intersection.
We have $n$ lines and therefore at most $\binom n2$ collisions.

The general case, with balls of radius $r>0$ can be reduced to the case above.
To do this, one has to exclude the space occupied by the balls.
In other words, if the $i$-th ball is centered at $x$, then we assume that its coordinate is
$x-2\cdot r\cdot i$.
With these new coordinates, the balls behave exactly as the balls with vanishing radii.

\parit{Comment.}
Generalization of this problem for different masses is discussed in some popular books \cite{galperin-zemlyakov,tabachnikov}.
Let us also recommend a video by Grant Sanderson \cite{sanderson} about a peculiar connection between $\pi$ and the number of collisions in 1-dimensional billiards.

\parbf{\ref{ex:collision}.}
Assume that a trajectory $\gamma$ is defined on the interval $[a,b)$
and the collisions accumulate at $b$.
Without loss of generality, we may assume that $\gamma$ meets each wall infinitely many times.
Consider the infinite puff pastry $\spc{R}^\gamma$ for $\gamma$
and let $\bar\gamma$ be its lift.

Pass to the completion $\bar{\spc{R}}^\gamma$ of $\spc{R}^\gamma$;
observe that $\bar{\spc{R}}^\gamma$ is $\CAT(0)$.
Define $\bar\gamma(b)$ to be the limit point of $\bar\gamma(t)$ as $t\to b$.
Notice that $\bar\gamma\:[a,b]\to \bar{\spc{R}}_\gamma$ is a minimizing geodesic.

Show that $\bar\gamma(b)$ lies in the lift of the intersection of all walls;
in other words, $\bar\gamma(b)$ belongs to the intersection of all the levels of the puff pastry.
Conclude that $\bar\gamma$ completely lies on the lowest level, and arrive at a contradiction. 

Finally, by taking two tangent discs as the walls of a billiard table, we obtain the needed example.
Indeed, a trajectory that starts near the common point of the discs in the direction perpendicular to their common tangent line will have to bounce intensively for quite a while.

\parbf{\ref{ex:convex-dist:major}.}
Observe that the statement holds in the Euclidean plane and apply the majorization theorem to the quadrangle $[\gamma_1(0)\gamma_1(1)\gamma_2(1)\gamma_2(0)]$.

\parbf{\ref{ex:devel-comp-CAT}.} Apply \ref{ex:convex-distfun}.

\parbf{\ref{ex:fenchel}.}
Suppose $D$ is a convex figure that majorizes $[p_1\ldots p_n]$.
Show that $D$ is an $n$-gon and, for each $i$, the majorization sends
a vertex of $D$, say $\tilde p_i$, to $p_i$.
Show that the external angle of $D$ at $\tilde p_i$ cannot exceed the external angle at $p_i$;
make the conclusion.

\begin{wrapfigure}{r}{25 mm}
\vskip-0mm
\centering
\includegraphics{mppics/pic-416}
\end{wrapfigure}

\parbf{\ref{ex:FM}.}
The claimed lower bound on the sum of external angles is a theorem of the first author and Richard Bishop \cite{alexander-bishop-1998}.
The required polygon is shown on the diagram;
it lies in the product space of the real line and a \index{tripod}\emph{tripod}; that is, three line segments glued together at one end.

There are no such examples in the Euclidean space.
Indeed, the original Fáry--Milnor theorem implies that 
if $\beta$ is a nontrivial knot in Euclidean space, then the
sum of its external angles is \textit{strictly} larger than $4\cdot\pi$.

\parbf{\ref{ex:arm-lemma}.}
Note that the majorization theorem for polygons is already proved, and we can use it.

By the majorization theorem, $P=[x_1\ldots x_{n}]$ can be majorized by a convex plane polygon, say $\dot P=[\dot x_1\ldots \dot x_{n}]$.
Note that we have
\[
\dist{\dot x_i}{\dot x_{i-1}}{\EE^2}
=
\dist{x_i}{x_{i-1}}{\spc{U}}
\quad \text{and}\quad
\mangle\hinge{\dot x_i}{\dot x_{i-1}}{\dot x_{i+1}}_{\EE^2}
\ge
\mangle\hinge{x_i}{x_{i-1}}{x_{i+1}}_{\spc{U}}
\]
for all $i$, and $\dist{\dot x_1}{\dot x_n}{\EE^2}
=
\dist{x_1}{x_{n}}{\spc{U}}$ as well.%
\footnote{We also get that
\[\mangle\hinge{\dot x_1}{\dot x_2}{\dot x_n}
\ge
\mangle\hinge{x_1}{x_2}{x_n}
\quad\text{and}\quad
\mangle\hinge{\dot x_{n}}{\dot x_1}{\dot x_{n-1}}
\ge
\mangle\hinge{x_n}{x_1}{x_{n-1}},\] but we will not need it.}

It follows that
\[
\dist{\dot x_i}{\dot x_{i-1}}{\EE^2}
=
\dist{\tilde x_i}{\tilde x_{i-1}}{\EE^2}
\quad \text{and}\quad
\mangle\hinge{\dot x_i}{\dot x_{i-1}}{\dot x_{i+1}}
\ge
\mangle\hinge{\tilde x_i}{\tilde x_{i-1}}{\tilde x_{i+1}}
\]
for all $i$.
By the classical arm lemma,
\[\dist{\dot x_1}{\dot x_{n}}{\EE^2}
\ge
\dist{\tilde x_1}{\tilde x_{n}}{\EE^2}\]
--- hence the result.

\parit{Remarks.}
In the equality case, the map $\tilde x_i\mapsto x_i$ can be extended to a distance-preserving map from the convex hull of $\tilde P$ to $\spc{U}$; see \cite[9.63]{alexander-kapovitch-petrunin-2025}.

The classical arm lemma has the following close relative in differential geometry.

\parbf{Bow lemma.}
\textit{Suppose $\tilde\gamma$ is a smooth arc of a convex planar curve, that is, it runs along the boundary of some convex domain in $\EE^2$.
Further, let $\gamma$ be a smooth curve in $\EE^m$.
Assume that both curves come with arc-length parametrization by the interval $[a,b]$ and at each time the curvature of $\gamma$ does not exceed the curvature of $\tilde\gamma$.
Then the distance between the endpoints of $\gamma$ is at least the distance between the endpoints of $\tilde\gamma$; that is,
\[\dist{\gamma(b)}{\gamma(a)}{\EE^m} \ge \dist{\tilde\gamma(b)}{\tilde\gamma(a)}{\EE^2}.\]
Moreover, in case of equality, the map $\tilde\gamma(t)\mapsto \gamma(t)$ can be extended to a distance-preserving map $\EE^2\to\EE^m$.}

\medskip

This lemma can be formulated for $\CAT(\kappa)$ spaces using the notion of extrinsic curvature introduced by the first author and Richard Bishop \cite{alexander-bishop-2006,alexander-bishop-2010}.
A special case of this $\CAT(\kappa)$ version was proved earlier by the same authors \cite{alexander-bishop-1998};
as far as we know, the proof in  the general case has not been written anywhere.

\parbf{\ref{ex:null-homotopic}.}
Note that the existence of a null-homotopy is equivalent to the following.
There are two continuous one-parameter families of paths $\alpha_\tau$ and $\beta_\tau$, $\tau\in[0,1]$
such that: 
\begin{itemize}
\item $\length\alpha_\tau$, $\length\beta_\tau<\pi$ for any~$\tau$.
\item $\alpha_\tau(0)=\beta_\tau(0)$ and $\alpha_\tau(1)=\beta_\tau(1)$ for any~$\tau$.
\item $\alpha_0(t)=\beta_0(t)$ for any~$t$.
\item $\alpha_1(t)=\alpha(t)$ and $\beta_1(t)=\beta(t)$ for any~$t$.
\end{itemize}

By Corollary~\ref{cor:discrete-paths},
$\gamma_{\alpha_\tau}=\gamma_{\beta_\tau}$ for any $\tau$ --- hence the result.

\parbf{\ref{ex:geod-circle}.}
By the globalization theorem there is a nontrivial homotopy class of closed curves.

Consider a shortest noncontractible closed curve $\gamma$ in $\spc{X}$;
due to compactness of $\spc{X}$ such a curve exists.
Indeed, let $\ell$ be the infimum of lengths of all noncontractible closed curves in $\spc{X}$.
Geodesic homotopy construction implies that two sufficiently close closed curves in $\spc{X}$ are homotopic.
Then choosing a sequence of unit speed noncontractible curves whose lengths converge to $\ell$, an Arzel\'{a}--Ascoli-type argument shows that these curves subconverge to a noncontractible curve $\gamma$ of length~$\ell$.

Assume that $\gamma$ is not a geodesic circle,
that is, there are two points $p$ and $q$ on $\gamma$ such that the distance $\dist{p}{q}{}$
is shorter than the lengths of the arcs, say $\alpha_1$ and $\alpha_2$, of $\gamma$ from $p$ to~$q$.
Consider the product curves, say $\gamma_1$ and $\gamma_2$,
of $[qp]$ with $\alpha_1$ and~$\alpha_2$.
Then $\gamma_1$ or $\gamma_2$ is noncontractible,
\[\length\gamma_1<\ell
\quad\text{and}\quad
\length\gamma_2<\ell\]
--- a contradiction.

{

\begin{wrapfigure}{r}{25 mm}
\vskip-10mm
\centering
\includegraphics{mppics/pic-11}
\end{wrapfigure}

In the $\CAT(1)$ case we also have a geodesic circle.
The proof is done in nearly  the same way, but we need to consider the homotopy classes of closed curves shorter than $2\cdot \pi$.
One also needs to apply \ref{ex:null-homotopic}, to show that curves $\gamma_1$ and $\gamma_2$ are not contractible in the class of curves shorter than $2\cdot \pi$.

}

\parit{Comment.}
Note that the surface of revolution of the graph of $y\z=e^x$ around the $x$-axis is locally $\CAT(0)$ but has no closed geodesics.
Therefore, in this exercise, one cannot trade compactness of $\spc{X}$ for properness.

\parbf{\ref{ex:branching-cover}.}
Consider a closed $\eps$-neighborhood $A_\eps$ of the geodesic $\gamma$.
Note that $A_\eps$ is convex.
By the Reshetnyak gluing theorem, the double $\spc{W}_\eps$ of $\spc{U}$ along $A_\eps$ is $\CAT(0)$.

Consider the space $\spc{W}'_\eps$ obtained by
taking the double cover $\spc{U}\setminus A_\eps$ and gluing back~$A_\eps$.

Observe that $\spc{W}'_\eps$ is locally isometric to $\spc{W}_\eps$. 
That is, for any point $p'\in\spc{W}'_\eps$ there is a point $p\in\spc{W}_\eps$ such that the $\delta$-neighborhood of $p'$ is isometric to the $\delta$-neighborhood of $p$ for all small $\delta>0$.

Further observe that $\spc{W}'_\eps$ is simply connected since it admits a deformation retraction onto $A_\eps$ (following the nearest-point projection), which is contractible.
By the globalization theorem, $\spc{W}'_\eps$ is $\CAT(0)$.

It remains to note that $\tilde {\spc{U}}$ can be obtained as the limit of $\spc{W}'_\eps$ as $\eps\to 0$, and apply~\ref{prop:cat-limit}.

\parbf{\ref{ex:unique-geod=CAT}.}
By the globalization theorem (\ref{thm:hadamard-cartan}),
it is sufficient to show that $\spc{P}$ is locally $\CAT(0)$.
Assume the contrary.
Then by \ref{thm:PL-CAT}, a link $\Sigma$ of some simplex contains a closed geodesic $\alpha$ with length $4\cdot\ell<2\cdot\pi$.
We can assume that $\Sigma$ has minimal possible dimension;
so, by \ref{thm:PL-CAT}, $\Sigma$ is locally $\CAT(1)$.

Divide $\alpha$ into two equal arcs $\alpha_1$ and $\alpha_2$.

Assume $\alpha_1$ and $\alpha_2$ are length minimizing;
parameterize them by $[-\ell,\ell]$.
Fix a small $\delta>0$ and 
consider two curves in $\Cone\Sigma$ written in polar coordinates as 
\[\gamma_i(t)=(\alpha_i(\arctan \tfrac t\delta),\sqrt{\delta^2+t^2}).\]
Observe that both curves $\gamma_1$ and $\gamma_2$ are geodesics in $\Cone\Sigma$ and have common ends.

Note that a small neighborhood of the tip of $\Cone\Sigma$ admits an isometric embedding into~$\spc{P}$.
Hence we can construct two geodesics $\gamma_1$ and $\gamma_2$ in $\spc{P}$ with common endpoints which contradicts the assumption on $\spc{P}$.

It remains to consider the case when $\alpha_1$ (and therefore $\alpha_2$) is not length minimizing.

Pass to a maximal length minimizing arc $\bar\alpha_1$ of $\alpha_1$.
Since $\Sigma$ is locally $\CAT(1)$, \ref{cor:discrete-paths} implies that 
there is another geodesic $\bar\alpha_2$ in $\Sigma$ that shares endpoints with $\bar\alpha_1$.
It remains to repeat the above construction for the pair $\bar\alpha_1$, $\bar\alpha_2$.

\parit{Comments.}
Thanks to \ref{ex:CAT-geodesic}, we can say now that \textit{a polyhedral space $\spc{P}$ is $\CAT(0)$ if and only if any two points in $\spc{P}$ can be joined by a unique geodesic}.

In the proof, one could also apply \ref{ex:geod-circle};
in this case the last part of the argument is not needed.

\parbf{\ref{ex:S3}.}
Note that it is sufficient to construct a polyhedral space $\spc{P}$ homeomorphic to the 3-disc such that (1) $\spc{P}$ is locally $\CAT(0)$ in its interior and (2) the boundary of $\spc{P}$ is locally concave; in particular, each edge on the boundary of $\spc{P}$ has angle at least~$\pi$.

Indeed, once $\spc{P}$ is constructed, taking the double of $\spc{P}$ along its boundary produces the needed metric on~$\mathbb{S}^3$.

The construction of $\spc{P}$ goes along the same lines as the construction of a Riemannian metric on the 3-disc with concave boundary and negative sectional curvature. 
This construction is given by Joel Hass in \cite{hass}.

\parit{Comments.} By the globalization theorem (\ref{thm:hadamard-cartan}) the obtained metric on $\mathbb{S}^3$ is not locally $\CAT(0)$.

This problem originated from a discussion 
in Oberwolfach
between
Brian Bowditch,
Tadeusz Januszkiewicz,
Dmitri Panov
and 
the third author.
Another solution was given by Karim Adiprasito \cite{adiprasito};
he proved that an example can be found among spaces that admit a cubulation into unit cubes.

\parbf{\ref{ex:baricenric-flag}.}
Checking the flag condition is straightforward once we know the following description of the barycentric subdivision.

Each vertex $v$ of the barycentric subdivision 
corresponds to a simplex $\triangle_v$ of the original triangulation.
A set of vertices forms a simplex in the subdivision 
if it can be ordered, say as $v_1,\ldots,v_\kay$,
so that the corresponding simplices form a nested sequence
\[\triangle_{v_1}\subset\ldots\subset\triangle_{v_\kay}.\]

\parit{Comment.}
There is a compact metrizable, contractible, and locally contractible space that does not admit a $\CAT(1)$ length metric \cite{petrunin-2024}.
According to this exercise, such space cannot be triangulated. 
%\footnote{ \color{red} I am not sure why this comment is here. is the point that such spaces can not be triangulated?}

\parbf{\ref{ex:flag>=pi/2}.}
Use induction on the dimension to prove that if in a spherical simplex $\triangle$ every edge is at least $\tfrac\pi2$, then
all dihedral angles of $\triangle$ are at least~$\tfrac\pi2$.

The rest of the proof goes along the same lines as the proof of the flag condition (\ref{thm:flag}).
The only difference is that a geodesic spends \textit{at least} $\pi$ on each visit to $\Star_v$.

\parit{Comment.}
Note that it is not sufficient to assume that all dihedral angles of the simplices are at least~$\tfrac\pi2$.
Indeed, the 2-dimensional sphere with the interior of a small rhombus removed is a spherical polyhedral space glued from four triangles with all the angles at least~$\tfrac\pi2$.
On the other hand the boundary of the rhombus is a closed local geodesic of length less than $2\cdot\pi$ in this space.
Therefore the space cannot be $\CAT(1)$.

\parbf{\ref{ex:polyhedron-glue}.}
Apply the globalization theorem (\ref{thm:hadamard-cartan}) with \ref{ex:cone+susp} and \ref{ex:flag>=pi/2}. 

\parbf{\ref{ex:tree}.}
The space $\spc{T}_n$ has a natural cone structure with the vertex formed by the completely degenerate tree --- all its edges have zero length.
Note that the space $\Sigma$
over which the cone is taken comes naturally with a triangulation 
with all-right spherical simplices.

The link of any simplex of this triangulation satisfies the no-triangle condition (\ref{def:flag}).
Indeed, fix a simplex $\triangle$ of the complex;
it can be described by the combinatorics of a possibly degenerate tree.
A triangle in the link of $\triangle$ can be described by three ways to resolve a degeneracy by adding one edge of positive length,
such that (1) any pair of these resolutions can be done simultaneously, but (2) all three cannot be done simultaneously.
Direct inspection shows that this is impossible.

Therefore, by Proposition~\ref{prop:no-trig} our complex is flag.
It remains to apply the flag condition (\ref{thm:flag}) and then \ref{ex:cone+susp}.

\parbf{\ref{ex:flag-aspherical}.}
If the complex $\mathcal{S}$ is flag, then its cubical analog $\square_{\mathcal{S}}$ is locally $\CAT(0)$ and therefore aspherical.

Assume now that the complex $\mathcal{S}$ is not flag. 
Extend it to a flag complex $\mathcal{T}$ by gluing a simplex in every clique (that is, a complete subgraph) of its one-skeleton.

Note that the cubical analog $\square_{\mathcal{S}}$ is a proper subcomplex in $\square_{\mathcal{T}}$.
Since $\mathcal{T}$ is flag,
$\tilde\square_{\mathcal{T}}$,
the universal cover of $\square_{\mathcal{T}}$, is $\CAT(0)$.
Let $\tilde\square_{\mathcal{S}}$ be the inverse image of $\square_{\mathcal{S}}$ in $\tilde\square_{\mathcal{T}}$.

Choose a cube $Q$ with minimal dimension in $\tilde\square_{\mathcal{T}}$ which is not present in $\tilde\square_{\mathcal{S}}$.
By Exercise~\ref{ex:locally-convex}, $Q$ is a convex set in $\tilde\square_{\mathcal{T}}$.
The nearest-point projection $\tilde\square_{\mathcal{T}}\to Q$ is a retraction.
It follows that the boundary $\partial Q$ is not contractible in $\tilde\square_{\mathcal{T}}\setminus\Int Q$.
Therefore the spheroid $\partial Q$ is not contractible in $\tilde\square_{\mathcal{S}}$.
That is, $\tilde\square_{\mathcal{S}}$ is not aspherical, and therefore $\square_{\mathcal{S}}$ is not aspherical either.

\parbf{\ref{ex:example-pi_infty-new}.}
The solution is very similar to the proof of Lemma~\ref{lem:example-pi_infty}, but a few changes are needed.

The cycle $\gamma$ is taken in the complement $\mathcal{S}\setminus\{v\}$ (or, alternatively, in the link of $v$ in $\mathcal{S}$).
Instead of a vertex, one has to take the edge $e$ in $\tilde Q$ that corresponds to $v$;
so we have to show the existence of a large cycle in $\tilde Q$ that is not contractible in $\tilde Q\setminus e$.
Consider $G$,
which is made from the squares parallel to the squares of the cubical complex and which meet the edges of the complex orthogonally at their midpoints.
Note that formally speaking, $G$ is not a subcomplex of the cubical analog.

\parbf{\ref{ex:funny-S}.}
In the proof we apply the following lemma. 
It follows from the disjoint discs property;
see \cite{edwards, daverman}.

\medskip

\parbf{Lemma.}
\textit{Let $\spc{S}$ be a finite simplicial complex that
is homeomorphic to an $m$-dimensional homology manifold for some $m\ge 5$.
Assume that all vertices of
$\spc{S}$ have simply connected links.
Then $\spc{S}$ is a topological manifold.}

\medskip

Note that it is sufficient to construct a simplicial complex $\spc{S}$
such that 
\begin{itemize}
\item $\spc{S}$ is a closed $(m-1)$-dimensional homology manifold;
\item $\pi_1(\spc{S}\setminus\{v\})\ne0$ for some vertex $v$ in $\spc{S}$;
\item $\spc{S}\sim \mathbb{S}^{m-1}$; that is, $\spc{S}$ is homotopy equivalent to~$\mathbb{S}^{m-1}$.
\end{itemize}

Indeed, assume such $\spc{S}$ is constructed.
Then the suspension
$\spc{R}\z=\Susp\spc{S}$
is an $m$-dimensional homology manifold with a natural triangulation coming from~$\spc{S}$.
According to the lemma,
$\spc{R}$ is a topological manifold.
According to the generalized Poincar\'{e} conjecture,
$\spc{R}\simeq\mathbb{S}^m$;
that is, $\spc{R}$ is homeomorphic to~$\mathbb{S}^m$.
Since $\Cone \spc{S}\simeq \spc{R}\setminus\{s\}$, where $s$ denotes the south pole of the suspension
and $\EE^m\simeq \mathbb{S}^m\setminus\{p\}$
for any point $p\in \mathbb{S}^m$,
we get 
\[\Cone \spc{S}\simeq\EE^m.\]

It remains to construct~$\spc{S}$.
Fix an $(m-2)$-dimensional homology sphere $\Sigma$ with a triangulation such that $\pi_1\Sigma\ne0$.
An example of that type exists for any $m\ge 5$; a proof is given in \cite{kervaire}.

Remove from $\Sigma$ the interior of one $(m-2)$-simplex.
Denote the resulting complex by~$\Sigma'$.
Since $m\ge 5$, we have $\pi_1\Sigma=\pi_1\Sigma'$.

Consider the product $\Sigma'\times [0,1]$. 
Attach to it the cone over its boundary $\partial (\Sigma'\times [0,1])$.
Denote by $\spc{S}$ the resulting simplicial complex
and by $v$ the tip of the attached cone.

Note that $\spc{S}$ is homotopy equivalent to the spherical suspension over $\Sigma$.
This suspension is a simply connected homology sphere and hence is homotopy equivalent to~$\mathbb{S}^{m-1}$.
 Hence $\spc{S}\sim\mathbb{S}^{m-1}$.

The complement $\spc{S}\setminus\{v\}$ is homotopy equivalent to~$\Sigma'$.
Therefore 
\[
\pi_1(\spc{S}\setminus\{v\})
=\pi_1\Sigma'
=\pi_1\Sigma\ne 0.
\]
That is, $\spc{S}$ satisfies the conditions above.

\parbf{\ref{ex:cube-infty=>cube-2}}; \ref{SHORT.cube-infty} $\Rightarrow$ \ref{SHORT.cube-2}.
Since any closed curve can be considered as a short map from the boundary of a disc with some metric, it can be extended to a short map from a disc.
Therefore any injective space is simply connected.

Therefore the globalization theorem and flag condition (\ref{thm:hadamard-cartan} and \ref{thm:flag}) imply that it is sufficient to show that each link in $Q$ is flag.
Further, by \ref{prop:no-trig} it is sufficient to show that the link of each cube in $Q$ satisfies the no-triangle condition.

Arguing by contradiction, we can assume that the no-triangle condition does not hold at a vertex $v$; that is, a 0-dimensional cube.
In this case $v$ is a vertex of three edges $e_x$, $e_y$, and $e_z$;
each pair of edges belongs to one of the squares $s_x$, $s_y$, and $s_z$ with complementary index, but the squares $s_x$, $s_y$, $s_z$ do not belong to one cube.
For higher dimensional cubes we have a product of this configuration with a cube.

\begin{wrapfigure}{r}{25 mm}
\vskip-0mm
\centering
\includegraphics{mppics/pic-13}
\end{wrapfigure}

Let $m_x$, $m_y$, and $m_z$ be the midpoints of $e_x$, $e_y$, and $e_z$ respectively.
Consider 3 balls with centers $m_x$, $m_y$, and $m_z$ and radius $\tfrac14$.
Observe that each pair of balls have a common point;
but all three together have no points of intersection.
However, in an injective space, the former implies the latter.
Indeed, we may think that $m_x$, $m_y$, and $m_z$ are images of a distance-preserving map from a 3-point subset in a 4-point metric space with the fourth point at distance $\tfrac14$ from the rest,
and this map cannot be extended to a short map on the whole 4-point space.
(For more on injective spaces, see \cite{petrunin-2023}.)
That is, $(Q,\ell^\infty)$ is not an injective space --- a contradiction.

\parit{\ref{SHORT.cube-1} $\Rightarrow$ \ref{SHORT.cube-2}.}
Observe that median point $m(x,y,z)$ depends continuously on the triple of points $(x,y,z)$ and $m(x,x,y)=x$.

Given a loop $\gamma\:[0,1]\to Q$ with base at $p=\gamma(0)=\gamma(1)$,
consider the map $(a,b)\mapsto m(p,\gamma(a),\gamma(b))$ of the triangle $\triangle$ defined by $0\z\le a \z\le b\z\le 1$.
Note that the boundary of triangle runs along $\gamma$.
It follows that $\gamma$ is null homotopic and therefore $Q$ is simply connected.

\begin{wrapfigure}{r}{25 mm}
\vskip-3mm
\centering
\includegraphics{mppics/pic-15}
\end{wrapfigure}

It remains to check that all links of $Q$ satisfy the no-triangle condition.

Assume that a link of $Q$ does not satisfy the no-triangle condition.
The same way as in the previous problem, we can assume that it is a link of a vertex;
so we have a configuration of three squares $s_x$, $s_y$, and $s_z$, 
three edges $e_x$, $e_y$, and $e_z$, and one common vertex $v$ as above.
Let $x$, $y$, and $z$ be the squares $s_x$, $s_y$, and $s_z$, respectively.
Observe that the geodesics $[xy]_{\ell^1}$, $[xz]_{\ell^1}$, and $[yz]_{\ell^1}$ are uniquely defined and they have no common point.
It follows that the triple $(x,y,z)$ does not have a median; 
that is, $(Q,\ell^1)$ is not a median space --- a contradiction.

\parbf{\ref{ex:chopping-triangle}.} 
Observe that the triangle $\trig pqx$ is degenerate, in particular it is thin.
It remains to apply the inheritance lemma (\ref{lem:inherit-angle}).

\parbf{\ref{ex:concave-triangle}.}
By approximation,
it is sufficient to consider the case when $S$ has polygonal sides.

The latter case can be done by induction on the number of sides.
The base case of triangle is evident.

To prove the induction step, apply Alexandrov's lemma (\ref{lem:alex}) 
together with the construction in the inheritance lemma (\ref{lem:inherit-angle}).

\parbf{\ref{ex:bishop-sphere}.}
From Exercise~\ref{ex:hemisphere}, it follows that if the image in $\mathbb{S}^2$ of a lune in $K$ has perimeter smaller than $2\cdot\pi$, then it contains a closed hemisphere in its interior
or is contained in an open hemisphere.
The same holds for a triangular region with concave sides.

Since the interior of $K$ does not contain a closed hemisphere, the first case cannot occur.
It remains to apply the argument in \ref{thm:bishop-plane}.

We still have to check that the induced length metric on $K$ is finite.

This can be done by applying Smirnov’s theorem \cite[3.12]{duren} to a conformal parametrization of $K$.

\parbf{\ref{ex:two-planes}.}
The space $\tilde K$ is a cone over the branched covering $\Sigma$ of $\mathbb{S}^3$ infinitely branching along two great circles.

If the planes are not orthogonal, then the minimal distance between the circles is less than~$\tfrac\pi2$.
Assume that this distance is realized by a geodesic $[\xi\zeta]$.
The polygonal line made by four liftings of $[\xi\zeta]$ forms a closed
local geodesic in~$\Sigma$.
By Proposition~\ref{cor:loc-geod-are-min}
(or Corollary~\ref{cor:closed-geod-cat}),
$\Sigma$ is not $\CAT(1)$.
Therefore, by \ref{ex:cone+susp}, $\tilde K$ is not $\CAT(0)$.

If the planes are orthogonal, then the corresponding great circles in $\mathbb{S}^3$ are subcomplexes of a flag triangulation of $\mathbb{S}^3$ with all-right simplices.
The branched cover is also flag.
It remains to apply the flag condition~\ref{thm:flag}.

\parit{Comments.}
In \cite{charney-davis-1993}, Ruth Charney and Michael Davis
gave a complete answer to the analogous question for three planes.
In particular, they show that if a covering space of $\EE^4$
branching at three planes intersecting at the origin is $\CAT(0)$, then these all are complex planes for some complex structure on~$\EE^4$.

\parbf{\ref{ex:hemisphere}.}
Let $\alpha$ be a closed curve in $\mathbb{S}^2$ of length $2\cdot\ell$.

\begin{wrapfigure}{r}{35 mm}
\vskip-4mm
\centering
\includegraphics{mppics/pic-722}
\end{wrapfigure}

Assume $\ell<\pi$.
Let $\check\alpha$ be a subarc of $\alpha$ of length $\ell$, with endpoints $p$ and~$q$. 
Since $\dist{p}{q}{}\le\ell<\pi$, there is a unique geodesic $[pq]$ in~$\mathbb{S}^2$.
Let $z$ be the midpoint of $[pq]$. 

We claim that $\alpha$ lies in the open hemisphere $H$ centered at~$z$.

Assume the contrary; that is, $\alpha$ meets the equator $\partial H$ at a point $r$.
Without loss of generality we may assume that $r\in\check\alpha$.

The arc $\check\alpha$ together with its reflection in $z$ form a closed curve of length $2\cdot \ell$ which meets $r$ and its antipodal point~$r'$.
Thus $\ell\z=\length \check\alpha\ge \dist{r}{r'}{}=\pi$ --- a contradiction.

\parit{Solution via the Crofton formula.}
Let $\alpha$ be a closed curve in $\mathbb{S}^2$ of length $< 2\cdot\pi$.
We wish to prove that $\alpha$ is contained in a hemisphere in~$\mathbb{S}^2$.
By approximation it suffices to prove this for smooth curves $\alpha$ of length $< 2\cdot\pi$ with transverse self-intersections.

Given $v\in \mathbb{S}^2$, denote by $v^\perp$ the equator in $\mathbb{S}^2$ with the pole at~$v$.
Further, $\# X$ will denote the number of points in the set~$X$.

Obviously, if $\#(\alpha\cap v^\perp) =0$, then $\alpha$ is contained in one of the hemispheres determined by~$v^\perp$.
Note that $\#(\alpha\cap v^\perp)$ is even for almost all~$v$ (specifically, for every $v$ such that $\alpha$ is transverse to $v^\perp$).

Therefore, if $\alpha$ does not lie in a hemisphere, then
$\#(\alpha\cap v^\perp) \ge 2$ for almost all $v\in\mathbb{S}^2$.

By the Crofton formula we have that
\begin{align*}
\length(\alpha)
&=\frac 1 4\cdot \int\limits_{v\in \mathbb{S}^2}\#(\alpha\cap v^\perp)\ge 2\cdot\pi.
\end{align*}

\parbf{\ref{ex:inner-support}.}
Since $\Omega$ is not two-convex,
we can fix a simple closed curve $\gamma$ that lies in the intersection of a plane $W_0$ and $\Omega$, 
and is contractible in $\Omega$ but not contractible in $\Omega\cap W_0$.

Let $\phi\:\DD\to \Omega$ be a disc that shrinks~$\gamma$.
Applying the loop theorem (arguing as in the proof of Proposition~\ref{prop:3d-strong-2-convexity}), we can assume that $\phi$ is an embedding and $\phi(\DD)$ lies on one side of~$W_0$.

Let $Q$ be the bounded closed domain cut from $\EE^3$ by $\phi(\DD)$ and~$W_0$. 
By assumption it contains a point that is not in~$\Omega$. 
Changing $W_0,\gamma$ and $\phi$ slightly, we can assume that such a point lies in the interior of~$Q$.

Fix a circle $\Gamma$ in $W_0$ that surrounds $Q\cap W_0$.
Since $Q$ lies in a half-space with boundary $W_0$, there is a
smallest spherical dome with boundary $\Gamma$ that includes the set $R=Q\setminus\Omega$.

\begin{wrapfigure}{r}{43mm}
\vskip-0mm
\centering
\includegraphics{mppics/pic-19}
\end{wrapfigure}

The dome has to touch $R$ at some point~$p$.
The plane $W$ tangent to the dome at $p$ has the required property --- the point $p$ is an isolated point of the complement $W\setminus \Omega$.
Further, by construction a small circle around $p$ in $W$ is contractible in $\Omega$.

\parbf{\ref{ex:convex+saddle+broken=>PL}.}
The proof is simple and visual, but it is hard to write it down formally in a non-tedious way.
Read the following, wait until you feel the geometric picture, and then read it again.
If it still does not work, try reading the original \cite[§ 2]{shefel-1964}.

\parit{\ref{SHORT.ex:convex+saddle+broken=>PL:a}.}
By \ref{clm:F'} on page \pageref{clm:F'}, $\partial S$ is a finite union of line segments and points.

Choose a point $p\in S$.
Since $S$ is locally concave, there is a plane $\Pi$ that locally supports $S$ at $p$.
Choose a closed $\eps$-neighborhood $B\ni p$ such that $B\setminus \Omega'$ is convex.
Show that (1)
$V=\Pi\cap B\cap S$ is convex, (2)
$\Pi$ locally supports $S$ at any point of $V$,
and (3)
\[W=\Omega'\cap B\cap \Pi=(B\cap\Pi)\setminus V\]
has convex connected components.
(In particular, $p$ is not an isolated point in $S\cap \Pi$.)

It remains to apply this observation with an open-closed argument.

\parit{\ref{SHORT.ex:convex+saddle+broken=>PL:b}.}
Choose a point $p\in S$.
Given a local supporting plane $\Pi$ at $p$, denote by $A_\Pi$ the set provided by part \ref{SHORT.ex:convex+saddle+broken=>PL:a};
so, $A_\Pi$ is a line segment or a polygon lying in $\bar S$.

Assume that $p$ is not an interior point of $A_\Pi$ for some (and therefore any) supporting plane $\Pi$ at $p$.
Note that in this case $p$ lies on a line segment in $\bar S$ with the endpoints on $\partial S$;
denote this segment by $[xy]$.
In particular, all locally supporting planes at $p$ contain $[xy]$.
Consider two outermost supporting planes at $p$, say $\Pi_1$ and $\Pi_2$ (if the supporting plane is unique, then $\Pi_1=\Pi_2$).

Let us show that a neighborhood of $p$ in $S$ is covered by $A_{\Pi_1}\cup A_{\Pi_2}$.
Observe that this statement leads to a solution.

Assume that this is not the case.
Apply \ref{SHORT.ex:convex+saddle+broken=>PL:a} to show that there is a sequence of points $p_n\in S$ and line segments $[x_n\,y_n]\ni p_n$ with the same properties as above such that $p_n\to p$ and lines containing $[x_ny_n]$ approach the line of $[xy]$ as $n\to\infty$.
Moreover for one plane, say $\Pi_1$, we have $b_n>0$ for all $n$, and $b_n\z=o(a_n)$, where $a_n$ and $b_n$ denote the distance from $p_n$ to $[xy]$ and $\Pi_1$, respectively.

Recall that $\partial S$ is a union of finitely many line segments and points.
Therefore, after passing to a subsequence, we can choose two line segments $X$ and $Y$ in $\partial S$ (possibly degenerate) such that $x_n\in X$ and $y_n\in Y$ for any $n$.
Show that $X\subset \Pi_1$ or $Y\subset \Pi_1$.
Conclude that $[x_n\,y_n]\subset \Pi_1$ and arrive at a contradiction.

\parbf{\ref{ex:CAT=>two-convex}.}
Suppose $K$ is not two-convex;
let $\gamma$ be a closed simple curve in a plane $W$ that violates Definition~\ref{def:two-convex}.
Note that $W$ is not a vertical plane;
denote by $V$ the 3-dimensional subspace that is spanned by $W$ and the vertical direction.

Note that $\gamma$ is contractible in $V\cap K$.
Argue as in \ref{ex:inner-support} to show that there exists a plane triangle $\triangle\subset V$ whose sides lie completely in $K$,
but whose interior contains points from the complement $\EE^m\setminus K$.
This triangle is not thin in $K$.

Now let us show that the converse does not hold in $\EE^4$.
Let $K$ be the subgraph of the function
\[f(x,y,z)=10^{10}\cdot \min\{\,|y|+x,|z|-x\,\}.\]

Note that $K$ is two-convex, since it is the intersection of two two-convex sets (the subgraphs of $10^{10}(|y|+x)$ and $10^{10}(|z|-x)$).
Furthermore, $K$ is the complement of two open dihedral angles $D_1$ and $D_2$.
Therefore, $K=\Cone\Sigma$, where $\Sigma$
is the complement of two open lenses, say $L_1$ and $L_2$, that correspond to $D_1$ and $D_2$.
(Each lens is an intersection of two open hemispheres.)

Show that $\partial L_1\cap\partial L_2$ is a closed local geodesic of length smaller than $2\cdot\pi$ in $\Sigma$.
By \ref{cor:closed-geod-cat}, $\Sigma$ is not $\CAT(1)$.
Therefore, by \ref{ex:cone+susp}, $K$ is not $\CAT(0)$.

\parbf{\ref{ex:finite-action-CAT}}; \ref{SHORT.ex:finite-action-CAT:finite}
Suppose $G$ is finite.
Choose its orbit $\{p_1,\ldots,p_n\}$.
Consider the barycenter $z$ of the array $\bm{p}=(p_1,\ldots,p_n)$ with equal masses;
in other words, $z\df\spx{\bm{p}}(\tfrac1n,\ldots,\tfrac1n)$.
Observe that $z$ is a fixed point of the action.

\parit{\ref{SHORT.ex:finite-action-CAT:compact}.}
Let $\mu$ be the probability Haar measure on $G$.
Choose a point $p\in \spc{U}$ and consider the function
\[f=\tfrac12\cdot\int_{g\in G} \distfun_{g\cdot p}^2\cdot\mu.\]
Show that $f$ is 1-convex.
By \ref{lem:argmin(convex)}, $f$ has a unique minimum point, say~$z$.
Observe that $f$ is $G$-invariant;
therefore $z$ is a fixed point.

\parbf{\ref{ex:bary-jensen}}; \ref{SHORT.ex:bary-jensen:moment}.
Let $\gamma$ be geodesic from $z$ to $q$ and let $r=\dist{q}{z}{}$.
Consider the function
\[\phi(t)=\tfrac12\cdot \sum_i \mu_i\cdot \dist[2]{\gamma(t)}{p_i}{}.\]
Observe that $\phi$ is 1-convex, and $\phi^+(0)\ge 0$.
Conclude that
$\phi(r)\ge \tfrac12\cdot r^2$
and observe that it implies the statement.

\parit{\ref{SHORT.ex:bary-jensen:moment+}.}
Given two tangent vectors $v,w$, define their \emph{scalar product} as
\[\langle v,w\rangle\df|v|\cdot |w|\cdot\cos\theta,\]
where $\theta$ is the angle between $v$ and $w$;
if one of the vectors is zero, then you can choose $\theta$ arbitrary, say $\theta=0$.

Let $v_i$ be the tangent vector at $z$ of the geodesic path from $z$ to $p_i$.
Note that $\dist{z}{p_i}{}=|v_i|$ for each $i$.

Show that
\[\mu_0\cdot\langle w,v_0\rangle+\ldots+\mu_k\cdot\langle w,v_k\rangle\le 0 \leqno{(*)}\]
for any $w\in \T_z$.

Use comparison to show that
\[\dist[2]{p_i}{p_j}{}\ge \dist[2]{z}{p_i}{}+\dist[2]{z}{p_j}{}-2\cdot \langle v_i,v_j\rangle.\]
It remains to apply your algebraic skills.

\parit{\ref{SHORT.ex:bary-jensen:jensen}.}
Let $v_i\in\T_z$ be as above.

Define the differential $\dd_z h\:\T_z\to\RR$ of $h$ at $z$.
Observe that
\[h(p_i)\ge h(z)+\dd_z h(v_i)\]
for each $i$.
Further, show and use that there is a vector $v\in \T_z$ such that
\[\dd_z h(w)+\langle v,w\rangle\ge 0
\leqno{({*}{*})}\]
for any $w\in \T_z$.
In fact, one can assume that $v$ is \emph{gradient} of $-h$ at $z$ (briefly, $v=\nabla_z(-h)$),
which is uniquely defined as a vector that meets $({*}{*})$ and such that $\dd_zh(v)+|v|^2=0$.

It remains to do some algebra based on $({*})$ and $({*}{*})$.

\parit{Comment.}
The differential and gradient are defined in our book~\cite{alexander-kapovitch-petrunin-2025}, but the definitions and needed properties are quite straightforward.
Two more proofs of Jensen's inequality, which is part \ref{SHORT.ex:bary-jensen:jensen}, were presented by Karl-Theodor Sturm \cite{sturm}.

\parbf{\ref{ex:barysimple}.}
The first part follows directly from the definitions.
For the second part consider the cone over a circle with length bigger than $2\cdot\pi$,
or the product of a tripod with the real line.

\parbf{\ref{lem:nondeg-test-with-balls}}; \ref{SHORT.lem:nondeg-test-with-ball:U}.
Choose $q\in \bigcap_i B_i$.
Assume $s\notin \bigcup_i B_i$.
Observe that $\dist{p_i}{s}{}>\dist{p_i}{q}{}$ for any $i$ and apply \ref{thm:up-convex:bry}.
Conclude that $s\notin\spx{\bm{p}}(\triangle^k)$.

\parit{\ref{SHORT.lem:nondeg-test-with-balls:nondeg}.}
Show that there is a point, say $z$, that minimizes 
\[\ell=\max_i\{\dist{p_i}{z}{}-r_i\}.\]
Note that $\ell>0$.

Let us show that $\ell=\dist{p_i}{z}{}-r_i$ for any $i$.
Assume the contrary; that is, that $\ell>\dist{p_j}{z}{}-r_j$ for some $j$.
Choose a point $q_j\in\bigcap_{i\ne j} B_i$.
Note that $\ell>\dist{p_j}{x}{}-r_j$
for any point $x\in \left]zq_j\right[$.
Therefore $\ell$ is not minimal --- a contradiction.

Observe that $z\in \spx{\bm{p}}(\triangle^k)$.
By \ref{SHORT.lem:nondeg-test-with-ball:U}, $\spx{\bm{p}}(\partial\triangle^k)\subset \bigcup_{i} B_i$.
Hence $\spx{\bm{p}}(\triangle^k)\setminus\spx{\bm{p}}(\partial\triangle^k)\ne \emptyset$.

\parit{\ref{SHORT.lem:nondeg-test-with-balls:nondeg+}.}
Choose $z\in \spx{\bm{p}}(\triangle^k)\setminus\spx{\bm{p}}(\partial\triangle^k)$.
Show that for each $j$, there is a point $z_j$ such that
\[\dist{p_i}{z_j}{}<\dist{p_i}{z}{}\]
for any $i\ne j$.
Choose $r_i=\min_{j\ne i}\{\dist{p_i}{z_j}{}\}$.

\parbf{\ref{ex:helly}.}
We can assume that $n=m+2$ and $K_1\cap\ldots\cap K_{n}\z=\emptyset$.
Choose a point array $\bm{p}=(p_1,\dots,p_n)$ such that
\[p_i\in \bigcap_{j\ne i} K_j\]
for each $i$.
Observe that $\spx{\bm{p}}(\partial\triangle^k)\subset \bigcup_iK_i$.

By \ref{thm:bihoelder}, the simplex $\spx{\bm{p}}$ is degenerate.
Therefore, $\spx{\bm{p}}(\triangle^k)\subset \bigcup_iK_i$.
Apply Sperner's lemma \cite{sperner-lemma} to show that $\spx{\bm{p}}(\triangle^k)(\bm{\mu})\z\in \bigcap_iK_i$ for some $\bm{\mu}\in \triangle^k$ and each $i$.
Arrive at a contradiction.

\parbf{\ref{ex:dim-top-haus-CAT}.}
Let $T$ be a binary rooted tree.
Choose a metric on $T$ such that each edge from level $n$ to $n+1$ has length $\lambda_n$ for some $0<\lambda<1$.
Passing to completion of the obtained space, we add to it a \textit{crown} which is homeomorphic to the Cantor set.
Observe that the completion has topological dimension 1.

Calculate the Hausdorff dimension of the crown for given $\lambda$ and make a conclusion.
To do the last step, pass to the product of the obtained example with $\RR^{m-1}$.

{
\newpage
\phantomsection
\sloppy

\input{invitation-CAT.ind}

\def\emph{\textit}

\printbibliography[heading=bibintoc]
\fussy

}

\end{document}